\documentclass[11pt]{article}
\usepackage{amssymb,amsmath,fullpage}
\numberwithin{equation}{section}
\def\vP{\Phi}
\def\vp{\varphi}\def\vt{\vartheta}
\def\a{\alpha}\def\ta{\tilde\alpha}
\def\A{\cal A}
\def\b{\beta}
\def\c{\gamma}
\def\d{\delta}

\def\s{\sigma}
\def\t{\tau}
\def\x{\xi}
\def\y{\eta}
\def\z{\zeta}
\def\sgn{\hbox{\rm sgn}}

\def\diag{\hbox{\rm diag}}
\def\Res{\mathop{\hbox{\rm Res}}}
\def\8{\infty}
\def\<{\langle}
\def\>{\rangle}
\def\raisebullet{{\raise 0.5mm\hbox{$\sssize\bullet$}}}
\def\wtil{\widetilde}
\def\wt{\widetilde}
\def\what{\widehat}
\def\X21{{\textstyle\frac{x_2}{x_1}}}

\def\a{\alpha}
\def\b{\beta}

\def\s{\sigma}
\def\t{\tau}
\def\z{\zeta}
\def\8{\infty}
\def\sgn{\mbox{{\rm sgn}}}
\def\A{{\cal A}}

\def\qed{$\square$}
\def\Vec#1{\mbox{\boldmath $#1$}}

\newtheorem {thm}{Theorem}[section]
\newtheorem {lem}[thm]{Lemma}
\newtheorem {cor}[thm]{Corollary}
\newtheorem {prop}[thm]{Proposition}
\newtheorem {df}[thm]{Definition}
\newtheorem {rmk}[thm]{Remark}

\begin{document}
\title{
Gauss decomposition and 
$q$-difference equations \\\
for Jackson integrals of symmetric Selberg type
}
\author{Masahiko I{\small TO}
\\
\footnotesize  
Graduate School of Mathematics, Nagoya University, 
Chikusa-ku, Nagoya 464-8602, Japan 
\footnote{Current address: Department of mathematical sciences, University of the Ryukyus, Okinawa 903-0213, Japan
email: {\tt mito@sci.u-ryukyu.ac.jp}
}
\\
}
\date{
}
\maketitle
\begin{abstract} 
We provide explicit expressions for two types of first order $q$-difference systems for the Jackson integral of symmetric Selberg type. One is the $q$-difference system known to be the $q$-KZ equation and 
the other is the $q$-difference system for parameters different from the $q$-KZ equation. 
We use a basis of the systems introduced by Matsuo in his study of the $q$-KZ equation. 
As a result, the similarity of these two systems is discussed by concrete calculations. 
Intermediate calculations are made use of the {\em Riemann-Hilbelt method for $q$-difference equation from connection matrix} established by Aomoto.

\end{abstract}



\section{Introduction}

Let $q=e^{2\pi\sqrt{-1}\tau}$, Im $\!\tau>0$, be the elliptic modulus.
Let $\vP(t)$ be a $q$-multiplicative function on the algebraic torus $(\mathbb{C}^*)^n$ 
defined by
\begin{equation}
\label{eq:0.1}
\vP(t)=
\vP_{n,m}(t):=
t_1^{\a_1}\cdots t_n^{\a_n}
\prod_{j=1}^n 
\prod_{k=1}^m 
             \frac{(t_j/x_k)_\8}{(t_j q^{\b_k}/x_k)_\8}
       \prod_{1\le i<j\le n}
             \frac{(q^{1-\c}t_j/t_i)_\8}{(q^\c t_j/t_i)_\8},
\end{equation}
where $\a_j=\ta -n+j -2(j-1)\c$.
In the papers \cite{ao2,ak1}, 
Aomoto and Aomoto--Kato introduced 
the notion of the {\it symmetric part} 
$H^n((\Bbb C^*)^n,\vP,\nabla)_{sym}$ 
of $q$-analog of the twisted de Rham cohomology attached to $\vP(t)$,
whose dimension is ${n+m-1\choose m-1}$ 
if parameters $\ta$, $\b_k$, $x_k$ and $\c$ are generic.
If we take a basis $\{\psi_l(t)\}_{l\in L}$ of the cohomology $H^n((\Bbb C^*)^n,\vP,\nabla)_{sym}$,  
we can construct a solution of a system of holonomic $q$-difference equations on $(\Bbb C^*)^m$ 
by using the Jackson integrals $\wtil\psi_l=\int \vP \psi_l \varpi$ 
(see \eqref{eq:1.1.2} for the definition of Jackson integral). 
When we denote the $q$-shift $x_j \to qx_j$ for a value $x_j$ by $T_j$ 
the $q$-difference system is expressed by a suitable matrix $K_j(x)$ 
of rank ${n+m-1\choose m-1}$ given by 
$T_j(\wtil\psi_l)_{l\in L}=(\wtil\psi_l)_{l\in L}K_j(x)$.
\par 
In the papers \cite{ma1,ma2}, 
Matsuo claimed
that by taking a suitable basis 
the Jackson integrals for it 
give a solution of the quantized Knizhnik--Zamolodchikov difference equations
for the matrix coefficients of the product of 
intertwining operators $R_{ij}$ called $R$-matrices for the quantum affine group $U_q(\widehat{sl}_2)$
in the sense of Frenkel and Reshetikhin \cite{fr}, 
i.e., the matrix $K_j(x)$ is expressed as  
$$\textstyle
K_j(x)=R_{j,j+1}(\frac{x_j}{x_{j+1}})\cdots R_{j,m}(\frac{x_j}{x_{m}})
  \ D_j
  \ R_{j,1}(\frac{qx_j}{x_1})\cdots R_{j,j-1}(\frac{qx_j}{x_{j-1}})
\qquad (1\le j\le m),
$$
where $D_j$ is some diagonal matrix. 
Varchenko \cite{v} 
made Matsuo's work complete. 
According to the result of \cite{ma2,v}, 
the system of $q$-difference equations 
for the Jackson integrals $\wtil\psi_{l}(t)$ with respect to the $q$-shift $x_j \to qx_j$ 
eventually reduces to that of the cases $m=2$ 
because $R$-matrix $R_{ij}$ is decomposed into a direct sum of these of the case $m=2$.
Therefore we consider the problem 
of finding an explicit expression of the matrix $K_j$ or $R_{ij}$ for $m=2$. 
\par
On the other hand, 
this problem has already been studied by Mimachi \cite{mi1,mi2}. 
He introduced one of expressions of the matrix $K_j(x)$ 
by using values of certain Schur polynomials 
and evaluated it explicitly when
$n=1,2$ and $3$ (see \cite{mi2}). 
Aomoto and Kato \cite{ak6} 
also gave another approach to express it 
in terms of the Gauss matrix decomposition \cite{g} 
and evaluated it explicitly 
in the case where $n=1$ and $2$, 
but in the form $K_j(0)K_j(x)^{-1}$ .
They used a method which they call the Riemann--Hilbert problem for $q$-difference equation
from connection matrix \cite{ao4, ao5}. 
It is a surprising method 
because, under some assumptions, the matrix which represents $q$-difference system
is exactly determined only from the information of connection matrix 
between asymptotic behaviors of its fundamental solutions.
However, for evaluating the matrix $K_j(x)$, 
we do not need such method if we can evaluate the $R$-matrix $R_{ij}$. 
Actually, the explicit form of the $R$-matrix as the Gauss decomposition 
is not so difficult to find out.
This is one of the aims of the present paper and will lead us to Theorems \ref{thm:1.3.3} and \ref{thm:1.3.4} in Section \ref{section:1}. 
(However we will confirm in Subsection \ref{subsection:4.2}
that the explicit form of the matrix $K_j(x)$ is also obtained from Aomoto's Riemann--Hilbert method.)
\par
When $n=1$ the Jackson integral \cite{an} associated with $\vP_{n,2}(t)$ is 
equivalent to 
Heine's hypergeometric series 
\begin{align*}
{}_2\vp_1(q^\a,q^\b;q^\c;x)
&=
\sum_{\nu=0}^\8
\frac{(q^\a)_\nu(q^\b)_\nu}{(q^\c)_\nu(q)_\nu}x^\nu
=\frac{1}{1-q}
    \frac{(q^{\a})_\8(q^\b)_\8}{(q^{\c})_\8(q)_\8}
     \int_0^1 t^X\frac{(q^\c t)_\8(qt)_\8}{(q^{\a}t)_\8(q^{\b}t)_\8}{d_q t\over t},
\end{align*}
%
%
where $x=q^X$.
Heine's hypergeometric series 
satisfies 
the following transformation formula \cite[p.13, Eq.(1.4.1)]{gr} 
\begin{equation}
\label{eq:0.2}
{}_2\vp_1(q^\a,q^\b;q^\c;x)=\frac{(q^\a)_\8(q^\b x)_\8}{(q^\c)_\8(x)_\8}{}_2\vp_1(x,q^{\c-\a};q^\b x;q^\a).
\end{equation}
One of the reason 
the transformation \eqref{eq:0.2} holds is the equality
$$
(q^{B})_\8\int_0^1 t^A\frac{(tq^{X+C})_\8(qt)_\8}{(tq^X)_\8(tq^{B})_\8}{d_q t\over t}
=(q^{C})_\8\int_0^1 t^X\frac{(tq^{A+B})_\8(qt)_\8}{(tq^A)_\8(tq^{C})_\8}{d_q t\over t},
$$
which changes $q^A$ and $q^X$. 
It seems interesting to treat the $q$-difference system 
with respect to the parameter $q^{\ta}$ 
as that of $x$. 
We state the $q$-difference system 
for the same basis as Matsuo's case 
with respect to the parameter shift $T_{\ta}:\ta\to \ta +1$ for $m=2$:
\begin{equation}
\label{eq:0.3}
T_{\ta}
(\wtil\psi_l)_{l\in L}=(\wtil\psi_l)_{l\in L}A(q^{\ta}),
\end{equation}
where $A(q^{\ta})$ is a suitable rational matrix of rank $n+1$. 
As we see in Theorem \ref{thm:5.0.2} in Section \ref{section:5}, although we do not know each element of 
the matrix $A(q^{\ta})$, we have its Gauss decomposition. 
The Gauss decomposition of the matrix $A(q^{\ta})$ 
is very similar to that of the matrix $K_j(x)$
(see Theorem \ref{thm:5.0.2} in Section \ref{section:5} and compare it with Theorem \ref{thm:1.3.3}). 
Furthermore, the $R$-matrix $R_{ij}$ is determined 
from asymptotic behaviors of certain special solutions for the $q$-difference system \eqref{eq:0.3}
through \eqref{eq:3.1}, \eqref{eq:3.3}, \eqref{eq:3.4}, \eqref{eq:9.15} and \eqref{eq:9.16}. 
In particular, the upper and lower triangular matrices 
of the Gauss matrix decomposition of $R_{ij}$
are determined from the matrices $A(0)$ or $A(\8)$ 
via \eqref{eq:2.2}, \eqref{eq:2.3}, \eqref{eq:3.6}, \eqref{eq:3.7} and \eqref{eq:9.11}  
(see Remark \ref{rmk:9.4} 
in Appendix \ref{section:9} and Examples in Appendix \ref{section:10}). 
\par
In Section \ref{section:2}, we review the Riemann--Hilbert problem 
for $q$-difference equation from connection matrix. 
In order to evaluate $A(q^{\ta})$, 
we use Aomoto's Riemann--Hilbert method for it 
because $A(q^{\ta})$ is determined 
only from the data of the principal connection matrix $G$ 
which has been studied in \cite{ao3,ak2-1,ak2-2,ak3,ak4,ak5,akm}. 
The principal connection matrix $G$ can also be deduced from 
the hypergeometric pairing studied by Tarasov and Varchenko \cite{tv}, 
which is related to this problem.\\
\par
\noindent
{\bf Remark}
This note was written in 1997 when the author was a graduate student at Nagoya University. At that time, the Internet was not yet popular enough, and there were page restrictions on paper publication. Although he compiled it in notebook form, the intermediate calculations were too long, so he did not publish it and only gave printed copies to a limited number of people involved. 
Therefore, there were several papers \cite{ao5,ak6} in the bibliographic list at that time that had the title of this note. 
On the other hand, time has passed, and recently the author has found 
that the main results (Theorems \ref{thm:1.3.3} and \ref{thm:5.0.2}) of this note can be derived relatively simply by a method different from the method in this note, so he has published another proof of these theorems in \cite{ito}. 
Furthermore, he heard from Prof.~Y.~Yamada that there was an application of these theorems, 
and the results of \cite{ito} were cited in the recent paper \cite{AHKOSSY}, 
so he decided to publish this note here. 
Lastly it should also be noted 
that in \cite{AHKOSSY} the similarity 
between the matrices $A(q^{\ta})$ and $K_j(x)$ is explained as consequence of the so-called {\it base-fiber duality} of the gauge theory. 

\tableofcontents


\section{Quantized Knizhnik--Zamolodchikov difference equations 
and $R$-matrix.}
\label{section:1}
\subsection{Notation}
Let $q=e^{2\pi\sqrt{-1}\tau}$, ${\rm Im}\,\tau>0$, be the elliptic modulus.
For an arbitrary number $c\in \Bbb C$, $0<|c|<1$, we use the notations 
$$(x;c)_{\8}:=\prod_{\nu=0}^\infty (1-c^\nu x), 
\quad (x;c)_\nu :=\frac{(x;c)_{\8}}{(xc^\nu ;c)_{\8}},
\quad {}_r (x;c)_s=\frac{(x;c)_r}{(x;c)_{r-s}(x;c)_s}.
$$
If $c=q$, we simply write $(x)_\8:=(x;q)_\8$ and $(x)_\nu:=(x;q)_\nu$. 
Let $\vt(x)$ be  the {\it Jacobi elliptic theta function} defined by 
$$
\vt(x):=(x)_\8(q/x)_\8(q)_\8.
$$
We also use the notations 
$$\vt(x)_r:=\vt(x)\cdot\vt(xq^\c)\cdots \vt(xq^{(r-1)\c}) 
\quad\hbox{ and }\quad 
{}_r\vt(x)_s:=\frac{\vt(x)_r}{\vt(x)_{r-s}\cdot \vt(x)_{s}},
$$
which have the relations 
$$
\lim_{q\to 0}\vt(x)_r=(x;q^\c)_r
\quad\hbox{ and }\quad 
\lim_{q\to 0}{}_r\vt(x)_s={}_r(x;q^\c)_s, 
$$
if we fix $q^\c$ as a single character $c$ such as a number that does not depend on $q$. 
Let $\vP(t)=\vP_{n,m}(t)$ be the same function as \eqref{eq:0.1}. 
The function $\vP(t)$ satisfies a quasi-symmetric property with respect to the symmetric group 
$\frak S_n$ of $n$th order such that 
\begin{equation}
\label{eq:1.1.1}
\s \vP(t)=U_\s(t) \ \vP(t),\quad \s\in \frak S_n,
\end{equation}
with a $q$-periodic function $U_{\s}(t)$ as  
$$
U_\s(t):=
\!\!\!\!\!
\prod_{\substack{i<j \\ \s^{-1}(i)>\s^{-1}(j)}}
\!\!\!\!\!
\Big({t_j\over t_i}\Big)^{\!2\c -1}
\frac{\vt(q^{\c}t_j/t_i)}{\vt(q^{1-\c}t_j/t_i)}
=\sgn\,\s
\!\!\!\!\!\!
\prod_{\substack{i<j \\ \s^{-1}(i)>\s^{-1}(j)}}
\!\!\!\!\!
q^{-\c}
\Big({t_j\over t_i}\Big)^{\!2\c}
\frac{\vt(q^{\c}t_j/t_i)}{\vt(q^{-\c}t_j/t_i)}
$$
where $\{U_\s(t)\}_{\s\in \frak S_n}$ satisfies the {\it one cocycle} condition 
$
U_{\s\s'}(t)=U_\s(t)\cdot \s U_{\s'}(t).
$ 
\begin{df}{\rm
For an arbitrary point $\x=(\x_1,\ldots,\x_n)\in (\Bbb C^*)^n$ 
and a function $f(t)$ of $t=(t_1,\ldots,t_n)\in (\Bbb C^*)^n$, 
we define the {\it Jackson integral of $f(t)$ over the lattice $\<\x\>$} as follows:
\begin{equation}
\label{eq:1.1.2}
\int_{\<\x\>}
  f(t)\ \frac{d_q t_1}{t_1}\wedge\cdots\wedge\frac{d_q t_n}{t_n}
  :=(1-q)^n
  \sum_{(\nu_1,\ldots,\nu_n)\in\Bbb Z^{n}}
  f(q^{\nu_1}\xi_1,\ldots,q^{\nu_n}\xi_n),
\end{equation}
where
$$\<\x\>:=
\{(\x_1 q^{\nu_1},\ldots,\x_n q^{\nu_n})\in (\Bbb C^*)^n
\,;\, \nu_i\in \Bbb Z\ (i=1,\ldots,n)\}. $$
}\end{df}
We simply write 
$\varpi=\frac{d_q t_1}{t_1}\wedge\cdots\wedge\frac{d_q t_n}{t_n}$.
For any function $\vp(t)$ we use $\wtil\vp$ 
for the Jackson integral defined as 
\begin{equation*}
\wtil\vp=\wtil\vp(\x):=
\int_{\<\x\>}\vP(t)\vp(t)\varpi.
\end{equation*}
\subsection{Matsuo's basis and $q$-KZ equation} 
\par
Let $L$ denote a set of multi-indeces as 
$$
L=\{
(l_1,\ldots,l_m)\in \Bbb Z_{\ge 0};\ l_1+\cdots+l_m=n\}.
$$
For $l=(l_1,\ldots,l_m)\in L$, let $\psi_{l}(t)$ be the rational functions 
introduced by Matsuo in \cite{ma2} 
as follows:
\begin{align*}
\psi_{l}(t)
&=
\psi_{(l_1,\ldots,l_m)}
\Big[
  \begin{matrix}
      x_{1},\ldots,x_{m}\\
      \b_{1},\ldots,\b_{m}
  \end{matrix}
\Big]
(t_1,\ldots,t_n)
\\&
=\A\left[
                   \prod_{j=1}^m
                   \Bigg[
                   \prod_{k=1}^{l_j+\cdots+l_m}
                   \frac{1-q^{\b_{j-1}} t_k/x_{j-1}}{1-t_k/x_j}
                   \Bigg]
                   \cdot\prod_{1\le i<j\le n}(t_i-q^{-\c}t_j)
               \right]
\end{align*}
where $q^{\b_0}=0$ and 
$\cal A$ is an alternating sum such that 
$${\cal A} g(t):=\sum_{\sigma\in \frak S_n}\sgn \sigma\cdot \sigma g(t).$$
Let $T_j$ denote the $q$-shift operator defined by
$$
T_jF(x_1,\ldots,x_m)=F(x_1,\ldots,qx_j,\ldots,x_m).
$$
We can consider a system of linear ordinary $q$-difference equations
for a vector $${\Vec y}(x_1,\ldots,x_m)=(\wt\psi_l)_{l\in L}$$ 
in a tensor coordinate $(z_l)_{l\in L}$ 
as follows:
\begin{equation*}
T_j{\Vec y}(x_1,\ldots,x_m)={\Vec y}(x_1,\ldots,x_m)K_j(x_1,\ldots,x_m) \qquad (1\le j\le m), 
\end{equation*}
where $K_j(x_1,\ldots,x_m)$ is a suitable matrix function of order $\#L={n+m-1\choose m-1}$.  
\begin{thm}
[Matsuo {\rm\cite{ma2}}, Varchenko {\rm\cite{v}}
]
\label{thm:1.2.1}
The matrix $K_j(x_1,\ldots,x_m)$ is expressed as 
\begin{align*}
&K_j(x_1,\ldots,x_m)\nonumber\\
&\quad=R_{j,j+1}(\textstyle{x_j\over x_{j+1}}\displaystyle)
  R_{j,j+2}(\textstyle{x_j\over x_{j+2}}\displaystyle)
  \cdots 
  R_{j,m}(\textstyle{x_j\over x_{m}}\displaystyle)
  \ 
  D_j
  \ 
  R_{j,1}(\textstyle{qx_j\over x_{1}}\displaystyle)
  R_{j,2}(\textstyle{qx_j\over x_{2}}\displaystyle)
  \cdots 
  R_{j,j-1}(\textstyle{qx_j\over x_{j-1}}\displaystyle)
\end{align*}
where $D_j$ is a diagonal matrix defined by $$(z_l)_{l\in L}D_j:=(q^{(\ta-(n-1)\c)l_j }z_l)_{l\in L}$$
and 
$R_{i,j}(\textstyle{x_i\over x_{j}}\displaystyle)
=\bigoplus_{\nu=1}^{n} 
 \big[R_{i,j}^{(\nu)}(\textstyle{x_i\over x_{j}}\displaystyle)\big]^{{\nu+m-1\choose \nu}}$ 
changes the basis $\{\psi_l\}_{l\in L}$ according to the following rule: 
\begin{align*}
&
\Bigg(
\psi_{(l_{\s(1)},...,l_{\s(i)},l_{\s(i+1)},...,l_{\s(m)})}
\Big[
  \begin{matrix}
       x_{\s(1)},..., x_{\s(i)}, x_{\s(i+1)},..., x_{\s(m)}\\
      \b_{\s(1)},..., \b_{\s(i)}, \b_{\s(i+1)},..., \b_{\s(m)}
  \end{matrix}
\Big]
(t)
\Bigg)
_{\!\!\!l_{\s(i)}+l_{\s(i+1)}=\nu}
\hskip -15pt
R_{\s(i),\s(i+1)}^{(\nu)}(\textstyle{{x_{\s(i)}\over x_{\s(i+1)}}}\displaystyle)
\\[5pt]
&=
\Bigg(
\psi_{(l_{\s(1)},...,l_{\s(i+1)},l_{\s(i)},...,l_{\s(m)})}
\Big[
  \begin{matrix}
      x_{\s(1)},..., x_{\s(i+1)}, x_{\s(i)},..., x_{\s(m)}\\
      \b_{\s(1)},..., \b_{\s(i+1)}, \b_{\s(i)},..., \b_{\s(m)}
  \end{matrix}
\Big]
(t)
\Bigg)
_{\!\!\!l_{\s(i+1)}+l_{\s(i)}=\nu}
\end{align*}
for $\s\in \frak S_m$.
\end{thm}

\begin{rmk}
{\rm
\label{rmk:1.2.2}
The matrix 
$R_{i,j}^{(\nu)}(\textstyle{x_i\over x_{j}}\displaystyle)$
of rank $\nu+1$ coincides with the matrix 
$R_{i,j}(\textstyle{x_i\over x_{j}}\displaystyle)$ for $\vP_{\nu,2}(t)$ 
(see \cite[Theorem 3.5.10]{v}
). 
Therefore, in order to get an explicit expression of $q$-difference equations for $\vP_{n,m}(t)$
it suffices to know the matrix $R_{i,j}(\textstyle{x_i\over x_{j}}\displaystyle)$ in the case $m=2$.
From now on we will consider $m=2$. 
}
\end{rmk}
\subsection{Gauss decomposition of $R$-matrix }
\par
When $m=2$ we have 
$$
\vP(t)=
\vP_{n,2}(t):=
t_1^{\a_1}\cdots t_n^{\a_n}
\prod_{j=1}^n 
             \frac{(t_j/x_1)_\8}{(t_j q^{\b_1}/x_1)_\8}
             \frac{(t_j/x_2)_\8}{(t_j q^{\b_2}/x_2)_\8}
       \prod_{1\le i<j\le n}
             \frac{(q^{1-\c}t_j/t_i)_\8}{(q^\c t_j/t_i)_\8}
$$
where $\a_j=\ta -n+j -2(j-1)\c$.
Let $\tau$ be the operation which exchanges $x_1$, $\b_1$ for $x_2$, $\b_2$ respectively.
We set 
\begin{equation}
\label{eq:1.3.1}
\psi_s(t)
:=\psi_{(s,n-s)}
\Big[
  \begin{matrix}
      x_{1},x_{2}\\
      \b_{1},\b_{2}
  \end{matrix}
\Big]
(t),\quad
\vp_s(t)
:=
\t\psi_{(s,n-s)}
=\psi_{(s,n-s)}
\Big[
  \begin{matrix}
      x_{2},x_{1}\\
      \b_{2},\b_{1}
  \end{matrix}
\Big]
(t)
\end{equation}
for $0\le s\le n$. 
These two bases are connected by the matrix $R_{i,j}(\textstyle{x_i\over x_j}\displaystyle)$ via 
\begin{equation*}
(\psi_n(t),\psi_{n-1}(t),\ldots,\psi_0(t))
R_{1,2}(\textstyle{x_1\over x_{2}}\displaystyle)=
(\vp_0(t),\vp_1(t),\ldots,\vp_n(t))
\end{equation*}
and
\begin{equation*}
(\vp_0(t),\vp_1(t),\ldots,\vp_n(t))
R_{2,1}(\textstyle{x_2\over x_1}\displaystyle)=
(\psi_n(t),\psi_{n-1}(t),\ldots,\psi_0(t))
,
\end{equation*}
so that 
\begin{equation}
\label{eq:1.3.4}
R_{1,2}(\textstyle{x_1\over x_{2}}\displaystyle)R_{2,1}(\textstyle{x_2\over x_{1}}\displaystyle)=I
\quad \hbox{and}\quad
R_{1,2}(\textstyle{x_1\over x_{2}}\displaystyle)=J\ \t R_{2,1}(\textstyle{x_2\over x_{1}}\displaystyle)\ J,
\end{equation}
where $I$ is the identity matrix and $J$ is the matrix $(\d_{i,n-j})_{i,j=0}^n$. 
The $q$-KZ equations are 
\begin{align}
T_1(\wtil\psi_n,\wtil\psi_{n-1},\ldots,\wtil\psi_0)
&=(\wtil\psi_n,\wtil\psi_{n-1},\ldots,\wtil\psi_0)
K_1(x_1,x_2)
,
\label{eq:1.3.5}
\\
T_2(\wtil\psi_n,\wtil\psi_{n-1},\ldots,\wtil\psi_0)
&=(\wtil\psi_n,\wtil\psi_{n-1},\ldots,\wtil\psi_0)
K_2(x_1,x_2).
\label{eq:1.3.6}
\end{align}

\begin{thm}[Matsuo {\rm \cite{ma2}}
]
\label{thm:1.3.1}
The matrices $K_1(x_1,x_2)$ and $K_2(x_1,x_2)$ are expressed as 
$$K_1(x_1,x_2)=R_{1,2}(\textstyle{x_1\over x_{2}}\displaystyle)D_1 
\quad \hbox{and}\quad 
  K_2(x_1,x_2)=D_2 R_{2,1}(\textstyle{qx_2\over x_{1}}\displaystyle),$$ 
where $D_1=\diag[q^{(\ta-(n-1)\c)(n-s)}]_{s=0}^n$ and $D_2=J D_1J$. 
\end{thm}
\par
By \eqref{eq:1.3.1}, the expression \eqref{eq:1.3.5} is equal to the following:
\begin{align}
T_2(\wtil\vp_0,\wtil\vp_{1},\ldots,\wtil\vp_n)
&=\t T_1(\wtil\psi_n,\wtil\psi_{n-1},\ldots,\wtil\psi_0) J 
 =\t (\wtil\psi_n,\wtil\psi_{n-1},\ldots,\wtil\psi_0)\t K_1(x_1,x_2) J 
\nonumber\\
&=(\wtil\vp_0,\wtil\vp_{1},\ldots,\wtil\vp_n)\ J \t R_{1,2}(\textstyle{x_1\over x_{2}}\displaystyle)J\ J D_1J
\nonumber\\
&=(\wtil\vp_0,\wtil\vp_{1},\ldots,\wtil\vp_n)R_{2,1}(\textstyle{x_2\over x_1}\displaystyle) D_2.
\label{eq:1.3.7}
\end{align}
\begin{rmk}
\label{rmk:1.3.2}
{\rm
The expressions \eqref{eq:1.3.5} and \eqref{eq:1.3.6} are essentially the same 
because, by using \eqref{eq:1.3.7}, the equation \eqref{eq:1.3.6} is deduced from \eqref{eq:1.3.5}. 
}
\end{rmk}
\par
Therefore we now take a basis $\{\vp_s(t);0\le s \le n\}$ with
\begin{equation}
\label{eq:1.3.8}
\vp_s(t)
=\A
           \left[
                   \prod_{k=1}^n
                   \frac{1}{1-t_k/x_2}
                   \prod_{k=1}^{n-s}
                   \frac{1-q^{\b_2} t_k/x_2}{1-t_k/x_1}
                   \prod_{1\le i<j\le n}(t_i-q^{-\c}t_j)
           \right] 
\end{equation}
and a $q$-difference system
\begin{equation}
\label{eq:1.3.9}
T_2(\wtil\vp_0,\wtil\vp_{1},\ldots,\wtil\vp_n)
=(\wtil\vp_0,\wtil\vp_{1},\ldots,\wtil\vp_n)
K(x_1,x_2).
\end{equation}
From \eqref{eq:1.3.7}, we have $K(x_1,x_2)=R_{2,1}(\textstyle{x_2\over x_{1}}\displaystyle)D_2$
and $D_2=\diag[q^{(\ta-(n-1)\c)s}]_{s=0}^n$. 
%
\begin{thm}
\label{thm:1.3.3}
The matrix $R_{2,1}(\textstyle{x_2\over x_{1}}\displaystyle)$ 
admits the following Gauss decomposition:
$$R_{2,1}(\textstyle{x_2\over x_{1}}\displaystyle)
=U_R\cdot D_R\cdot L_R ,$$
where
$U_R=(u_{\!R,rs})$, $s\ge r$, is an upper triangular matrix,
\begin{align*}
u_{\!R,rs}&=(-1)^{s-r}(\X21 q^{-\b_2})^{(s-r)}q^{-(s-r)(s+r-1)\c/2}\\
        &\quad\cdot\frac{(q^\c;q^\c)_s}{(q^\c;q^\c)_{s-r}\cdot (q^\c;q^\c)_r}\cdot
        \frac{(q^{\b_1+(n-s)\c};q^\c)_{s-r}}{(\X21 q^{\b_1-\b_2+(n-2s+1)\c};q^\c)_{s-r}},
\end{align*}
\par
\noindent
$D_R=\diag[d_{\!R,0},\ldots,d_{\!R,n}]$,
\begin{align*}
d_{\!R,r}    
&=q^{-(n-r)\b_2-r(n-r)\c}
      \cdot 
      \frac{(\X21 q^{\b_1};q^\c)_{n-r}}{(\X21 q^{-\b_2-(r-1)\c};q^\c)_{r}}
      \cdot 
      \frac{(\X21 q^{\b_1-\b_2+(n-2r+1)\c};q^\c)_{r}}{(\X21 q^{\b_1-\b_2-r\c};q^\c)_{n-r}},
\end{align*}
\par
\noindent
and 
$L_R=(l_{\!R,rs})$, $r\ge s$, is a lower triangular matrix,
\begin{align*}
l_{\!R,rs}&=(-1)^{r-s}
           q^{\b_2(s-r)}q^{-(r-s)(r+s-1)\c/2}\\
        &\quad\cdot\frac{(q^\c;q^\c)_{n-s}}{(q^\c;q^\c)_{r-s}\cdot (q^\c;q^\c)_{n-r}}\cdot
        \frac{(q^{\b_2+s\c};q^\c)_{r-s}}{(\X21 q^{\b_1-\b_2+(n-2r+1)\c};q^\c)_{r-s}}.
\end{align*}
\end{thm}
\par We have another Gauss decomposition expression. 
\begin{thm}
\label{thm:1.3.4}
The matrix $R_{2,1}(\textstyle{x_2\over x_{1}}\displaystyle)$ 
admits the following Gauss decomposition:
$$R_{2,1}(\textstyle{x_2\over x_{1}}\displaystyle)
=L'_R\cdot D'_R\cdot U'_R ,$$
where 
$L'_R=(l'_{\!R,rs})$, $r\ge s$, is a lower triangular matrix,
\begin{align*}
l'_{\!R,rs}
&=(-1)^{r-s}
            q^{-(r-s)(r+s-1)\c/2}
          \cdot\frac{(q^\c;q^\c)_{n-s}}{(q^\c;q^\c)_{r-s}\cdot (q^\c;q^\c)_{n-r}}\cdot
        \frac{(q^{\b_2+s\c};q^\c)_{r-s}}{(\X21 q^{-(n-2s-1)\c};q^\c)_{r-s}},
\end{align*}
\par
\noindent
$D'_R=\diag[d'_{\!R,0},\ldots,d'_{\!R,n}]$,
\begin{align*}
d'_{\!R,r}    
&=q^{-(n-r)(\b_2+r\c)}
      \cdot 
      \frac{(\X21 q^{\b_1};q^\c)_{r}}{(\X21 q^{-\b_2-(n-r-1)\c};q^\c)_{n-r}}
      \cdot 
      \frac{(\X21 q^{-(n-2r-1)\c};q^\c)_{n-r}}{(\X21 q^{-(n-r)\c};q^\c)_{r}},
\end{align*}
\par
\noindent
and 
$U'_R=(u'_{\!R,rs})$, $s\ge r$, is an upper triangular matrix,
\begin{align*}
u'_{\!R,rs}&=q^{-(s-r)r\c}
        \cdot\frac{(q^\c;q^\c)_s}{(q^\c;q^\c)_{s-r}\cdot (q^\c;q^\c)_r}\cdot
        \frac{(q^{\b_1+(n-s)\c};q^\c)_{s-r}}{(\textstyle{x_1\over x_2}\displaystyle q^{(n-r-s)\c};q^\c)_{s-r}}.
\end{align*}
\end{thm}
\par\noindent{\bf Proof.} See Appendix \ref{section:9}. \hfill\qed
%
\begin{rmk}
\label{rmk:1.3.5}
{\rm
We can also derive Theorem \ref{thm:1.3.3} 
by using the method of the Riemann--Hilbert problem 
for $q$-difference equation \eqref{eq:1.3.9} from a connection matrix (see Subsection \ref{subsection:4.2}).
}
\end{rmk}

\begin{rmk}
\label{rmk:1.3.6}
{\rm
From \eqref{eq:1.3.4}, we have two kind of expressions 
\begin{equation}
\label{eq:1.3.10}
R_{1,2}(\textstyle{x_1\over x_{2}}\displaystyle)
   =J \t U_R J \cdot J \t D_R J \cdot J \t L_R J
   =J \t L'_R J \cdot J \t D'_R J \cdot J \t U'_R J.
\end{equation}
}
\end{rmk}

\begin{rmk}
\label{rmk:1.3.7}
{\rm
From \eqref{eq:1.3.4} and \eqref{eq:1.3.10}, we have 
$$
L_R^{-1}=J \t U_R J,\quad D_R^{-1}= J \t D_R J,\quad U_R^{-1}=J \t L_R J,
$$
\vskip -4mm
\noindent 
and
$$
U'_R{}^{-1}=J \t L'_R J,\quad D'_R{}^{-1}= J \t D'_R J,\quad L'_R{}^{-1}=J \t U'_R J.
$$
}
\end{rmk}
\begin{cor}
\label{cor:1.3.8}
The matrix $K(x_1,x_2)$ 
that represents the $q$-difference system \eqref{eq:1.3.9} is given by 
\begin{align}
K(x_1,x_2)
&=U_R\cdot D_R\cdot L_R \cdot \diag[q^{(\ta-(n-1)\c)s}]_{s=0}^n
\label{eq:1.3.11}
\\
&=L'_R\cdot D'_R\cdot U'_R \cdot \diag[q^{(\ta-(n-1)\c)s}]_{s=0}^n.
\nonumber
\end{align}
\end{cor}
\begin{cor}
\label{cor:1.3.9}
The determinant of the matrices $R_{2,1}(\textstyle{x_2\over x_{1}}\displaystyle)$ and $K(x)$ are given by the following expressions:
\begin{align}
\det R_{2,1}(\textstyle{x_2\over x_{1}}\displaystyle)
              &=\prod_{r=0}^n
                   q^{-(n-r)(\b_2+r\c)}
                  \frac{(\X21 q^{\b_1};q^\c)_{r}}{(\X21 q^{-\b_2-(r-1)\c};q^\c)_{r}},
\nonumber\\[5pt]
\det K(x_1,x_2)&=q^{(\ta-(n-1)\c)n(n+1)/2}\det R_{2,1}(\textstyle{x_2\over x_{1}}\displaystyle).
\label{eq:1.3.12}
\end{align}
\end{cor}
\par\noindent{\bf Proof.} By Theorem \ref{thm:1.3.3} or Theorem \ref{thm:1.3.4}, we have  
\begin{align*}
\det R_{2,1}(\textstyle{x_2\over x_{1}}\displaystyle)
=\det D_R \ (\hbox{or }=\det D'_R ).
\end{align*}
The result now follows from the following identity:
\begin{equation*}
\prod_{r=0}^n  \frac{(\X21 q^{-(n-2r-1)\c};q^\c)_{n-r}}{(\X21 q^{-(n-r)\c};q^\c)_{r}}
=\prod_{r=0}^n  \frac{(\X21 q^{\b_1-\b_2+(n-2r+1)\c};q^\c)_{r}}{(\X21 q^{\b_1-\b_2-r\c};q^\c)_{n-r}}
=1.
\end{equation*}
\hfill $\square$
\begin{rmk}
\label{rmk:1.3.10}
{\rm
Eq.\,\eqref{eq:1.3.12} 
was conjectured by Mimachi in \cite{mi2} 
and another proof for it was given 
by Aomoto and Kato in \cite{ak6}. 
}
\end{rmk}

%
\section{Review of Riemann--Hilbert problem for $q$-difference equation from connection matrices}
\label{section:2}
\par
We recall the notion of the Riemann--Hilbert problem for $q$-difference equation 
following \cite{ao1,ao4, ao5, ak6}. 
\par 
We consider a linear ordinary $q$-difference equation for a vector function 
${\Vec y}(z)= (y_0(z),\ldots,y_n(z))$, $z\in \Bbb C^*$, satisfying
\begin{equation*}
{\Vec y}(qz)={\Vec y}(z)A(z),
\end{equation*}
where $A(z)$ is a suitable rational matrix function of order $n+1$. 
We now assume the following conditions for the matrix $A(z)$;
\begin{itemize}
\item[P1)]
The matrix $A(z)$ is holomorphic at $z=0$ and $z=\8$.
\item[P2)]
$A(0)=\lim\limits_{z\to 0}A(z)$ and $A(\8)=\lim\limits_{z\to \8}A(z)$ are diagonalizable, i.e., 
there exist matrices $C^+$, $C^-$ and 
diagonal matrices 
$D^+=\diag[\mu_0,\ldots,\mu_n]$,
$D^-=\diag[\mu^*_0,\ldots,\mu^*_n]$ 
such that 
\begin{align}
A(0)&=(C^+)^{-1}q^{D^+}C^+,
\label{eq:2.2}\\
A(\8)&=(C^-)^{-1}q^{D^-}C^-,
\label{eq:2.3}
\end{align}
where $q^{D^+}=\diag[q^{\mu_0},\ldots,q^{\mu_n}]$
and $q^{D^-}=\diag[q^{\mu^*_0},\ldots,q^{\mu^*_n}]$.
\item[P3)]
The diagonal elements of $D^+$ and $D^-$ satisfy the following {\it non-resonance} condition:%
\footnote{\ P3) is a condition that should be referred to as H2) in the reference \cite{ao4}, but the symbol H2) does not appear in \cite{ao4} due to a typo. It should be noted that equations (2) and (3) in \cite{ao4} 
actually correspond to the condition H2).}
\begin{align*}
q^{\mu_i-\mu_j}&\not=q^{\pm 1}, q^{\pm 2},\ldots,\\
q^{\mu^*_i-\mu^*_j}&\not=q^{\pm 1}, q^{\pm 2},\cdots.
\end{align*}
\item[P4)]
The matrix $A(z)$ does not depend on $q$.
\end{itemize}
Under these conditions P1), P2) and P3), 
from the classical theorem due to G.\,D.\,Birkhoff (see \cite{ao4, ao5}), 
we know that 
there exists a unique solution $Y_0(z)$ of the equation 
\begin{equation}
\label{eq:2.4}
Y(qz)=Y(z)A(z)
\end{equation}
such that $Y_0(z)$ satisfies the asymptotic behavior
\begin{align*}
Y_0(z)&\sim (C^+)^{-1}z^{D^+}C^+  \quad \hbox{at} \quad z=0,
\end{align*}
and we also know that there exists a unique solution $Y_\8(z)$ of the equation \eqref{eq:2.4}
such that $Y_\8(z)$ satisfies the asymptotic behavior
\begin{align*}
Y_\8(z)&\sim (C^-)^{-1}z^{D^-}C^- \quad \hbox{at} \quad z=\8.
\end{align*}
We call $Y_0(z)$ and $Y_\8(z)$ the {\it fundamental solutions} 
of \eqref{eq:2.4} at $z=0$ and $z=\8$ respectively. 
The connection matrix $P(z)$ between the fundamental solutions $Y_0(z)$ and $Y_\8(z)$ is defined by 
$$
P(z):=Y_0(z)Y_\8(z)^{-1}.
$$
For an arbitrary matrix $X(q)$ depending on $q$, we denote the limit for $q\to 0$ as 
$$
(X)_0 =\lim_{q\to 0}X(q).
$$
\begin{thm}[Aomoto's lemma]
\label{thm:2.1}
In addition to {\rm P1), P2)} and {\rm P3)}, under the condition {\rm P4)}, the following limit formula holds:
\begin{equation}
\label{eq:2.5}
\big(P(z)\big)_{\!0}=A(0)A(z)^{-1}.
\end{equation}
\end{thm}
\par\noindent{\bf Proof.}
See \cite{ao4, ao5}. 
\hfill$\square$\\[5pt]
\indent
This theorem will play a crucial role in calculating $A(z)$ explicitly in the following sections.
\begin{rmk}
\label{rmk:2.2}
{\rm
In the condition P2), if we can choose the matrices 
$C^+=\big(c^+_{rs}\big)_{r,s=0}^n$ and $C^-=\big(c^-_{rs}\big)_{r,s=0}^n$ 
as triangular matrices, 
the unipotent matrices 
$\diag[(c^+_{rr})^{-1}]_{r=0}^n C^+$ and $\diag[(c^-_{rr})^{-1}]_{r=0}^n C^-$ also satisfy P2). 
Thus we assume that the matrices $C^+$ and $C^-$ are unipotent if they are triangular. 
}
\end{rmk}
%
\section{Characteristic cycles and fundamental solutions}
\par
Let $F_r^{n-r}$, $0\le r\le n$, denote the partition of the set $\{1,\ldots, n\}$ into subsets 
$\{1,\ldots, n-r\}$ and $\{n-r+1,\ldots, n\}$.
Let $\x_{F^{n-r}_r}=(\x_1,\ldots,\x_n)$, $\y_{F^{n-r}_r}=(\y_1,\ldots,\y_n)$, 
    $\z_{F^{n-r}_r}=(\z_1,\ldots,\z_n)$ and $\d_{F^{n-r}_r}=(\d_1,\ldots,\d_n)$
be the four points in $(\Bbb C^*)^n$ defined by
\begin{eqnarray*}
\x_{F^{n-r}_r} 
 &:& 
 \bigg\{
                   \begin{matrix}
                      \x_1=x_1,&  \x_2=x_1q^{\c},& \ldots,&\x_{n-r}=x_1q^{(n-r-1)\c}, \\
                      \x_{n-r+1}=x_2, & \x_{n-r+2}=x_2 q^{\c},& \ldots,&\x_{n}=x_2 q^{(r-1)\c}, 
                   \end{matrix}
\\[3pt]
\y_{F^{n-r}_r} 
 &:& 
 \bigg\{
                   \begin{matrix}
                      \y_1=x_1 q^{-\b_1},&  \y_2=x_1 q^{-\b_1-\c},& \!\!\ldots,&\y_{n-r}=x_1 q^{-\b_1-(n-r-1)\c}, \\
                      \y_{n-r+1}=x_2 q^{-\b_2}, & \y_{n-r+2}=x_2 q^{-\b_2-\c},& \!\!\ldots,&\y_{n}=x_2 q^{-\b_2-(r-1)\c}, 
                   \end{matrix}
\\[3pt]
\z_{F^{n-r}_r} 
 &:& 
 \bigg\{
                   \begin{matrix}
                      \z_1=x_1 q^{-\b_1},&  \z_2=x_1 q^{-\b_1-\gamma},& \ldots,&\z_{n-r}=x_1 q^{-\b_1-(n-r-1)\gamma}, \\
                      \z_{n-r+1}=x_2, & \z_{n-r+2}=x_2 q^{\gamma},& \ldots,&\z_{n}=x_2 q^{(r-1)\gamma}, 
                   \end{matrix}
\\[3pt]
\d_{F^{n-r}_r} 
 &:& 
 \bigg\{
                   \begin{matrix}
                      \d_1=x_1,&  \d_2=x_1 q^{\c},& \ldots,&\d_{n-r}=x_1 q^{(n-r-1)\gamma}, \\
                      \d_{n-r+1}=x_2 q^{-\b_2}, & \d_{n-r+2}=x_2 q^{-\b_2-\c},& \ldots,&\d_{n}=x_2 q^{-\b_2-(r-1)\c}. 
                   \end{matrix} 
\end{eqnarray*}
We call the lattices $\<\x_{F^{n-r}_r} \>$, $\<\eta_{F^{n-r}_r} \>$, $\<\z_{F^{n-r}_r} \>$ 
and $\<\d_{F^{n-r}_r} \>$, $0\le r\le n$, the {\it characteristic cycles}.
Since the ordinary Jackson integrals over 
the cycles $\<\eta_{F^{n-r}_r} \>$, $\<\z_{F^{n-r}_r} \>$ and $\<\d_{F^{n-r}_r} \>$ 
diverge, 
we have to define the {\it regularized Jackson integrals} for them as follows:
\begin{align*}
&\int_{\<\eta_{F^{n-r}_r}\>}
\vP(t)\vp(t)\varpi
:=(1-q)^n\sum_{(\nu_1,\ldots,\nu_n)\in\Bbb Z^n}
\Res_{\ 
       \substack{ 
         \substack{ t_1=\eta_1q^{\nu_1}, \\[1.5pt] \cdots\cdots,}
         \\
          t_{n}=\eta_{n}q^{\nu_{n}} 
       }
     }
\vP(t_1,\ldots,t_n) \vp(t_1,\ldots,t_n) 
\ \frac{d t_1}{t_1}\wedge\cdots\wedge\frac{d t_n}{t_n},\\[3pt]
&\int_{\<\z_{F^{n-r}_r}\>}
  \vP(t)\vp(t)\varpi
:=(1-q)^n \sum_{(\nu_1,\ldots,\nu_n)\in\Bbb Z^n}
\Res_{\ 
        \substack{ 
         \substack{ t_1=\z_1q^{\nu_1}, \\[1.5pt] \cdots\cdots,}
         \\
          t_{n-r}=\z_{n-r}q^{\nu_{n-r}} 
       }
     }
\vP(t_1,\ldots,t_{n-r},q^{\nu_{n-r}}\xi_{n-r},\ldots,q^{\nu_n}\z_n)\\
&\hskip 190pt \cdot 
\vp(t_1,\ldots,t_{n-r},q^{\nu_{n-r}}\xi_{n-r},\ldots,q^{\nu_n}\z_n)
\ \frac{d t_1}{t_1}\wedge\cdots\wedge\frac{d t_{n-r}}{t_{n-r}},
\end{align*}
and
\begin{align*}
&\int_{\<\d_{F^{n-r}_r}\>}
  \vP(t)\vp(t)\varpi\\
&\hspace{50pt}:=(1-q)^n \sum_{(\nu_1,\ldots,\nu_n)\in\Bbb Z^n}
\Res_{\
       \substack{ 
       \substack{t_{n-r+1}=\d_{n-r+1}q^{\nu_{n-r+1}}, 
             \\[1.5pt]
        \!\!\!\!\!\!\!\!\!\!\!\!\!\!\!\!
                    \cdots\cdots,}
         \\ 
        \!\!\!\!\!\!\!\!\!\!\!\!\!\!\!\!
          t_{n}=\d_n q^{\nu_{n}} 
       }
     }
\vP(q^{\nu_1}\d_{1},\ldots,q^{\nu_{n-r}}\d_{n-r},t_{n-r+1},\ldots,t_{n})\\
&\hspace{160pt} \cdot 
\vp(q^{\nu_1}\d_{1},\ldots,q^{\nu_{n-r}}\d_{n-r},t_{n-r+1},\ldots,t_{n})
\ \frac{d t_{n-r+1}}{t_{n-r+1}}\wedge\cdots\wedge\frac{d t_n}{t_n}.
\end{align*}
Let $\vp_s(t)$ be the function defined by \eqref{eq:1.3.8} in Section \ref{section:1}.
We consider a linear ordinary $q$-difference equation for a vector function 
${\Vec y}(q^{\ta})=(\wt \vp_0,\wt \vp_1,\ldots,\wt \vp_n)$
satisfying 
\begin{equation*}
{\Vec y}(q^{\ta+1})={\Vec y}(q^{\ta})A(q^{\ta}).
\end{equation*}
where $A(q^{\ta})$ is a suitable rational matrix function of order $n+1$. 
We define the following two matrices
\begin{equation}
\label{eq:3.1}
Y_\x:=\bigg(\int_{\<\x_{F^{n-r}_r}\>}\vP(t)\vp_s(t)\varpi\bigg)_{r,s=0}^n,
\quad
Y_\y:=\bigg(\int_{\<\y_{F^{n-r}_r}\>}\vP(t)\vp_s(t)\varpi\bigg)_{r,s=0}^n,
\end{equation}
which are solutions of the matrix equation 
\begin{equation}
\label{eq:3.2}
Y(q^{\ta+1})=Y(q^{\ta})A(q^{\ta}).
\end{equation}
We set $(y)^{\ta}:=y_1^{\ta}\cdots y_n^{\ta}$ for $y=(y_1,\ldots, y_n)$. 
The solutions $Y_\x$ and $Y_\y$ have the following asymptotic behaviors:
\begin{align}
Y_\x&\sim (q^{\ta})^{D_{\!A}^+}C_{\!A}^+ \quad\hbox{ at }\quad\ta\to +\8,
\label{eq:3.3}
\\
Y_\y&\sim (q^{\ta})^{D_{\!A}^-}C_{\!A}^- \quad\hbox{ at }\quad\ta\to -\8
\label{eq:3.4}
\end{align}
where 
$$
(q^{\ta})^{D_{\!A}^+}
=\diag[(\x_{F^{n-r}_r} )^{\ta}]_{r=0}^n,\quad
(q^{\ta})^{D_{\!A}^-}
=\diag[(\y_{F^{n-r}_r} )^{\ta}]_{r=0}^n,
$$
and
$C_{\!A}^+=(c_{\!A,rs}^+)_{r,s=0}^n$ and $C_{\!A}^-=(c_{\!A,rs}^-)_{r,s=0}^n$ 
are matrices not depending on $q^{\ta}$ defined by 
\begin{equation}
\label{eq:3.5}
c_{\!A,rs}^+
:=(1-q)^n \frac{\vP(\x_{F^{n-r}_r} )\vp_s(\x_{F^{n-r}_r} )}
              {(\x_{F^{n-r}_r} )^{\ta}},\quad
c_{\!A,rs}^-
:=(1-q)^n \!\!\!
\Res_{t=\y_{F^{n-r}_r}}\!\!\!
         \frac{\vP(t)\vp_s(t)}
              {(t )^{\ta}}
\frac{d t_1}{t_1}\wedge\cdots\wedge\frac{d t_n}{t_n}. 
\end{equation}
Since 
\begin{align*}
q^{D_{\!A}^+}&=\diag[x_1^{n-r}x_2^r\ q^{[r(r-1)+(n-r)(n-r-1)]\c/2}]_{r=0}^n,\\[5pt]
q^{D_{\!A}^-}&=\diag[x_1^{n-r}x_2^r\ q^{-(n-r)\b_1-r\b_2-[r(r-1)-(n-r)(n-r-1)]\c/2}]_{r=0}^n,
\end{align*}
the condition P3) is satisfied for generic parameters $x_1$, $x_2$, $\b_1$, $\b_2$  and $\c$.
\begin{prop}
\label{prop:3.1}
The matrix $C_{\!A}^+$ is lower triangular. 
The nonzero elements of the matrix $(C_{\!A}^+)_0$ are the following:
\begin{align*}
(c_{\!A,rr}^+)_0 
&=q^{-[r(r-1)+(n-r)(n-r-1)]\c/2}
\frac
     { (\X21 q^{-(n-r)\c};q^\c)_{r}\cdot (\X21 q^{-(n-r-1)\c};q^\c)_{r}}
     {(q^{\b_1};q^\c)_{n-r}\cdot (\X21 q^{\b_1};q^\c)_r\cdot (q^{\b_2};q^\c)_r \cdot (\X21 q^{\c};q^\c)_r},
\\[3pt]
\Big(\frac{c_{\!A,rs}^+}{c_{\!A,rr}^+}\Big)_0 
&=
\frac{c_{\!A,rs}^+}{c_{\!A,rr}^+}
=q^{-(r-s)(n-r)\c}
\frac{(q^\c;q^\c)_{n-s}}
     {(q^\c;q^\c)_{r-s}(q^\c;q^\c)_{n-r}}
\frac{(q^{\b_2+s\c};q^\c)_{r-s}}
     {(\X21 q^{(r+s-n)\c};q^\c)_{r-s}}
\quad \hbox{for}\quad r\ge s.
\end{align*}
\end{prop}
\par\noindent
{\bf Proof.} See Appendix \ref{section:7}.\hfill$\square$
\begin{prop}
\label{prop:3.2}
The matrix $C_{\!A}^-$ is upper triangular. 
The non-zero elements of the matrix $(C_{\!A}^-)_0$ are the following:
\begin{align*}
(c_{\!A,rr}^-)_0\
&=(-1)^n\frac{(q^{-\c};q^{-\c})_{n-r}\cdot(q^{-\c};q^{-\c})_{r}}{(1-q^{-\c})^{n}}
\cdot
  \frac{(\X21 q^{-\b_2-(r-1)\c};q^{\c})_{r}}{(\X21 q^{\b_1-\b_2-(r-1)\c};q^{\c})_{r}}
\\&\hskip 160pt
\cdot\frac{1}{{}_n(\X21 q^{\b_1-\b_2-(r-1)\c};q^\c)_{r}\cdot {}_n(\X21 q^{\b_1-\b_2-r\c};q^\c)_{r}},\\
\Big(
     \frac{c_{\!A,rs}^-}{c_{\!A,rr}^-}
\Big)_0
&=
     \frac{c_{\!A,rs}^-}{c_{\!A,rr}^-}\\
&=
(\X21 q^{-\b_2-r\c})^{s-r}
  \frac{(q^{\b_1+(n-s)\c};q^{\c})_{s-r}}
       {(\X21 q^{\b_1-\b_2+(n-s-r)\c};q^{\c})_{s-r}}
\cdot
     \frac{(q^{\c};q^{\c})_s}{(q^{\c};q^{\c})_r\cdot (q^{\c};q^{\c})_{s-r}}
\quad \hbox{ for }\quad r\le s.
\end{align*}
\end{prop}
\par\noindent
{\bf Proof.} See Appendix \ref{section:8}. \hfill$\square$

If we define 
$$
Y_0:=(C_{\!A}^+)^{-1} \ Y_\x,\quad 
Y_\8:=(C_{\!A}^-)^{-1} \ Y_\y,
$$
then the matrices $Y_0=Y_0(q^{\ta})$ and $Y_\8=Y_\8(q^{\ta})$ are also solutions of the equation \eqref{eq:3.2} 
and satisfy the following asymptotic behaviors:
\begin{align*}
Y_0(q^{\ta})&\sim (C_{\!A}^+)^{-1}(q^{\ta})^{D_{\!A}^+}C_{\!A}^+ \quad\hbox{at}\quad\ta\to +\8,\\
Y_\8(q^{\ta})&\sim (C_{\!A}^-)^{-1}(q^{\ta})^{D_{\!A}^-}C_{\!A}^- \quad\hbox{at}\quad\ta\to -\8.
\end{align*}
This implies that the matrix $A(q^{\ta})$ satisfies the condition P1) and P2) 
for the unipotent matrices 
$$
C^+=\Bigg(\frac{c_{\!A,rs}^+}{c_{\!A,rr}^+}\Bigg)_{r,s=0}^n,
\quad  
C^-=\Bigg(\frac{c_{\!A,rs}^-}{c_{\!A,rr}^-}\Bigg)_{r,s=0}^n
$$
in \eqref{eq:2.2} and \eqref{eq:2.3}, i.e., 
\begin{align}
A(0)&= (C_{\!A}^+)^{-1}q^{D_{\!A}^+}C_{\!A}^+ = (C^+)^{-1}q^{D_{\!A}^+}C^+, 
\label{eq:3.6}
\\
A(\8)&= (C_{\!A}^-)^{-1}q^{D_{\!A}^-}C_{\!A}^-=(C^-)^{-1}q^{D_{\!A}^-}C^-. 
\label{eq:3.7}
\end{align}
\par
Let $G=G(\ta,\b_1,\b_2,x_1,x_2)$ be the connection matrix between $Y_\x$ and $Y_\y$ defined by 
\begin{equation}
\label{eq:3.8}
G:=Y_\x \ Y_\y^{-1},
\end{equation}
which coincides with what Aomoto--Kato called the {\it principal connection matrix} (see Section \ref{section:4}), 
and then 
\begin{equation}
\label{eq:3.9}
Y_0(q^{\ta})Y_\8(q^{\ta})^{-1}
=(C_{\!A}^+)^{-1}\ Y_\x \ Y_\y^{-1}\ C_{\!A}^-
=(C_{\!A}^+)^{-1}\ G\ C_{\!A}^-.
\end{equation}
\par
Since the conditions P1), P2) and P3) are satisfied, 
by \eqref{eq:3.9} and Aomoto's lemma \eqref{eq:2.5}, 
we have 
$$
A(0)A(q^{\ta})^{-1}=\Big((C_{\!A}^+)^{-1}\ G\ C_{\!A}^-\Big)_{\!0}
$$
if the condition P4) holds for $A(q^{\ta})$. 
In Appendix \S \ref{section:6}, 
we will see that the condition P4) holds for $A(q^{\ta})$ .
From \eqref{eq:3.6}, we have 
\begin{equation}
\label{eq:3.10}
A(q^{\ta})=\Big((C_{\!A}^-)^{-1}\ G^{-1}\ q^{D_{\!A}^+}\ C_{\!A}^+\Big)_{\!0}.
\end{equation}

%
\section{Principal connection matrix}
\label{section:4}
\subsection{Principal connection matrix}
The set $\{\<\x_{F^{n-r}_r}\>;0\le r\le n\}$ makes up a basis of 
the dual space of $H^n(\bar X,\vP,\nabla)_{sym}$ (see the definition in \cite{ak5}). 
For an arbitrary $\x\in {\Bbb C}^*$, $\<\x\>$ is expressed uniquely as 
a linear combination of $\<\x_{F^{n-r}_r}\>$, $0\le r\le n$, in such a way that
\begin{equation}
\label{eq:4.1}
\<\x\>=\sum_{r=0}^n c_r\cdot \<\x_{F^{n-r}_r}\> 
\end{equation}
for some {\it pseudo-constants} $c_r$ 
in the sense that they do not change under the displacement
$\ta \mapsto \ta +1$, $\b_k\mapsto \b_k+1$, $\c\mapsto \c+1$ and $x\mapsto qx$. 
We denote the coefficient $c_r$ by $\big(\<\x\>;\<\x_{F^{n-r}_r}\>\big)_\vP$. 
Namely, \eqref{eq:4.1} means that 
\begin{equation}
\label{eq:4.2}
\int_{\<\x\>}\vP(t)\vp(t)\varpi
=\sum_{r=0}^n\ \big(\<\x\>;\<\x_{F^{n-r}_r}\>\big)_\vP
 \cdot \int_{\<\x_{F^{n-r}_r}\>}\hskip -4mm\vP(t)\vp(t)\varpi.
\end{equation}
For the other bases 
$\{\<\y_{F^{n-r}_r}\>;0\le r\le n\}$, 
$\{\<\z_{F^{n-r}_r}\>;0\le r\le n\}$ 
and 
$\{\<\d_{F^{n-r}_r}\>;0\le r\le n\}$,
we define the coefficients 
$\big(\<\x\>;\<\y_{F^{n-r}_r}\>\big)_\vP$, 
$\big(\<\x\>;\<\z_{F^{n-r}_r}\>\big)_\vP$ 
and 
$\big(\<\x\>;\<\d_{F^{n-r}_r}\>\big)_\vP$
in the same manner as above. 
By \eqref{eq:4.2} and \eqref{eq:3.8}, the elements of the principal connection matrix $G=(g_{rs})_{r,s=0}^n$ are written as follows:
$$
g_{rs}
=\big(\<\x_{F_r^{n-r}}\>;\<\y_{F_s^{n-s}}\>\big)_\vP
=\sum_{i=0}^n\big(\<\x_{F_r^{n-r}}\>;\<\z_{F_i^{n-i}}\>\big)_\vP
        \cdot\big(\<\z_{F_i^{n-i}}\>;\<\y_{F_s^{n-s}}\>\big)_\vP.
$$
\begin{thm}[Gauss decomposition {\rm \cite{ak5}}
]
\label{thm:4.1}
The principal connection matrix $G$ admits the following Gauss decomposition:
\begin{equation}
\label{eq:4.3}
G=(H_{\<\z;\x\>})^{-1}\ H_{\<\z;\y\>},
\end{equation}
where the matrices $H_{\<\z;\x\>}=\big(h^{{}^{\!++}}_{rs}\big)_{r,s=0}^n$ 
and 
$H_{\<\z;\y\>}=\big(h^{{}^{\!+-}}_{rs}\big)_{r,s=0}^n$ 
are defined by 
$$h^{{}^{\!++}}_{rs}:=\big(\<\z_{F_r^{n-r}}\>;\<\x_{F_s^{n-s}}\>\big)_\vP
, \quad
  h^{{}^{\!+-}}_{rs}:=\big(\<\z_{F_r^{n-r}}\>;\<\y_{F_s^{n-s}}\>\big)_\vP.
$$ 
The matrices $H_{\<\z;\x\>}$ and $H_{\<\z;\y\>}$ 
become an upper triangular matrix and a lower one respectively. 
An arbitrary element of the matrices 
$H_{\<\z;\x\>}$, $H_{\<\z;\y\>}$, $(H_{\<\z;\x\>})^{-1}$ and $(H_{\<\z;\y\>})^{-1}$
is expressed in a theta product form.
\end{thm}
\par\noindent{\bf Proof.} See Theorem 8.3 in \cite{ak5}. 
\hfill\qed\\
\par
Since we have already known the explicit form of $H_{\<\z;\x\>}$ and $H_{\<\z;\y\>}$ 
(see Theorem 8.1, 8.2 in \cite{ak5} 
and Lemma 13 in \cite{ak6}), 
in particular, we have 
\begin{align}
(h^{{}^{\!++}}_{rr})_0
&=(-1)^{n-r}
\frac
     { (q^\c;q^\c)_{n-r}
       \cdot (q^{\b_1};q^\c)_{n-r}
       \cdot ({\textstyle{x_1\over x_2}\displaystyle} q^{\b_2};q^\c)_{n-r}
      }
     {(1-q^\c)^{n-r}
       \cdot ({\textstyle{x_1\over x_2}\displaystyle};q^\c)_{n-r}
      }\nonumber\\
&\hskip 5mm \cdot
\frac
     {({\textstyle{x_1\over x_2}\displaystyle} q^{-\b_1-(n-r-1)\c};q^\c)_{n-r}
       \cdot (q^{\ta+\b_2-2(n-r-1)\c};q^\c)_{n-r}
      }
     {({\textstyle{x_1\over x_2}\displaystyle} q^{\b_2-\b_1-(n-r-1)\c};q^\c)_{n-r}
       \cdot 
      (q^{\ta+\b_2+\b_1-(n-r-1)\c};q^\c)_{n-r}
      }\nonumber\\
&\hskip 5mm \cdot
\frac
     {{}_n({\textstyle{x_1\over x_2}\displaystyle} q^{-(r-1)\c};q^\c)_{r}
        \cdot {}_n({\textstyle{x_1\over x_2}\displaystyle} q^{-r\c};q^\c)_{r}
      }
     {{}_n({\textstyle{x_1\over x_2}\displaystyle} q^{-(n-r-1)\c};q^\c)_{r}
        \cdot {}_n({\textstyle{x_1\over x_2}\displaystyle} q^{-(n-r)\c};q^\c)_{r}
      },
\label{eq:4.4}
\\[5pt]
\Big(
{h^{{}^{\!++}}_{rs}\over h^{{}^{\!++}}_{rr}}
\Big)_{\!0}
&=q^{-r(s-r)\c}
\frac
     {{}_n({\textstyle{x_1\over x_2}\displaystyle} q^{-(n-r-1)\c};q^\c)_{r}
        \cdot {}_n({\textstyle{x_1\over x_2}\displaystyle} q^{-(n-r)\c};q^\c)_{r}
        \cdot {}_n({\textstyle{x_1\over x_2}\displaystyle} q^{-(s-1)\c};q^\c)_{s}
      }
     {{}_n({\textstyle{x_1\over x_2}\displaystyle} q^{-(n-s-1)\c};q^\c)_{s}
        \cdot {}_n({\textstyle{x_1\over x_2}\displaystyle} q^{-(n-s)\c};q^\c)_{s}
        \cdot {}_n({\textstyle{x_1\over x_2}\displaystyle} q^{-(r-1)\c};q^\c)_{r}
      }\nonumber\\
&\hskip 5mm \cdot
\frac
     {({\textstyle{x_1\over x_2}\displaystyle} q^{-r\c};q^\c)_{n-r}}
     {({\textstyle{x_1\over x_2}\displaystyle} q^{-s\c};q^\c)_{n-s}}
\cdot
\frac
     {(q^{\b_2+r\c};q^\c)_{s-r}}
     {({\textstyle{x_1\over x_2}\displaystyle} q^{\b_2+(n-s)\c};q^\c)_{s-r}}\nonumber\\
&\hskip 5mm \cdot
\frac
     {(\X21 q^{-\ta-\b_2+(n-1)\c};q^\c)_{s-r}
         \cdot {}_s(q^\c;q^\c)_r 
      }
     {(q^{-\ta-\b_2+(2n-r-s-1)\c};q^\c)_{s-r}
         \cdot (\X21 q^{(r+s-n)\c};q^\c)_{s-r}
      }
\quad\hbox{ for }\quad s\ge r,
\label{eq:4.5}
\end{align}
and
\begin{align}
(h^{{}^{\!+-}}_{rr})_0
&=(-1)^r
\frac{(1-q^\c)^r\cdot(q^{\ta+\b_2-(2n-r-1)\c};q^\c)_r}
     {(q^\c;q^\c)_r \cdot (q^{\ta-2(n-1)\c};q^\c)_r \cdot(q^{\b_2};q^\c)_r},
\label{eq:4.6}
\\[5pt]
\Big(
{h^{{}^{\!+-}}_{rs}\over h^{{}^{\!+-}}_{rr}}
\Big)_{\!0}
&=q^{-s(r-s)\c}\frac{(q^{\ta+\b_1-(n-1)\c}\X21;q^\c)_{r-s}\cdot(q^{\ta+\b_2-(2n-r-1)\c};q^\c)_s}
                    {(q^{\ta+\b_2-(2n-r-1)\c};q^\c)_r}\nonumber\\
&\hskip 25mm \cdot
\frac{(q^{\b_2};q^\c)_r\cdot (q^{\b_1-\b_2+(n-r-s)\c}\X21;q^\c)_{s}}
                    {(q^{\b_2};q^\c)_s\cdot (q^{\b_1+(n-r)\c}\X21;q^\c)_{r-s}}
{}_r(q^\c;q^\c)_s 
\quad\hbox{ for }\quad r\ge s.
\label{eq:4.7}
\end{align}
\par
Let $\t$ be the operation which exchanges $x_1$, $\b_1$ for $x_2$, $\b_2$. 
\begin{thm}[Quasi-symmetry of second kind {\rm \cite{ak5}}
]
\label{thm:4.2}
Under the action of $\t$, the principal connection matrix 
$G=G(\ta,\b_1,\b_2,x_1,x_2)$ changes as follows:
\begin{align}
\t G(\ta,\b_1,\b_2,x_1,x_2)&=G(\ta,\b_2,\b_1,x_2,x_1)\nonumber\\
                     &=S(x_2/x_1)\ G(\ta,\b_1,\b_2,x_1,x_2)\ {}^t\!S(q^{\b_2-\b_1}x_1/x_2),
\label{eq:4.8}
\end{align}
where we put $S(x):=\big(a_{r,s}(x)\,\d_{r,s-n}\big)_{r,s=0}^n$
and 
\begin{equation*}
a_{r,n-r}(x):=x^{2r(n-r)\c}q^{-r(n-r)\c+r(n-r)(n-2r)\c^2}
                    \cdot{}_n\vt(xq^{-(r-1)\c})_r \cdot{}_n\vt(xq^{-r\c})_r.
\end{equation*}
\end{thm}
\par\noindent
{\bf Proof.} See Theorem 5.2 in \cite{ak5}. 
\hfill\qed
\begin{prop}
\label{prop:4.3}
The principal connection matrix $G$ admits the following Gauss decomposition:
$$G=(H_{\<\d;\x\>})^{-1}\ H_{\<\d;\y\>}$$
where $H_{\<\d;\x\>}$ and $H_{\<\d;\y\>}$ are a lower triangular matrix and an upper 
one respectively defined by
\begin{align*}
H_{\<\d;\x\>}&:=S'(x_2/x_1)\ \t H_{\<\z;\x\>}\ {}^t\!S(x_2/x_1)\\
H_{\<\d;\y\>}&:=S'(x_2/x_1)\ \t H_{\<\z;\y\>}\ {}^t\!S(q^{\b_1-\b_2}x_2/x_1),
\end{align*} 
where $S'(x):=\big(a_{r,s}(xq^{-\b_2-(s-1)\c})\cdot\d_{r,s-n}\big)_{r,s=0}^n$.
\end{prop}
\par\noindent
{\bf Proof.}
From \eqref{eq:4.8}, we have 
$$
\tau G=S(x_2/x_1)\ G\ {}^{t}\!S(q^{\b_2-\b_1}x_1/x_2).
$$
Then, 
\begin{equation}
\label{eq:4.10}
G=\tau\tau G=\tau S(x_2/x_1)\ \tau G\ \tau {}^{t}\!S(q^{\b_2-\b_1}x_1/x_2)
=S(x_1/x_2)\ \tau G\ {}^{t}\!S(q^{\b_1-\b_2}x_2/x_1).
\end{equation}
On the other hand, by \eqref{eq:4.3}, we have 
\begin{align}
(H_{\<\d;\x\>})^{-1}\ H_{\<\d;\y\>}
&=\Big(S'(x_2/x_1)\ \t H_{\<\z;\x\>}\ {}^t\!S(x_2/x_1)\Big)^{\!\!-1}
   S'(x_2/x_1)\ \t H_{\<\z;\y\>}\ {}^t\!S(q^{\b_1-\b_2}x_2/x_1)\nonumber\\
&={}^tS(x_2/x_1)^{-1}\ (\t H_{\<\z;\x\>})^{-1}\ S'(x_2/x_1)^{-1}
  \ S'(x_2/x_1)\ \t H_{\<\z;\y\>}\ {}^t\!S(q^{\b_1-\b_2}x_2/x_1)\nonumber\\
&=S(x_1/x_2)\ \t (H_{\<\z;\x\>})^{-1}\ \t H_{\<\z;\y\>}\ {}^t\!S(q^{\b_1-\b_2}x_2/x_1)\nonumber\\
&=S(x_1/x_2)\ \tau G\ {}^{t}\!S(q^{\b_1-\b_2}x_2/x_1).
\label{eq:4.11}
\end{align}
The proposition now follows from \eqref{eq:4.10} and \eqref{eq:4.11}.
\hfill$\square$\\
\par
From Proposition \ref{prop:4.3}, the elements of the matrices $H_{\<\d;\y\>}
=\big(h^{{}^{\!--}}_{rs}\big)_{r,s=0}^n$ 
and 
$H_{\<\d;\x\>}=\big(h^{{}^{\!-+}}_{rs}\big)_{r,s=0}^n$ 
are written as follows:
\begin{equation}
\label{eq:4.12}
h^{{}^{\!--}}_{rs}
= \frac{
\t h^{{}^{\!+-}}_{n-r,n-s}
\cdot a_{s,n-s}(q^{\b_1-\b_2}x_2/x_1)}{a_{n-r,r}(q^{-\b_2-(r-1)\c}x_2/x_1)}
,\quad
h^{{}^{\!-+}}_{rs}
= 
\frac{\t h^{{}^{\!++}}_{n-r,n-s}
\cdot a_{s,n-s}(x_2/x_1)}{a_{n-r,r}(q^{-\b_2-(r-1)\c}x_2/x_1)}.
\end{equation}

\begin{rmk}
\label{rmk:4.4}
{\rm
In \cite{ak5}, 
it is proved that the element $h^{{}^{\!--}}_{rs}$ is equal to 
$\big(\<\d_{F_r^{n-r}}\>;\<\y_{F_s^{n-s}}\>\big)_\vP$. 
Since the elements of the matrix $G=(g_{rs})_{r,s=0}^n$ are expressed as 
$$
g_{rs}
=\big(\<\x_{F_r^{n-r}}\>;\<\y_{F_s^{n-s}}\>\big)_\vP
=\sum_{i=0}^n\big(\<\x_{F_r^{n-r}}\>;\<\d_{F_i^{n-i}}\>\big)_\vP
        \cdot\big(\<\d_{F_i^{n-i}}\>;\<\y_{F_s^{n-s}}\>\big)_\vP,
$$
we finally have $h^{{}^{\!-+}}_{rs}=\big(\<\d_{F_r^{n-r}}\>;\<\x_{F_s^{n-s}}\>\big)_\vP$. 
}
\end{rmk}
%
\subsection{A remark on Aomoto--Kato case}
\label{subsection:4.2}
\par
In \cite{ak6}, 
Aomoto and Kato have also studied the $q$-difference system \eqref{eq:1.3.9}
from the viewpoint of the Riemann--Hilbert method. 
In this section we will show how to derive Corollary \ref{cor:1.3.8} by using it.
\par
Since the functions $\vP_{n,2}(t)$ and $\vp_s(t)$ are depending on $x_1$ and $x_2$, 
we denote $\vP_{n,2}(t)$ and $\vp_s(t)$ by $\vP(x_1,x_2;t)$ and $\vp_s(x_1,x_2;t)$ respectively.
By the transformation 
$$
\int \vP(x_1,x_2;t)\vp_s(x_1,x_2;t)\varpi
=x_1^{\a_1+\cdots+\a_n}\int \vP(1,x_2/x_1;t)\vp_s(1,x_2/x_1;t)\varpi,
$$
it suffices to consider the $q$-difference system \eqref{eq:1.3.9} in the case $x_1=1$ and $x_2=x$. 
For a vector function ${\Vec y}(x):=(\wt\vp_0,\wt\vp_{1},\ldots,\wt\vp_n)$
where 
$$
\vp_s(t)
=\A
           \left[
                   \prod_{k=1}^n
                   \frac{1}{1-t_k/x}
                   \prod_{k=1}^{n-s}
                   \frac{1-q^{\b_2} t_k/x}{1-t_k}
                   \prod_{1\le i<j\le n}(t_i-q^{-\c}t_j)
           \right], 
$$
the $q$-difference system \eqref{eq:1.3.9} is written as 
$$
{\Vec y}(qx)={\Vec y}(x)K(x).
$$
where $K(x):=K(1,x)$. 
We define the following two matrices
$$
Y_\z:=\bigg(\int_{\<\z_{F^{n-r}_r}\>}\vP(t)\vp_s(t)\varpi\bigg)_{r,s=0}^n,
\quad
Y_\d:=\bigg(\int_{\<\d_{F^{n-r}_r}\>}\vP(t)\vp_s(t)\varpi\bigg)_{r,s=0}^n,
$$
where $\z_{F^{n-r}_r}$ and $\d_{F^{n-r}_r}$ is defined in Section 3, and 
these matrices are solutions of the matrix equation 
\begin{equation}
\label{eq:4.13}
Y(qx)=Y(x)K(x).
\end{equation}
Aomoto and Kato studied in \cite{ak6} 
asymptotic behaviors of the solutions $Y_\z$ and $Y_\d$:
\begin{align*}
Y_\z&\sim V_+(x)\ x^{D_{\!K}^+}C_{\!K}^+ \quad\hbox{at}\quad x\to 0,\\
Y_\d&\sim V_-(x)\ x^{D_{\!K}^-}C_{\!K}^- \quad\hbox{at}\quad x\to\8.
\end{align*}
where $V_+(x)$ and $V_-(x)$ are pseudo-constant diagonal matrices evaluated in \cite{ak5} 
as 
\begin{align*} 
V_+(x)&=\diag[v_r(x)]_{r=0}^n,\\
V_-(x)&=\diag[v^*_r(x)]_{r=0}^n,
\end{align*}
and $D_{\!K}^+$ and $D_{\!K}^-$ are diagonal matrices
\begin{align*} 
D_{\!K}^+&=\diag[r\ta-(n-r)\b_2-r(2n-r-1)\c]_{r=0}^n,\\[3pt]
D_{\!K}^-&=\diag[r\ta+r\b_1-r(r-1)\c]_{r=0}^n.
\end{align*}
The matrices $C_{\!K}^+=(c_{\!K,rs}^+)_{r,s=0}^n$ and $C_{\!K}^-=(c_{\!K,rs}^-)_{r,s=0}^n$ 
are a lower triangular matrix and an upper one respectively.
\begin{lem}[Aomoto--Kato {\rm \cite{ak6}}
]
\label{lem:4.5}
The matrices $(V_+(x))_0$, $(V_-(x))_0$, $(C_{\!K}^+)_0$ and $(C_{\!K}^-)_0$
are evaluated as follows:
\begin{align}
(v_{r})_0
&=q^{-(n-r)\b_2}\cdot 
\frac{(xq^{\b_1};q^\c)_{n-r}}
     {(xq^{\b_1-\b_2};q^\c)_{n-r}},
\nonumber\\[7pt]
(v^*_{r})_0
&=
\frac{(-1)^{r(n-r)}q^{r(n-r)\c}\cdot (xq^{-\b_2-(r-1)\c};q^\c)_r}
     {(xq^{-\b_2-(r-1)\c};q^\c)_r
       \cdot{}_n(xq^{-\b_2-(n-2)\c};q^\c)_r \cdot {}_n(xq^{-\b_2-(n-1)\c};q^\c)_r},\nonumber\\[7pt]
(c_{\!K,rr}^+)_0
&=(-1)^{n-r}q^{(n-r)\b_2-(n-r)(n-r-1)\c/2-r(r-1)\c/2}
\frac{(q^{\c};q^\c)_{n-r}\cdot (q^{\ta+\b_2-(2n-r-1)\c};q^\c)_{r}}
     {(1-q^\c)^{n-r}\cdot (q^{\b_2};q^\c)_{r}\cdot (q^{\ta-2(n-1)\c};q^\c)_{r}}
\nonumber\\[7pt]
\Big({c_{\!K,rs}^+ \over c_{\!K,rr}^+}\Big)_0 
&={c_{\!K,rs}^+ \over c_{\!K,rr}^+}
=q^{(n-r)(s-r)\c}
        \frac{(q^\c;q^\c)_{n-s}}{(q^\c;q^\c)_{r-s}\cdot (q^\c;q^\c)_{n-r}}
        \frac{(q^{\b_2+s\c};q^\c)_{r-s}}{(q^{\ta+\b_2-(2n-r-s-1)\c};q^\c)_{r-s}},\nonumber\\[7pt]
(c_{\!K,rr}^-)_0
&=(-1)^{r(n-r+1)}q^{-n(n-1)\c/2}
\frac{(q^{\c};q^\c)_{r}\cdot (q^{\ta+\b_1-(n+r-1)\c};q^\c)_{n-r}}
     {(1-q^\c)^r\cdot (q^{\b_1};q^\c)_{n-r}\cdot (q^{\ta-2(n-1)\c};q^\c)_{n-r}},\nonumber\\[7pt]
\Big({c_{\!K,rs}^- \over c_{\!K,rr}^-}\Big)_0 
&={c_{\!K,rs}^- \over c_{\!K,rr}^-}\nonumber\\
&=q^{(s-r)[\ta-(n-1)\c]+r(r-s)\c}
\cdot        
        \frac{(q^\c;q^\c)_s}{(q^\c;q^\c)_{s-r}\cdot (q^\c;q^\c)_r}
        \frac{(q^{\b_1+(n-s)\c};q^\c)_{s-r}}{(q^{\ta+\b_1-(r+s-1)\c};q^\c)_{s-r}}.
\label{eq:4.14}
\end{align}
\end{lem}
If we define 
$$
Y'_0:=(C_{\!K}^+)^{-1}\ (V_+(x))^{-1} \ Y_\z,\quad 
Y'_\8:=(C_{\!K}^-)^{-1}\ (V_-(x))^{-1} \ Y_\d,
$$
then the matrices $Y'_0=Y'_0(x)$ and $Y'_\8=Y'_\8(x)$ are fundamental solutions of the equation \eqref{eq:4.13}
and satisfy 
\begin{align*}
Y'_0(x)
&\sim (C_{\!K}^+)^{-1} x^{D_{\!K}^+}C_{\!K}^+ \quad\hbox{at}\quad x\to 0,\\[3pt]
Y'_\8(x)
&\sim (C_{\!K}^-)^{-1} x^{D_{\!K}^-}C_{\!K}^- \quad\hbox{at}\quad x\to \8.\\
\end{align*}
\par
Let $\what P$ be the connection matrix between $Y_\z$ and $Y_\d$ defined by 
$$
\what P:=Y_\z \ Y_\d^{-1},
$$
In \cite{ak5}, 
$\what P$ is expressed in the form of Gauss decomposition as 
\begin{equation}
\label{eq:4.15}
\what P:=H_{\<\z;\y\>}\ (H_{\<\d;\y\>})^{-1}.
\end{equation}
By Aomoto's lemma \eqref{eq:2.5}, we have 
\begin{align*}
K(0)K(x)^{-1}
&=\big(Y'_0(x)Y'_\8(x)^{-1}\big)_0\\
&=\big((C_{\!K}^+)^{-1}\ (V_+(x))^{-1} \ \what P \ V_-(x)\ C_{\!K}^-\big)_0.
\end{align*}
From \eqref{eq:4.15} and $K(0)=(C_{\!K}^+)^{-1} q^{D_{\!K}^+}C_{\!K}^+$, we have 
$$K(x)
=\Big(
(C_{\!K}^-)^{-1}\ (V_-(x))^{-1}\ 
H_{\<\d;\y\>}\ H_{\<\z;\y\>}^{-1}\ 
V_+(x)\ q^{D_{\!K}^+}\ C_{\!K}^+
\Big)_0.
$$
Hence, by using \eqref{eq:4.6}, \eqref{eq:4.7}, \eqref{eq:4.12}, Lemma \ref{lem:4.5} and Lemma \ref{lem:5.2.1} in Section \ref{section:5}, 
we can evaluate $K(x)$ as \eqref{eq:1.3.11} in Corollary \ref{cor:1.3.8}.
%
\section{Main result for $q$-difference equations with parameter shift 
$\tilde{\alpha}\to 
\tilde{\alpha}+1$%
}
\label{section:5}
\par
From \eqref{eq:3.10} and Proposition \ref{prop:4.3}, it follows that 
\begin{equation}
\label{eq:5.0.1}
A(q^{\ta})=\Big((C_{\!A}^-)^{-1}\ 
(H_{\<\d;\y\>})^{-1}\ H_{\<\d;\x\>}\ 
q^{D_{\!A}^+}\ C_{\!A}^+\Big)_{\!0}.
\end{equation}
Since the matrices $C_{\!A}^-$ and $H_{\<\d;\y\>}$ are upper triangular 
and the matrices $C_{\!A}^+$ and $H_{\<\d;\x\>}$ are lower triangular, 
we can decompose the matrix $A(q^{\ta})$ 
as the product of lower and upper triangular matrices in the following form:
$$
\hskip -5mm
\left(
\begin{matrix}
        1&u_{01}&u_{02}&\cdots&u_{0n}\\
         &     1&u_{12}&\cdots&u_{1n}\\
         &     &\ddots&\ddots&\vdots\\
         &     &      &1     &u_{n\!-\!1,n}\\
         &     &      &      &1
\end{matrix}
\right)
\left(
\begin{matrix}
        d_{0}&     &      &      &      \\
         &     d_{1}&     &      &     \\
         &     &\ddots&      &      \\
         &     &      &d_{n\!-\!1}     &       \\
         &     &      &      &d_{n}
\end{matrix}
\right)
\left(
\begin{matrix}
        1&     &      &      &      \\
   l_{10}&     1&     &      &      \\
   \vdots&\ddots      &\ddots&      &      \\
   l_{n\!-\!1,0}
&\cdots&   l_{n\!-\!1,n\!-\!1}     &1     &      \\
   l_{n,0}&\cdots  &l_{n,n\!-\!2}&l_{n,n\!-\!1}      &1
\end{matrix}
\right).
$$
This expression is unique and we denote by $U_{\!A}$ and $L_{\!A}$ 
the above left and right matrices respectively, so that 
\begin{equation}
\label{eq:5.0.2}
A(q^{\ta})=U_{\!A}\ \diag [d_0,\ldots,d_n]\ L_{\!A}.
\end{equation}
%
\begin{thm}
\label{thm:5.0.1}
The elements of $U_{\!A}$, 
$\diag [d_{0},\ldots,d_{n}]$ and $L_{\!A}$ are expressed as follows:
\begin{align*}
u_{rs}&=(-1)^{s-r}q^{(s-r)[\ta-(n-1)\c]-(s-r)(s+r-1)\c/2}\\
        &\hskip 8mm \cdot\frac{(q^\c;q^\c)_s}{(q^\c;q^\c)_{s-r}\cdot (q^\c;q^\c)_r}
        \frac{(q^{\b_1+(n-s)\c};q^\c)_{s-r}}{(q^{\ta+\b_1-2(s-1)\c};q^\c)_{s-r}}
\quad \hbox{ for }\quad  r\le s,\\[5pt]
d_{r}&=q^{\mu_r}\frac{(q^{\ta+\b_1-2(r-1)\c};q^\c)_r\cdot (q^{\ta-2(n-1)\c};q^\c)_{n-r}}
                  {(q^{\ta+\b_1+\b_2-(r-1)\c};q^\c)_r\cdot (q^{\ta+\b_1-(n-1+r)\c};q^\c)_{n-r}},\\[5pt]
&\hskip 8mm 
\mbox{ where } q^{\mu_r}=x_1^{n-r}x_2^r\ q^{[r(r-1)+(n-r)(n-r-1)]\c/2},\\[5pt]
l_{rs}&=(-1)^{r-s}
         (x_2/x_1)^{s-r}q^{-(r-s)(r+s-1)\c/2}\\
        &\hskip 8mm \cdot\frac{(q^\c;q^\c)_{n-s}}{(q^\c;q^\c)_{r-s}\cdot (q^\c;q^\c)_{n-r}}
        \frac{(q^{\b_2+s\c};q^\c)_{r-s}}{(q^{\ta+\b_1-2(r-1)\c};q^\c)_{r-s}}
\quad \mbox{ for }\quad  r\ge s.
\end{align*}
\end{thm}
We prove Theorem \ref{thm:5.0.1} in the following sections. 
Moreover, separating the matrix depending on $x_1$ and $x_2$, we have 
\begin{thm}
\label{thm:5.0.2}
The matrix $A(q^{\ta})$ is expressed as follows:
$$
A(q^{\ta})=\bar{U}_{\!A}\ \bar{D}_{\!A}\ \bar{L}_{\!A}\ \diag [\ x_1^{n-r}x_2^r\ ]_{r=0}^n,
$$
where 
$\bar{U}_{\!A}=(\bar{u}_{rs})$, $s\ge r$ :upper triangular matrix,
\begin{align*}
\bar{u}_{rs}&=(-1)^{s-r}q^{(s-r)[\ta-(n-1)\c]-(s-r)(s+r-1)\c/2}\\
        &\quad\cdot\frac{(q^\c;q^\c)_s}{(q^\c;q^\c)_{s-r}\cdot (q^\c;q^\c)_r}
        \frac{(q^{\b_1+(n-s)\c};q^\c)_{s-r}}{(q^{\ta+\b_1-2(s-1)\c};q^\c)_{s-r}}\quad r\le s,
\end{align*}
\par
\noindent
$\bar{D}_{\!A}=\diag[\bar{d}_{0},\ldots,\bar{d}_{n}]$,
$$
\bar{d}_{r}=q^{[r(r-1)+(n-r)(n-r-1)]\c/2}
             \frac{(q^{\ta+\b_1-2(r-1)\c};q^\c)_r\cdot (q^{\ta-2(n-1)\c};q^\c)_{n-r}}
                  {(q^{\ta+\b_1+\b_2-(r-1)\c};q^\c)_r\cdot (q^{\ta+\b_1-(n-1+r)\c};q^\c)_{n-r}},
$$
\par
\noindent
$\bar{L}_{\!A}=(\bar{l}_{rs})$, $r\ge s$ :lower triangular matrix,
\begin{align*}
\bar{l}_{rs}&=(-1)^{r-s}
        q^{-(r-s)(r+s-1)\c/2}\\
        &\quad\cdot\frac{(q^\c;q^\c)_{n-s}}{(q^\c;q^\c)_{r-s}\cdot (q^\c;q^\c)_{n-r}}
        \frac{(q^{\b_2+s\c};q^\c)_{r-s}}{(q^{\ta+\b_1-2(r-1)\c};q^\c)_{r-s}}\quad r\ge s.
\end{align*}
\end{thm}
\begin{rmk}
\label{rmk:5.0.3}
{\rm
If we compare Theorem \ref{thm:5.0.2} with Corollary \ref{cor:1.3.8}, 
we find that the matrices $A(q^{\ta})$ and $K(x_1,x_2)$ are very similar to each other, 
especially the substitution of $q^{\a}$ into $(x_1/x_2)q^{-\b_2+(n-1)\c}$ transforms 
$\bar{U}_{\!A}$ and $\bar{L}_{\!A}$ into $U_R$ and $L_R$ respectively.
}\end{rmk}
%
\begin{rmk}
\label{rmk:5.0.4}
{\rm
The elements of $U_{\!A}^{-1}=\big(u_{rs}^{*}\big)_{r,s=0}^n$ coincides with 
the value ${c_{\!K,rs}^- /c_{\!K,rr}^-}$ of the coefficient matrix $C_{\!K}^-$ 
(Compare \eqref{eq:4.14} in Lemma \ref{lem:4.5} with (5.4.2) in Theorem \ref{thm:5.4}).
}\end{rmk}
%
\subsection{Diagonal elements and determinant of 
$A ( q^{\tilde{\alpha}})$ 
}
In this section, we first evaluate the diagonal matrix $\diag [d_0,\ldots,d_n]$ 
of $A(q^{\ta})$ in the expression of \eqref{eq:5.0.2}. 
%
\begin{thm}
\label{thm:5.1.1}
The elements $d_r$, $0\le r\le n$, of the diagonal matrix $\diag [d_0,\ldots,d_n]$ are evaluated as follows:
\begin{align*}
d_r   &=q^{\mu_r}\frac{(q^{\ta+\b_1-2(r-1)\c};q^\c)_r\cdot (q^{\ta-2(n-1)\c};q^\c)_{n-r}}
                   {(q^{\ta+\b_1+\b_2-(r-1)\c};q^\c)_r\cdot (q^{\ta+\b_1-(n-1+r)\c};q^\c)_{n-r}},
\quad\quad 0\le r \le n.
\end{align*}
\end{thm}
In particular,
%
\begin{cor}
\label{cor:5.1.2}
The determinant of the matrix $A(q^{\ta})$ is the following:
$$\det A(q^{\ta})
=(x_1x_2)^{n(n+1)/2}q^{(n-1)n(n+1)\c/3}
\prod_{r=0}^n \frac{(q^{\ta-2(n-1)\c};q^\c)_{r}}
                   {(q^{\ta+\b_1+\b_2-(r-1)\c};q^\c)_r}.
$$
\end{cor}
\par\noindent{\bf Proof.} By Theorem \ref{thm:5.1.1}, we have  
$$\det A(q^{\ta})
=q^{\mu_0+\cdots+\mu_n}
\prod_{r=0}^n \frac{(q^{\ta+\b_1-2(r-1)\c};q^\c)_r\cdot (q^{\ta-2(n-1)\c};q^\c)_{n-r}}
                   {(q^{\ta+\b_1+\b_2-(r-1)\c};q^\c)_r\cdot (q^{\ta+\b_1-(n-1+r)\c};q^\c)_{n-r}}.
$$
The result follows from the following identity:
$\displaystyle 
\prod_{r=0}^n \frac{(q^{\ta+\b_1-2(r-1)\c};q^\c)_r}
                     {(q^{\ta+\b_1-(n-1+r)\c};q^\c)_{n-r}}
=1.
$\hfill\qed
\par\noindent{\bf Proof of Theorem \ref{thm:5.1.1}.} 
From \eqref{eq:4.4}, \eqref{eq:4.6} and \eqref{eq:4.12}, we have the explicit forms of $(h^{{}^{\!--}}_{rr})_0$ 
and 
$(h^{{}^{\!-+}}_{rr})_0$ as follows:
\begin{align}
(h^{{}^{\!--}}_{rr})_0&=(-1)^{n-r}
\frac{(1-q^{\c})^{n-r}\cdot (q^{\ta+\b_1-(n+r-1)\c};q^\c)_{n-r}}
     {(q^{\c};q^\c)_{n-r}\cdot (q^{\ta-2(n-1)\c};q^\c)_{n-r}\cdot (q^{\b_1};q^\c)_r}\nonumber\\
&\hskip 38mm \cdot\frac{{}_n(\X21 q^{\b_1-\b_2-(r-1)\c};q^\c)_r \cdot{}_n(\X21 q^{\b_1-\b_2-r\c};q^\c)_r}
           {{}_n(\X21 q^{-\b_2-(n-2)\c};q^\c)_r \cdot{}_n(\X21 q^{-\b_2-(n-1)\c};q^\c)_r},
\label{eq:5.1.1}
\\[5pt]
(h^{{}^{\!-+}}_{rr})_0
&=(-1)^r 
\cdot\frac{(q^{\b_2};q^\c)_r\cdot (\X21  q^{\b_1};q^\c)_r \cdot (\X21 q^{-\b_2-(r-1)\c};q^\c)_r\cdot (q^{\c};q^\c)_r}
           {(\X21 ;q^\c)_r\cdot (\X21 q^{\b_1-\b_2-(r-1)\c};q^\c)_r \cdot (1-q^\c)^r}\nonumber\\
&\hskip 8mm 
 \cdot\frac{(q^{\ta+\b_1-2(r-1)\c};q^\c)_r }{(q^{\ta+\b_1+\b_2-(r-1)\c};q^\c)_r}
\cdot\frac{{}_n(\X21 q^{-(n-r-1)\c};q^\c)_r \cdot{}_n(\X21 q^{-(n-r)\c};q^\c)_r}
           {{}_n(\X21 q^{-\b_2-(n-2)\c};q^\c)_r \cdot{}_n(\X21 q^{-\b_2-(n-1)\c};q^\c)_r}.
\label{eq:5.1.2}
\end{align}
In Propositions \ref{prop:3.1} and \ref{prop:3.2}, 
we derived the following expressions:
\begin{align}
(c_{\!A,rr}^+)_0
&=q^{-[r(r-1)+(n-r)(n-r-1)]\c/2}
\frac
     { (\X21 q^{-(n-r)\c};q^\c)_{r}\cdot (\X21 q^{-(n-r-1)\c};q^\c)_{r}}
     {(q^{\b_1};q^\c)_{n-r}\cdot (\X21 q^{\b_1};q^\c)_r\cdot (q^{\b_2};q^\c)_r \cdot (\X21 q^{\c};q^\c)_r},
\label{eq:5.1.3}
\\[5pt]
(c_{\!A,rr}^-)_0
&=(-1)^n q^{-[r(r-1)+(n-r)(n-r-1)]\c/2}
\frac{(q^{\c};q^{\c})_{n-r}\cdot(q^{\c};q^{\c})_{r}}{(1-q^{\c})^{n}}
\cdot
  \frac{(\X21 q^{-\b_2-(r-1)\c};q^{\c})_{r}}{(\X21 q^{\b_1-\b_2-(r-1)\c};q^{\c})_{r}}
  \nonumber\\[-1pt]
&\hskip 50mm
\cdot
\frac{1}{{}_n(\X21 q^{\b_1-\b_2-(r-1)\c};q^\c)_{r}\cdot {}_n(\X21 q^{\b_1-\b_2-r\c};q^\c)_{r}}.
\label{eq:5.1.4}
\end{align}
Finally, from \eqref{eq:5.0.1} and \eqref{eq:5.0.2}, we have 
\begin{equation}
d_r={(c_{\!A,rr}^-)_0}^{\!\!-1}
      \cdot {(h^{{}^{\!--}}_{rr})_0}^{\!\!-1}
      \cdot (h^{{}^{\!-+}}_{rr})_0 
      \cdot (c_{\!A,rr}^+)_0
      \cdot q^{\mu_r}
\label{eq:5.1.5}
\end{equation}
and the result follows from \eqref{eq:5.1.1}--\eqref{eq:5.1.5}.\hfill $\square$

%
\subsection{$q$-Binomial lemmas}
\par
In this section, we prepare two $q$-binomial lemmas 
that will be useful in the following sections.
\begin{lem}\label{lem:5.2.1}
Let $z$, $y$ and $c$ be arbitrary numbers $\in \Bbb C$. Then
\begin{equation}
\label{eq:5.2.1}
\sum_{j=0}^k
(-1)^{j}z^j c^{3j(j-1)/2}
        \frac{(c;c)_{k}}{(c;c)_{k-j}\cdot(c;c)_{j}}
\cdot
\frac
     {(yc^{1-j};c)_{j}}
     {(zc^{j-1};c)_{j}}
\cdot
\frac
     {(zyc^j;c)_{k-j}}
     {(zc^{2j};c)_{k-j}}=1.
\end{equation}
\end{lem}
\par\noindent{\bf Remark} 
When $c\to 1$, the above formula reduces to the following well-known combinatorial formula:
$$(1-z)^k=\big((1-zy)-(z-zy)\big)^k=\sum_{j=0}^k(-1)^j{k\choose j}(z-zy)^j(1-zy)^{k-j}.$$
\par\noindent{\bf Proof.}
Multiply both sides of \eqref{eq:5.2.1} by $(z;c)_{2k-1}$. Then we have 
\begin{align}
(z;c)_{2k-1}
&=\sum_{j=0}^k
(-1)^{j}z^j c^{3j(j-1)/2}
        \frac{(c;c)_{k}}{(c;c)_{k-j}\cdot(c;c)_{j}}
\cdot
     (yc^{1-j};c)_{j}
\cdot
     (zyc^j;c)_{k-j}\nonumber\\
&\hskip 38mm\cdot
(z;c)_{j-1}
\cdot
(1-zc^{2j-1})
\cdot
(zc^{j+k};c)_{k-j-1}.
\label{eq:5.2.2}
\end{align}
We prove \eqref{eq:5.2.2} instead of \eqref{eq:5.2.1}. 
We denote by $ g(z)$ and $ g_j(z)$ the summation and summand of the RHS of \eqref{eq:5.2.2} respectively:
\begin{align*}
 g_j(z)
&:=
(-1)^{j}z^j c^{3j(j-1)/2}
        \frac{(c;c)_{k}}{(c;c)_{k-j}\cdot(c;c)_{j}}
\cdot
     (yc^{1-j};c)_{j}
\cdot
     (zyc^j;c)_{k-j}\\
&\hskip 33mm\cdot
(z;c)_{j-1}
\cdot
(1-zc^{2j-1})
\cdot
(zc^{j+k};c)_{k-j-1},
\\[-3pt]
 g(z)
&:=\sum_{j=0}^k g_j(z).
\end{align*}
In order to prove \eqref{eq:5.2.2} by the factor theorem,
we show that the polynomial $ g(z)$ of degree $2k-1$ equals $1$ at $z=0$ 
and vanishes at $2k-1$ points $z=c^{1-l}$, $1\le l \le 2k-1$. 
Since $ g(0)=1$ is easy to check, it is enough to show the following:
\begin{align*}
&
2 g(c^{1-l})
=2\sum_{j=0}^l  g_j(c^{1-l})
=\sum_{j=0}^l\Big( g_{j}(c^{1-l})+ g_{l-j}(c^{1-l})\Big)=0,\\
&
2 g(c^{2-l-k})
=2\sum_{j=l}^k  g_j(c^{2-l-k})
=\sum_{j=0}^{k-l}\Big( g_{l+j}(c^{2-l-k})+ g_{k-j}(c^{2-l-k})\Big)=0,
\end{align*}
which are confirmed from the following lemma. 
\hfill\qed
%
\begin{lem}\label{lem:5.2.3}
For $z=c^{1-l}$ or \ $c^{2-l-k}$, $1\le l\le k$, it follows that  
\vskip 2mm
\hskip 8mm
$
\left\{
\begin{array}{cc}
  g_{j}(c^{1-l})+ g_{l-j}(c^{1-l})=0 \quad\hbox{for}\quad 0\le j \le l, & \\[5pt]
         g_{j}(c^{1-l})=0 \quad\quad\hbox{for}\quad l< j \le k & 
\end{array}
\right.
$
\vskip 2mm
\noindent and
\vskip 2mm
\hskip 8mm
$
\left\{
\begin{array}{cc}
g_{l+j}(c^{2-l-k})+ g_{k-j}(c^{2-l-k})=0  \quad\hbox{for}\quad 0\le j \le k-l, & \\[5pt]
         g_{j}(c^{2-l-k})=0 \quad\quad\hbox{for}\quad 0\le j <l. & 
\end{array}
\right.
$
\end{lem}
\par\noindent{\bf Proof.} It is straightforward and left to the reader.
\hfill\qed
%
\begin{lem}\label{lem:5.2.4}
Let $z$ and $c$ be arbitrary numbers $\in \Bbb C$. Then
\begin{equation}
\label{eq:5.2.3}
\sum_{j=0}^k
\frac
     {(-1)^{j}c^{-j(2k-j-1)/2}}
     {(c;c)_{j}\cdot(c;c)_{k-j}\cdot
     (zc^{-2(j-1)};c)_{j}\cdot
      (zc^{-k-i+1};c)_{k-j}}
=0.
\end{equation}
\end{lem}
\par\noindent{\bf Remark} 
When $c\to 1$, the above formula reduces to the following well-known combinatorial formula:
$$\sum_{j=0}^k(-1)^j{k\choose j}=0.$$
\par\noindent{\bf Proof.}
By multiplying both sides of \eqref{eq:5.2.3} by $(zc^{-2k+2};c)_{2k-1}$, we have 
\begin{equation}
\label{eq:5.2.4}
\sum_{j=0}^k
(-1)^{j}c^{-j(2k-j-1)/2}
        \frac{(z;c^{-j+2})_{j-1}\cdot(1-zc^{-2j+1})\cdot(zc^{-2k+2};c)_{k-j-1}}
             {(c;c)_{k-j}\cdot(c;c)_{j}}
=0.
\end{equation}
We prove \eqref{eq:5.2.4} instead of \eqref{eq:5.2.3}. 
We denote by $\wt g(z)$ and $\wt g_j(z)$ the summation and summand of the RHS of \eqref{eq:5.2.2} respectively:
\begin{align*}
\wt g_j(z)
&:=(-1)^{j}c^{-j(2k-j-1)/2}
        \frac{(z;c^{-j+2})_{j-1}\cdot(1-zc^{-2j+1})\cdot(zc^{-2k+2};c)_{k-j-1}}
             {(c;c)_{k-j}\cdot(c;c)_{j}},
\\[-3pt]
\wt g(z)
&:=\sum_{j=0}^k\wt g_j(z).
\end{align*}
In order to prove \eqref{eq:5.2.4} by the factor theorem,
we show that the polynomial $\wt g(z)$ of degree $k-1$ 
vanishes at $k$ points $z=c^{l-1}$, $1\le l \le k$, i.e.,
\begin{equation*}
2\wt g(c^{l-1})
=2\sum_{j=0}^l \wt g_j(c^{l-1})
=\sum_{j=0}^l\Big(\wt g_{j}(c^{l-1})+\wt g_{l-j}(c^{l-1})\Big)=0,
\quad 
1\le l \le k , 
\end{equation*}
which follows from the following lemma. \hfill\qed
%
\begin{lem}\label{lem:5.2.6}
For $z=c^{l-1}$, $1\le l\le k$, it follows that  
\begin{align*} 
&\wt g_{j}(c^{l-1})+ \wt g_{l-j}(c^{l-1})=0 \quad\hbox{for}\quad 0\le j \le l,\\
&\wt g_{j}(c^{l-1})=0 \quad\quad\hbox{for}\quad l< j \le k.
\end{align*}
\end{lem}
\par\noindent{\bf Proof.} 
It is straightforward and left to the reader.
\hfill\qed\\
%

%
\subsection{Evaluation of $L_{\!A}$}
\begin{thm}
\label{thm:5.3}
The elements of the matrix $L_A$ are expressed in a product of binomials as follows:
\begin{align}
l_{rs}&=(-1)^{r-s}
         (\X21)^{s-r}q^{-(r-s)(r+s-1)\c/2}\nonumber\\
        &\cdot\frac{(q^\c;q^\c)_{n-s}}{(q^\c;q^\c)_{r-s}\cdot (q^\c;q^\c)_{n-r}}
        \frac{(q^{\b_2+s\c};q^\c)_{r-s}}{(q^{\ta+\b_1-2(r-1)\c};q^\c)_{r-s}}
\quad \hbox{ for }\quad  r\ge s.
\label{eq:5.3.1}
\end{align}
\end{thm}
\par\noindent{\bf Proof.} From \eqref{eq:4.5} and \eqref{eq:4.12}, we have  
\begin{align}
\Big(
{h^{{}^{\!-+}}_{ri}\over h^{{}^{\!-+}}_{rr}}
\Big)_{\!0}
&=q^{-i(r-i)\c}
\frac
     { (\X21 q^{\c};q^\c)_{i}\cdot (\X21 q^{-(n-r-1)\c};q^\c)_{r}}
     { (\X21 q^{\c};q^\c)_{r}\cdot (\X21 q^{-(n-i-1)\c};q^\c)_{i}}
\cdot
\frac
     {(\X21 q^{-(n-r)\c};q^\c)_{r}}
     {(\X21 q^{-(n-i)\c};q^\c)_{i}}\nonumber\\
&\hskip 8mm\cdot
\frac
     {(\X21  q^{\b_1};q^\c)_i }
     {(\X21  q^{\b_1};q^\c)_r }
\cdot
\frac
     {(q^{\b_1};q^\c)_{n-i} \cdot (q^\c;q^\c)_{n-i} }
     {(q^{\b_1};q^\c)_{n-r} \cdot (q^\c;q^\c)_{n-r} }\nonumber\\
&\hskip 8mm\cdot
\frac
     {(q^{-\ta -\b_1+(n-1)\c}/\X21 ;q^\c)_{r-i}}
     {(q^{-\ta -\b_1+(r+i-1)\c};q^\c)_{r-i}\cdot(q^{(n-r-i)\c}/\X21 ;q^\c)_{r-i}\cdot(q^\c;q^\c)_{r-i}}.
\label{eq:5.3.2}
\end{align}
By Proposition \ref{prop:3.1}, it follows that
\begin{align}
\Big(\frac{c_{\!A,is}^+}{c_{\!A,ii}^+}\Big)_{\!0}
&=q^{-(i-s)(n-i)\c}
\frac{(q^\c;q^\c)_{n-s}}
     {(q^\c;q^\c)_{i-s}(q^\c;q^\c)_{n-i}}
\frac{(q^{\b_2+s\c};q^\c)_{i-s}}
     {(\X21  q^{(i+s-n)\c};q^\c)_{i-s}},
\label{eq:5.3.3}
\\
\Big(\frac{c_{\!A,ii}^+}{c_{\!A,rr}^+}\Big)_{\!0}
&=q^{[-i(i-1)-(n-i)(n-i-1)+r(r-1)+(n-r)(n-r-1)]\c/2}\nonumber\\[-8pt]
&\hskip 8mm\cdot
\frac
     { (\X21 q^{-(n-i)\c};q^\c)_{i}\cdot (\X21 q^{-(n-i-1)\c};q^\c)_{i}}
     { (\X21 q^{-(n-r)\c};q^\c)_{r}\cdot (\X21 q^{-(n-r-1)\c};q^\c)_{r}}\nonumber\\
&\hskip 8mm\cdot
\frac
     {(q^{\b_1};q^\c)_{n-r}\cdot (\X21 q^{\b_1};q^\c)_r\cdot (q^{\b_2};q^\c)_r \cdot (\X21 q^{\c};q^\c)_r}
     {(q^{\b_1};q^\c)_{n-i}\cdot (\X21 q^{\b_1};q^\c)_i\cdot (q^{\b_2};q^\c)_i \cdot (\X21 q^{\c};q^\c)_i},
\label{eq:5.3.4}
\\[5pt]
\Big(\frac{q^{\mu_{i}}}{q^{\mu_{r}}}\Big)_{\!0}
&=(\X21)^{i-r}q^{[i(i-1)+(n-i)(n-i-1)-r(r-1)-(n-r)(n-r-1)]\c/2}.
\label{eq:5.3.5}
\end{align}
Comparing \eqref{eq:5.0.1} with \eqref{eq:5.0.2}, 
the elements $l_{rs}$ of the matrix $L_A$ are expressed as follows:
$$
l_{rs}
=
\sum_{i=s}^r
\Big(
{h^{{}^{\!-+}}_{ri}\over h^{{}^{\!-+}}_{rr}}
\Big)_{\!0}
\Big(\frac{c_{\!A,is}^+}{c_{\!A,ii}^+}\Big)_{\!0}
\Big(\frac{c_{\!A,ii}^+}{c_{\!A,rr}^+}\Big)_{\!0}
\Big(\frac{q^{\mu_{i}}}{q^{\mu_{r}}}\Big)_{\!0}.
$$
Then, from \eqref{eq:5.3.2}--\eqref{eq:5.3.5}, we have  
\begin{align}
l_{rs}
&=
        \frac{(q^\c;q^\c)_{n-s}}{(q^\c;q^\c)_{n-r}}
        (q^{\b_2+s\c};q^\c)_{r-s}\nonumber\\
&\hskip 8mm\cdot
\sum_{i=s}^r
(-1)^{r-i}(\X21)^{i-r}q^{-(r-i)(r+i-1)\c/2-(i-s)(n-i)\c}
        \frac{1}{(q^\c;q^\c)_{r-i}\cdot(q^\c;q^\c)_{i-s}}\nonumber\\
&\hskip 16mm\cdot
\frac
     {(\X21 q^{\ta +\b_1-(n+r-i-2)\c};q^\c)_{r-i}}
     {(q^{\ta +\b_1-2(r-1)\c};q^\c)_{r-i}
       \cdot(\X21 q^{(2i-n+1)\c};q^\c)_{r-i}\cdot(\X21  q^{(i+s-n)\c};q^\c)_{i-s}}.
\label{eq:5.3.6}
\end{align}
We have to show the following identity to prove that \eqref{eq:5.3.6} coincides with \eqref{eq:5.3.1}:
\begin{align}
&\sum_{i=s}^r
(-1)^{r-i}(\X21)^{i-r}q^{-(r-i)(r+i-1)\c/2-(i-s)(n-i)\c}
        \frac{1}{(q^\c;q^\c)_{r-i}\cdot(q^\c;q^\c)_{i-s}}\nonumber\\
&\hskip 8mm\cdot
\frac
     {(\X21 q^{\ta +\b_1-(n+r-i-2)\c};q^\c)_{r-i}}
     {(q^{\ta +\b_1-2(r-1)\c};q^\c)_{r-i}
       \cdot(\X21 q^{(2i-n+1)\c};q^\c)_{r-i}\cdot(\X21  q^{(i+s-n)\c};q^\c)_{i-s}}\nonumber\\[3pt]
&=
(-1)^{r-s}
         (\X21)^{s-r}q^{-(r-s)(r+s-1)\c/2}
        \frac{1}{(q^\c;q^\c)_{r-s}\cdot (q^{\ta+\b_1-2(r-1)\c};q^\c)_{r-s}}.
\label{eq:5.3.7}
\end{align}
Dividing \eqref{eq:5.3.7} by the RHS of \eqref{eq:5.3.7}, \eqref{eq:5.3.7} reads as follows:
\begin{align}
1
&=\sum_{i=s}^r
(-1)^{i-s}(\X21 q^{-(n-i)\c+(i+s-1)\c/2})^{i-s}
        \frac{(q^\c;q^\c)_{r-s}}{(q^\c;q^\c)_{r-i}\cdot(q^\c;q^\c)_{i-s}}\nonumber\\
&\hskip 8mm\cdot
\frac
     {(q^{\ta +\b_1-(r+i-2)\c};q^\c)_{i-s}}
     {(\X21  q^{(i+s-n)\c};q^\c)_{i-s}}
\frac
     {(\X21 q^{\ta +\b_1-(n+r-i-2)\c};q^\c)_{r-i}}
     {(\X21 q^{(2i-n+1)\c};q^\c)_{r-i}}.
\label{eq:5.3.8}
\end{align}
If we put 
$j=i-s$, $k=r-s$, $z=\X21 q^{2s-n+j}$, $y=q^{\ta+\b_1-2s-k+1}$ and $c=q^{\c}$ in Lemma \ref{lem:5.2.1},
then we see 
the RHS of \eqref{eq:5.3.8} is equal to $1$ and this concludes the proof of the theorem.  
\hfill\qed

%
\subsection{Evaluation of $U_{\! A}$}
\begin{thm}
\label{thm:5.4}
The elements of $U_A=\big(u_{rs}\big)_{r,s=0}^n$ 
and $U_A^{-1}=\big(u_{rs}^{*}\big)_{r,s=0}^n$ are expressed in a product of binomials as follows:
\begin{align}
u_{rs}&=(-1)^{s-r}q^{(s-r)[\ta-(n-1)\c]-(s-r)(s+r-1)\c/2}\nonumber\\
&\hskip 16mm \cdot\frac{(q^\c;q^\c)_s}{(q^\c;q^\c)_{s-r}\cdot (q^\c;q^\c)_r}
        \frac{(q^{\b_1+(n-s)\c};q^\c)_{s-r}}{(q^{\ta+\b_1-2(s-1)\c};q^\c)_{s-r}}
\quad s \ge r,
\label{eq:5.4.1}\\
u^*_{rs}&=q^{(s-r)[\ta-(n-1)\c]+r(r-s)\c}\nonumber\\
&\hskip 16mm \cdot        
        \frac{(q^\c;q^\c)_s}{(q^\c;q^\c)_{s-r}\cdot (q^\c;q^\c)_r}
        \frac{(q^{\b_1+(n-s)\c};q^\c)_{s-r}}{(q^{\ta+\b_1-(r+s-1)\c};q^\c)_{s-r}}
\quad s \ge r.
\label{eq:5.4.2}
\end{align}
\end{thm}
\par\noindent{\bf Proof.} We show \eqref{eq:5.4.2} first.
From \eqref{eq:4.12} and \eqref{eq:4.7}, we have 
\begin{align}
\Big(
{h^{{}^{\!--}}_{ri}
\over 
 h^{{}^{\!--}}_{rr}}
\Big)_{\!0}
&=q^{(i-r)(i+r-n)\c}
\frac{{}_n(\X21 q^{\b_1-\b_2-(i-1)\c};q^\c)_i \cdot {}_n(\X21 q^{\b_1-\b_2-i\c};q^\c)_i }
     {{}_n(\X21 q^{\b_1-\b_2-(r-1)\c};q^\c)_r \cdot {}_n(\X21 q^{\b_1-\b_2-r\c};q^\c)_r }\nonumber\\
&\hskip 4mm\cdot
q^{-(n-i)(i-r)\c}
\frac{(q^{\ta+\b_2-(n-1)\c}/\X21 ;q^\c)_{i-r}}
     {(q^{\b_2+r\c}/\X21 ;q^\c)_{i-r}}
\cdot
{}_{n-r}(q^{\b_2-\b_1+(r+i-n)\c}/\X21 ;q^\c)_{n-i}\nonumber\\
&\hskip 4mm\cdot
\frac{(q^{\a+\b_1-(n+r-1)\c};q^\c)_{n-i}}
{{(q^{\a+\b_1-(n+r-1)\c};q^\c)_{n-r}}}
\cdot
\frac{(q^{\b_1};q^\c)_{n-r}}
{{(q^{\b_1};q^\c)_{n-i}}}
\cdot
{}_{n-r}(q^\c;q^\c)_{n-i}\nonumber\\[5pt]
&=q^{(i-r)(i+r-n)\c}
\frac{{}_n(\X21 q^{\b_1-\b_2-(i-1)\c};q^\c)_i \cdot {}_n(\X21 q^{\b_1-\b_2-i\c};q^\c)_i }
     {{}_n(\X21 q^{\b_1-\b_2-(r-1)\c};q^\c)_r \cdot {}_n(\X21 q^{\b_1-\b_2-r\c};q^\c)_r }\nonumber\\
&\hskip 4mm\cdot
q^{(\ta-(n+r-1)\c)(i-r)}
\frac{(\X21 q^{-\ta-\b_2+(n+r-i)\c};q^\c)_{i-r}}
     {(\X21 q^{-\b_2-(i-1)\c};q^\c)_{i-r}}
\cdot
\frac{(\X21 q^{\b_1-\b_2-(i-1)\c};q^\c)_{i-r} }
     {(\X21 q^{\b_1-\b_2+(n-2i+1)\c};q^\c)_{i-r} }\nonumber\\
&\hskip 4mm\cdot
\frac{(q^{\b_1+(n-i)\c};q^\c)_{i-r}}
{{(q^{\a+\b_1-(i+r-1)\c};q^\c)_{i-r}}}
\cdot
\frac{(q^{\c};q^{\c})_{n-r}}
     {(q^{\c};q^{\c})_{i-r}\cdot(q^{\c};q^{\c})_{n-i}}.
\label{eq:5.4.3}
\end{align}
By Proposition \ref{prop:3.2}, we have 
\begin{align}
\Big(
     \frac{c_{\!A,is}^-}{c_{\!A,ii}^-}
\Big)_{\!0}
&=(\X21 q^{-\b_2-i\c})^{s-i}
\frac{(q^{\b_1+(n-s)\c};q^{\c})_{s-i}}
     {(\X21 q^{\b_1-\b_2+(n-s-i)\c};q^{\c})_{s-i}}
\cdot
     \frac{(q^{\c};q^{\c})_s}{(q^{\c};q^{\c})_i\cdot (q^{\c};q^{\c})_{s-i}},
\label{eq:5.4.4}\\[5pt]
\Big(
{{c_{\!A,ii}^-}\over {c_{\!A,rr}^-}}
\Big)_{\!0}
&=q^{-(i-r)(i+r-n)\c}
\frac{(q^{\c};q^{\c})_{n-i}\cdot(q^{\c};q^{\c})_{i}}
     {(q^{\c};q^{\c})_{n-r}\cdot(q^{\c};q^{\c})_{r}}\nonumber\\
&\hskip 16mm\cdot
  \frac{(\X21 q^{-\b_2-(i-1)\c};q^{\c})_{i}}{(\X21 q^{\b_1-\b_2-(i-1)\c};q^{\c})_{i}}
\cdot
  \frac{(\X21 q^{\b_1-\b_2-(r-1)\c};q^{\c})_{r}}{(\X21 q^{-\b_2-(r-1)\c};q^{\c})_{r}}
\nonumber\\
&\hskip 16mm\cdot
\frac{{}_n(\X21 q^{\b_1-\b_2-(r-1)\c};q^\c)_{r}\cdot {}_n(\X21 q^{\b_1-\b_2-r\c};q^\c)_{r}}
     {{}_n(\X21 q^{\b_1-\b_2-(i-1)\c};q^\c)_{i}\cdot {}_n(\X21 q^{\b_1-\b_2-i\c};q^\c)_{i}}.
\label{eq:5.4.5}
\end{align}
Comparing \eqref{eq:5.0.1} with \eqref{eq:5.0.2}, using \eqref{eq:5.4.3}--\eqref{eq:5.4.5},
the elements $u^*_{rs}$ of the matrix $U_A^{-1}$ are expressed as follows:
\begin{align*}
u^*_{rs}
&
=\sum_{i=r}^s
\Big(
{h^{{}^{\!--}}_{ri}
 c_{\!A,is}^-
\over 
 h^{{}^{\!--}}_{rr}
 c_{\!A,rr}^-}
\Big)_{\!0}
=\sum_{i=r}^s
\Big(
{h^{{}^{\!--}}_{ri}
\over 
 h^{{}^{\!--}}_{rr}}
\Big)_{\!0}
\Big(
     \frac{c_{\!A,is}^-}{c_{\!A,ii}^-}
\Big)_{\!0}
\Big(
{{c_{\!A,ii}^-}\over {c_{\!A,rr}^-}}
\Big)_{\!0}
\\
&=
q^{[\ta-(n-1-r)\c](s-r)}
        \frac{(q^\c;q^\c)_s}{(q^\c;q^\c)_{s-r}\cdot (q^\c;q^\c)_r}
        \frac{(q^{\b_1+(n-s)\c};q^\c)_{s-r}}{(q^{\ta+\b_1-(r+s-1)\c};q^\c)_{s-r}}
\\
&
\hskip 8mm
\cdot
\sum_{i=r}^s
(\X21 q^{-\ta-\b_2+(n+r-i-1)\c})^{s-i}
\frac{(q^{\ta+\b_1-(r+s-1)\c};q^\c)_{s-i}}
     {(\X21 q^{\b_1-\b_2+(n-s-i)\c};q^{\c})_{s-i}}\\
&
\hskip 16mm
\cdot
\frac{(\X21 q^{-\ta-\b_2+(n+r-i)\c};q^\c)_{i-r}}
     {(\X21 q^{\b_1-\b_2+(n-2i+1)\c};q^\c)_{i-r}}
\cdot
    {}_{s-r}(q^{\c};q^{\c})_{s-i}.
\end{align*}
Thus the statement \eqref{eq:5.4.2} in Theorem \ref{thm:5.4} follows from the following identity: 
\begin{align}
1
&=
\sum_{i=r}^s
(\X21 q^{-\ta-\b_2+(n+r-i-1)\c})^{s-i}\cdot
\frac{(q^{\ta+\b_1-(r+s-1)\c};q^\c)_{s-i}}
     {(\X21 q^{\b_1-\b_2+(n-s-i)\c};q^{\c})_{s-i}}
\nonumber\\
&
\hskip 8mm\cdot
\frac{(\X21 q^{-\ta-\b_2+(n+r-i)\c};q^\c)_{i-r}}
     {(\X21 q^{\b_1-\b_2+(n-2i+1)\c};q^\c)_{i-r}}
\cdot
    {}_{s-r}(q^{\c};q^{\c})_{s-i}.
\label{eq:5.4.6}
\end{align}
We put 
$i=r+j$, $s=r+k$, $y=q^{-\ta-\b_1+(2r+k-1)\c}$, $z=q^{\b_1-\b_2+(n-2r-2k+1)\c}$ 
and $c=q^{\c}$ in the RHS of \eqref{eq:5.4.6}. Then, we have 
\begin{align}
\hbox{The RHS of \eqref{eq:5.4.6}}
&=
\sum_{j=0}^k (zyc^{k-j-1})^{k-j}
\frac{(y^{-1};c)_{k-j}}{(z c^{k-j-1};c)_{k-j}}
\cdot
\frac{(zy c^{k-j};c)_{j}}{(z c^{2(k-j)};c)_{j}}
\cdot
{}_k(c;c)_j
\nonumber\\
&=
\sum_{j=0}^k 
(-1)^{k-j}
z^{k-j}c^{3(k-j)(k-j-1)/2}
\frac{(yc^{1-(k-j)};c)_{k-j}}{(z c^{k-j-1};c)_{k-j}}
\cdot
\frac{(zy c^{k-j};c)_{j}}{(z c^{2(k-j)};c)_{j}}
\cdot
{}_k(c;c)_j
\nonumber\\
&=
\sum_{j=0}^k 
(-1)^{j}
z^jc^{3j(j-1)/2}
\frac{(yc^{1-j};c)_j}{(z c^{j-1};c)_{j}}
\cdot
\frac{(zy c^{j};c)_{k-j}}{(z c^{2j};c)_{k-j}}
\cdot
{}_k(c;c)_j .
\label{eq:5.4.7}
\end{align}
By Lemma \ref{lem:5.2.1}, we have already known that \eqref{eq:5.4.7} is equal to $1$. 
Hence \eqref{eq:5.4.2} follows.
%

Next we show \eqref{eq:5.4.1} of Theorem \ref{thm:5.4}.
What we want to prove is 
the following:
\begin{align}
&\sum_{j=r}^s\Bigg[(-1)^{j-r}q^{(j-r)[\ta-(n-1)\c]-(j-r)(j+r-1)\c/2}
\frac{(q^\c;q^\c)_j}{(q^\c;q^\c)_{j-r}\cdot (q^\c;q^\c)_r}
        \frac{(q^{\b_1+(n-j)\c};q^\c)_{j-r}}{(q^{\ta+\b_1-2(j-1)\c};q^\c)_{j-r}}\Bigg]\cdot u^*_{js}\nonumber\\
&=
\left\{
\begin{matrix}
1&r=s\\[3pt]0& r\not =s
\end{matrix}
\right.
.
\label{eq:5.4.8}
\end{align}
The LHS of \eqref{eq:5.4.8} is equal to
\begin{align}
&\sum_{j=r}^s
\Bigg[(-1)^{j-r}q^{(j-r)[\ta-(n-1)\c]-(j-r)(j+r-1)\c/2}
\frac{(q^\c;q^\c)_j}{(q^\c;q^\c)_{j-r}\cdot (q^\c;q^\c)_r}
        \frac{(q^{\b_1+(n-j)\c};q^\c)_{j-r}}{(q^{\ta+\b_1-2(j-1)\c};q^\c)_{j-r}}\Bigg]
\nonumber\\
&\hskip 8mm \cdot\Bigg[
q^{(s-j)[\ta-(n-1)\c]+j(j-s)\c}
        \frac{(q^\c;q^\c)_s}{(q^\c;q^\c)_{s-j}\cdot (q^\c;q^\c)_j}
        \frac{(q^{\b_1+(n-s)\c};q^\c)_{s-j}}{(q^{\ta+\b_1-(j+s-1)\c};q^\c)_{s-j}}
\Bigg]\nonumber\\[5pt]
&=\sum_{j=r}^s
(-1)^{j-r}q^{(s-r)[\ta-(n-1)\c]+j(j-s)\c-(j-r)(j+r-1)\c/2}
        \frac{(q^\c;q^\c)_s}{(q^\c;q^\c)_r\cdot (q^\c;q^\c)_{j-r}\cdot (q^\c;q^\c)_{s-j}}
\nonumber\\
&\hskip 16mm 
\cdot 
        \frac{(q^{\b_1+(n-s)\c};q^\c)_{s-r}}
             {(q^{\ta+\b_1-2(j-1)\c};q^\c)_{j-r}\cdot(q^{\ta+\b_1-(j+s-1)\c};q^\c)_{s-j}}.
\label{eq:5.4.9}
\end{align}
If $r=s$, it is easy to see that \eqref{eq:5.4.9} is equal to $1$.
We assume $r\not =s$ and put $k:=s-r$. 
\eqref{eq:5.4.9} is equal to 
\begin{align}
&q^{k[\a-(n+r-1)\c]}
        \frac{(q^\c;q^\c)_{r+k}\cdot(q^{\b_1+(n-s)\c};q^\c)_{k}}{(q^\c;q^\c)_r}
\nonumber\\
&\cdot\Bigg[\sum_{i=0}^k
        \frac{(-1)^i q^{-i(2k-i-1)\c/2}}
             {(q^\c;q^\c)_{i}\cdot (q^\c;q^\c)_{k-i} \cdot 
              (q^{\ta+\b_1-2(i+r-1)\c};q^\c)_{i}\cdot(q^{\ta+\b_1-(i+2r+k-1)\c};q^\c)_{k-i}}
\Bigg].
\label{eq:5.4.10}
\end{align}
If we put 
$z=q^{\a+\b_1-2r\c}$ and $c=q^{\c}$ in Lemma \ref{lem:5.2.3},
then we see \eqref{eq:5.4.10} is equal to $0$.
Therefore the proof is complete.  
\hfill\qed


\appendix
%
\section{The condition P4) of Riemann--Hilbert problem for 
$A(q^{\tilde{\alpha}})$ }
\label{section:6}
\begin{prop}
\label{prop:6.1}
The matrix $A(q^{\ta})$ depends only on 
$x$, $q^{\ta}$, $q^{\b_1}$, $q^{\b_2}$, $q^\c$
and it satisfies the condition {\rm P4)}, i.e., it does not depend on $q$.
\end{prop}
This proposition was suggested to me by Prof.~K.~Aomoto.
Before proving Proposition \ref{prop:6.1}, we prove five lemmas.
\par
For $\chi=(\nu_1,\ldots,\nu_n)\in \Bbb Z^n$,
we denote by $T^\chi$ 
the shift operator $T^\chi f(t):=f(q^{\nu_1}t_1,\ldots,q^{\nu_n}t_n)$ 
for a function $f(t)$ of $t\in (\Bbb C^*)^n$.
The function $b_\chi(t)$ called the {\it $b$-function} is defined by 
\begin{equation}
\label{eq:6.1}
b_\chi(t):
  =\displaystyle{T^\chi\vP(t)\over \vP(t)}
.
\end{equation} 
Let $\nabla_{\!\chi}$ be the {\it covariant $q$-difference operator} defined by 
$$\nabla_{\!\chi}\vp(t):= \vp(t)-b_\chi(t)\!\cdot T^\chi\vp(t)$$
for a rational function $\vp(t)$.
\begin{lem}
\label{lem:6.2}
The following equation holds for any $\chi\in \Bbb Z^n$:
$$\int_{\<\x\>}
  \vP(t)\cdot \nabla_{\!\chi}\vp(t)\ \varpi
=0.$$
\end{lem}
\par\noindent{\bf Proof.}
By definition of the Jackson integral, it follows that 
$$
\int_{\<\x\>}
  \vP_R(t)\vp(t)\ \varpi
=
\int_{\<\x\>}
  T^\chi(\vP_R\ \vp)(t)\ \varpi.
$$
Therefore, from \eqref{eq:6.1}, we have  
$\displaystyle
\int_{\<\x\>}
  \vP(t)\vp(t)\ \varpi
=
\int_{\<\x\>}
  \vP(t)\cdot\big(b_\chi(t)\!\cdot T^\chi\vp(t)\big)\ \varpi.
$
\hfill\qed\\

The following lemma is easily deduced from Lemma \ref{lem:6.2} 
and its proof is left to the reader. 
\begin{lem}
\label{lem:6.3}
The following equation holds for any $\chi\in \Bbb Z^n$:
$$\int_{\<\x\>}
  \vP(t) \cdot \A\big(\nabla_{\!\chi}\vp(t)\big)\ \varpi
=0.$$
\end{lem}
\vskip 4mm
%
We set 
$$D(t):=\prod_{1\le i<j\le n} (t_i-\ t_j).$$
The following useful lemma was proved by Kadell in \cite{kad}: 
\begin{lem}[Kadell's lemma]
\label{lem:6.4}
Let $Q\in \Bbb C$ be an arbitrary number and $J$ be a subset of $\{1,2,\ldots,$\\$n\}$. 
Then 
\begin{equation}
\label{eq:6.2}
\A\Big\{
          \prod_{j\in J} t_j  \!\!\!
          \prod_{1\le i<j\le n} (t_i-Q\ t_j)
  \Big\}
=Q^{e(J)}\cdot
         \frac{(Q;Q)_{|J|}\cdot(Q;Q)_{n-|J|}}{(1-Q)^n}\cdot 
        e_{|J|}(t)\ D(t),
\end{equation}
where $|J|:=\#J$, $e(J):= \#\{(i,j)\,|\,1\le i < j \le n,\, i\not\in J,\, j\in J\}$ 
and $e_k(t)$ is the $k$-th elementary symmetric polynomial in variables $t_1,\ldots,t_n$.
\end{lem}
\par\noindent{\bf Proof.} See \cite{kad}. 
\hfill\qed\\
\par
We define $\A_{i,j}(t)$ by the following:
$$
\A_{i,j}(t):=D(t)\cdot{\cal S}(t_1^2 t_2^2\cdots t_i^2 \cdot t_{i+1}t_{i+2}\cdots t_j),
\quad 0 \le i\le j \le n$$
where $\cal S$ is a symmetric sum such that ${\cal S} g(t)=\sum_{\s\in \frak S_n}\s g(t)$.
\par 
In order to make the notations clear, let us write down some elements of $\A_{i,j}(t)$ explicitly.
\begin{align} 
\A_{0,r}(t)&=e_r(t)\cdot D(t),  \ \ 0\le r \le n,    
\label{eq:6.3}\\
\A_{1,1}(t)&=(t_1^2+t_2^2+\cdots +t_n^2)\cdot D(t),\nonumber\\ 
           &\hskip 1.7mm \vdots\nonumber\\
\A_{n,n}(t)&=t_1^2\cdots t_n^2 \cdot D(t),\nonumber\\
\A_{r,n}(t)&=e_n(t)\cdot\A_{0,r}(t),  \ \ 0\le r \le n.
\label{eq:6.4}
\end{align}

We define $\vp_{r,s}(t)$ by the following:
$$
\vp_{r,s}(t):
=\frac{\A_{r,s}(t)}{\prod_{i=1}^n(1-t_i/x_1)(1-t_i/x_2)}.
$$
\begin{lem}
\label{lem:6.5}
$\vp_s(t)$ is expressed as a linear combination of $\vp_{0,r}(t)$, $0\le r\le n$: 
$$
\vp_s(t)=\sum_{r=0}^n u_{rs}\ \vp_{0,r}(t),
$$
where the coefficient $u_{rs}$ does not depend on $q$ and $q^{\ta}$.
\end{lem}
\par\noindent{\bf Proof.} It is straightforward from Kadell's lemma \eqref{eq:6.2}.\hfill\qed
%
%
\begin{lem}
\label{lem:6.6}
$\vp_{s,n}(t)$, $0\le s\le n$, is cohomologous to a linear combination of 
$\vp_{0,r}(t)$, $0\le r\le n$: 
$$
\int \vP(t)\vp_{s,n}(t)\varpi 
= \sum_{r=0}^n c_{s,r} \int \vP(t)\vp_{0,r}(t)\varpi,
$$
where the coefficients $c_{s,r}, 0\le s\le r\le n$, do not depend on $q$. 
\end{lem}
\par\noindent{\bf Proof.} Taking $\chi=(1,0,\ldots,0)\in \Bbb Z^n$, the $b$-function \eqref{eq:6.1} is 
$$
b_1(t):=\frac{T_1\varPhi(t)}{\varPhi(t)}
=q^{\ta}
\frac{1-t_1 q^{\b_1}/x_1}
     {1-t_1/x_1}
  \cdot
   \frac{1-t_1 q^{\b_2}/x_2}
        {1-t_1/x_2}
     \cdot
      \prod_{j=2}^{n} 
       \frac{t_1-q^{-\c}t_j}
            {t_1-q^{\c-1}t_j}.
$$
We put 
$$
\psi(t)=\frac{(t_2\cdots t_{r})^2
                   \cdot t_{r+1} \cdots t_{s}}
             {\prod_{j=2}^n(1-t_j/x_1)(1-t_j/x_2)}
        \prod_{1\le i< j\le n}(t_i-q^\c t_j).
$$
Then $\nabla_{\!1}\psi(t):=\psi(t)-b_1(t)\cdot T_1\psi(t)$ is written as 
\begin{equation}
\label{eq:6.5}
\nabla_{\!1}\psi(t)=\frac{F(t)}{\prod_{j=1}^n(1-t_j/x_1)(1-t_j/x_2)},
\end{equation}
where
\begin{align*} 
&F(t)= (t_2\cdots t_{r})^2 \cdot 
        t_{r+1} \cdots t_{s} \cdot 
      \bigg[ 
     (1-t_1/x_1)(1-t_1/x_2)\prod_{1\le i< j\le n}(t_i-q^\c t_j)
      \\
&\hskip 40mm -q^{\ta}
        (1-t_1 q^{\b_1}/x_1)
        (1-t_1 q^{\b_2}/x_2)
        \prod_{j=2}^n(t_1-q^{-\c} t_j)
        \prod_{2\le i< j\le n}(t_i-q^\c t_j)
      \bigg].
\end{align*} 
We can expand $\prod_{1\le i< j \le n} (t_i-Q_{ij}t_j)$ as follows:
$$
\prod_{1\le i< j \le n} (t_i-Q_{ij}t_j)
=\sum_{\s\in \frak S_n}Q_\s \cdot t_{\s(1)}^{\ n-1}\ t_{\s(2)}^{\ n-2}\cdots t_{\s(n-1)}
$$
like as
$$
\prod_{1\le i< j \le n} (t_i-t_j)
=\sum_{\s\in \frak S_n}\sgn \s \cdot t_{\s(1)}^{\ n-1}\ t_{\s(2)}^{\ n-2}\cdots t_{\s(n-1)}
=\A( t_{1}^{\ n-1}\ t_{2}^{\ n-2}\cdots t_{n-1}).
$$
Thus monomials which appear in $F(t)$ are the following:
\begin{align}
t_1^2t_2^2\cdots t_{r}^2 \cdot t_{r+1} \cdots t_{s}&\cdot t_{\s(1)}^{\ n-1}\ t_{\s(2)}^{\ n-2}\cdots t_{\s(n-1)},\nonumber\\
t_1 t_2^2\cdots t_{r}^2 \cdot t_{r+1} \cdots t_{s}&\cdot t_{\s(1)}^{\ n-1}\ t_{\s(2)}^{\ n-2}\cdots t_{\s(n-1)},\nonumber\\
\ t_2^2\cdots t_{r}^2 \cdot t_{r+1} \cdots t_{s}&\cdot t_{\s(1)}^{\ n-1}\ t_{\s(2)}^{\ n-2}\cdots t_{\s(n-1)},\quad  \s \in \frak S_n.
\label{eq:6.6}
\end{align}
We now define a ordering $\prec$ for $l=(l_1,\ldots,l_n)$ and $m=(m_1,\ldots,m_n)$ as follows:
$$
l\prec m 
\quad\hbox{if and only if}\quad  
l_1+\cdots+l_i \le m_1+\cdots+m_i 
\quad\hbox{for all}\quad i\ge 1.   
$$
We set $t^l:=t_1^{l_1}\cdots t_n^{l_n}$ and also define $t^l\prec t^m$ by $l\prec m$.
The monomial 
$$
t_1^2t_2^2\cdots t_{r}^2 \cdot t_{r+1} \cdots t_{s}\cdot t_{1}^{\ n-1}\ t_{2}^{\ n-2}\cdots t_{n-1}
$$
has maximum order with respect to $\prec$ in \eqref{eq:6.6}.
And the monomials appearing in 
$$
\A(t_1^2t_2^2\cdots t_{r}^2 \cdot t_{r+1} \cdots t_{s}\cdot t_{1}^{\ n-1}\ t_{2}^{\ n-2}\cdots t_{n-1})
$$
are included in those in $\A_{rs}$.
Thus $\A F(t)$ can be expressed 
as a linear combination of $\A_{r,s}(t)$ and $\A_{i,j}(t)$, $(i,j)\prec(r,s)$ as follows:
\begin{align*}
\A F(t)&={1-q^{\ta+\b_1+\b_2}\over x_1x_2}\A_{r,s}(t) - c'_{r-1,s+1}\ \A_{r-1,s+1}(t)- c'_{r-2,s+2}\ \A_{r-2,s+2}(t)+\cdots\\
       &\hskip 15mm  - c'_{r-1,s}\ \A_{r-1,s}(t) - c'_{r-2,s+1}\ \A_{r-2,s+1}(t)+\cdots\\
       &\hskip 15mm  - c'_{r-1,s-1}\ \A_{r-1,s-1}(t) - c'_{r-2,s+1}\ \A_{r-2,s+1}(t)+\cdots\\
       &={1-q^{\ta+\b_1+\b_2}\over x_1x_2}\A_{r,s}(t)-\!\!\!\!\sum_{(i,j)\prec(r,s)}\!\!\!\! c'_{i,j}\ \A_{i,j}(t).
\end{align*}
From Lemma \ref{lem:6.3} and \eqref{eq:6.5}, 
$\displaystyle
\A\Big\{\frac{F(t)}{\prod_{j=1}^n(1-t_j/x_1)(1-t_j/x_2)}\Big\}
$
is cohomologous to $0$. 
This implies 
that $\vp_{r,s}(t)$ is cohomologous to a linear combination of $\vp_{i,j}(t)$'s of lower order:
\begin{equation}
\label{eq:6.7} 
\int \vP(t)\vp_{r,s}(t)\varpi 
=\!\!\!\!\sum_{(i,j)\prec(r,s)}\!\!\!\! c''_{i,j}
\int \vP(t)\vp_{i,j}(t)\varpi. 
\end{equation}
The statement of Lemma \ref{lem:6.6} follows from a recurrent use of \eqref{eq:6.7} 
and the fact that all $c''_{i,j}$'s do not depend on $q$. 
\hfill\qed\\
\par
\par\noindent{\bf Proof of Proposition \ref{prop:6.1}.}
From \eqref{eq:6.4} and Lemma \ref{lem:6.6}, we have 
\begin{equation}
\label{eq:6.8} 
T_{\ta}\wt\vp_{0,s}=\int \vP(t)e_n(t)\vp_{0,s}(t)\varpi 
                     =\wt\vp_{s,n}
                     =\sum_{r=0}^n c_{rs}\wt\vp_{0,r},
\end{equation}
where the coefficients $c_{rs}$ do not depend on $q$.
We set 
$${\Vec y}(q^{\ta})=(\wt \vp_0,\wt \vp_1,\ldots,\wt \vp_n),\quad
  {\Vec y'}(q^{\ta})=(\wt \vp_{00},\wt \vp_{01},\ldots,\wt \vp_{0n}).$$
By Lemma \ref{lem:6.5}, we have 
$${\Vec y}(q^{\ta})={\Vec y'}(q^{\ta})U,\quad
{\Vec y}(q^{\ta+1})={\Vec y'}(q^{\ta+1})U$$
where $U$ is the matrix $(u_{rs})_{rs=0}^n$ not depending on $q$. 
Then we have 
$${\Vec y'}(q^{\ta+1})={\Vec y'}(q^{\ta})U A(q^{\ta}) U^{-1},$$
and, by \eqref{eq:6.8}, the matrix $U A(q^{\ta}) U^{-1}$ does not depend on $q$. 
Since the matrices $U A(q^{\ta}) U^{-1}$ and $U$ do not depend on $q$, 
the condition P4) holds for $A(q^{\ta})$.\hfill\qed

%
\section{Proof of Proposition \ref{prop:3.1}}
\label{section:7}
We set 
\begin{align}
\vP'(t)&:=
       \prod_{j=1}^n 
             \frac{(q t_j/x_1)_\8}{(t_j q^{\b_1}/x_1)_\8}
             \frac{(q t_j/x_2)_\8}{(t_j q^{\b_2}/x_2)_\8}
       \prod_{1\le i<j\le n}
             \frac{(q^{1-\c}t_j/t_i)_\8}{(q^\c t_j/t_i)_\8}\nonumber\\
&=\frac{\vP(t)}{t_1^{\a_1}\cdots t_n^{\a_n}\prod_{j=0}^n(1-t_j/x_1)(1-t_j/x_2) },
\label{eq:7.1}
\end{align}
\begin{equation}
\label{eq:7.2}
\phi_s(t_1,\ldots,t_n):=
                   \prod_{k=1}^{n-s}
                   (1-q^{\b_2} t_k/x_2)
                   \prod_{k=n-s+1}^n
                   (1-t_k/x_1)
\end{equation}
and
\begin{equation}
\label{eq:7.3} 
D^\c_{(\nu)}(t_1,\ldots,t_\nu):=\prod_{1\le i<j\le \nu}(t_i-q^{-\c}t_j),
\quad \hbox{in particular}\quad D^0_{(n)}(t)=D(t),
\end{equation}
so that $\vP(t)\vp_s(t)={t_1^{\a_1}\cdots t_n^{\a_n}}\vP'(t)\A\{\phi_s(t) D^\c_{(n)}(t)\}$.
For simplicity we abbreviate 
$\x_{F^{n-r}_r}$ by $\x=(\x_1,\ldots,\x_n)$ only in this section. 
By the definition \eqref{eq:1.1.2}, we have 
$$
\int_{\<\x\>} \vP(t)\vp_s(t)\varpi 
=(1-q)^n\sum_{\nu_1=0}^\8\cdots \sum_{\nu_n=0}^\8
\vP(q^{\nu_1}\xi_1,\ldots,q^{\nu_n}\xi_n)\cdot\vp_s(q^{\nu_1}\xi_1,\ldots,q^{\nu_n}\xi_n),
$$
so that 
\begin{equation}
\label{eq:7.4}
\int_{\<\x\>} \vP(t)\vp_s(t)\varpi 
- (1-q)^n\vP(\x)\vp_s(\x)
=O(\x_1^{\ta}\cdots \x_n^{\ta}q^{\ta}), 
\end{equation}
because the factor included in the integrand $\vP(t)\vp_s(t)$ with respect to $\ta$ 
is only 
$(t_1\cdots t_n)^{\ta}.$
Let $c_{\!A,rs}^+$ be the constant not depending on $q^{\ta}$ defined by \eqref{eq:3.5}, i.e., 
$$
c_{\!A,rs}^+
=(1-q)^n\cdot\x_2^{-2\c}\x_3^{-4\c}\cdots \x_n^{-2(n-1)\c}\cdot 
         (\x_1^{n-1}\x_2^{n-2}\cdots \x_{n-1})^{-1}
\cdot\vP'(\x)\cdot\A \{\phi_s(\x)D^\c_{(n)}(\x)\}.
$$
Then $(1-q)^n\vP(\x)\vp_s(\x)$ is written as 
$
(\x_1\cdots \x_n)^{\ta}\cdot 
c_{\!A,rs}^+
.
$
From \eqref{eq:7.4}, the asymptotic behavior of the matrix $Y_\x$ at $\ta\to \8$ is the following:
$$Y_\x\sim (q^{\ta})^{D_{\!A}^+}C_{\!A}^+,$$
where $(q^{\ta})^{D_{\!A}^+}=(\x_1\cdots \x_n)^{\ta}$ and $C_{\!A}^+=\big(c_{\!A,rs}^+\big)_{r,s=0}^n$. \\
\par
Before proving Proposition \ref{prop:3.1}, we show four lemmas and two propositions.
\begin{lem}
\label{lem:7.2}
Let $\nu_i$, $1\le i\le n$, be integers satisfying 
$$\{\nu_1,\nu_2,\ldots,\nu_n\}=\{1,2,\ldots,n\},\quad \nu_1<\nu_2<\cdots <\nu_{n-r} 
\quad\hbox{and}\quad \nu_{n-r+1}<\nu_{n-r+2}<\cdots <\nu_{n}.$$ 
For $\s\in \frak S_n$, 
we assume $$\{\s(\nu_1),\ldots,\s(\nu_{n-r})\}=\{1,2,\ldots,n-r\}
\quad\hbox{and}\quad
\{\s(\nu_{n-r+1}),\ldots,\s(\nu_n)\}=\{n-r+1,\ldots,n\}.$$ 
Then we have
$$
n-r=\s(\nu_1)>\s(\nu_2)>\cdots>\s(\nu_{n-r})=1,
$$
$$
n=\s(\nu_{n-r+1})>\s(\nu_{n-r+2})>\cdots>\s(\nu_n)=n-r+1
$$
if and only if 
$$\s D_{(n)}^\c(\x)\ne 0.$$
\end{lem}
\par\noindent{\bf Proof.}
We assume that there exist $i$ and $j$, $i<j$ such that 
$1\le \s(i)<\s(j)\le n-r$ or $n-r+1\le \s(i)<\s(j)\le n$. 
When the former holds, we define the set $E:=\{(i,j);1\le i<j\le n \ \hbox{and}\  1\le \s(i)<\s(j)\le n-r\}$.
We take $(i,j)\in E$ such that $\s(j)-\s(i)$ is minimum.
There exists a number $k$ such that 
\begin{equation}
\label{eq:7.5}
\s(k)=\s(i)+1.
\end{equation}
Then we have 
\begin{equation}
\label{eq:7.6}
\s(i)<\s(k)\le \s(j).
\end{equation}
We now suppose $k\le i$. Then $k<j$ and $\s(k)< \s(j)$ by \eqref{eq:7.6}, so that  $(k,j)\in E$.
By using \eqref{eq:7.5}, we have 
$$
\s(j)-\s(k)=(\s(j)-\s(i))-1.
$$ 
This contradicts the fact $\s(j)-\s(i)$ is minimum. 
\par
Thus we have $i<k$. By using \eqref{eq:7.6}, we have $(i,k)\in E$. 
Since $\s(j)-\s(i)$ is minimum, it follows that 
\begin{equation}
\label{eq:7.7}
\s(k)-\s(i)\ge \s(j)-\s(i).
\end{equation}
From \eqref{eq:7.5}, \eqref{eq:7.6} and \eqref{eq:7.7} we have $\s(j)=\s(k)=\s(i)+1$. Then  
$$\x_{\s(i)}-q^{-\c}\x_{\s(j)}=\x_{\s(i)}-q^{-\c}\x_{\s(i)+1}=0$$
because of $q^\c\x_l=\x_{l+1}$ for $1\le l< n-r$ by definition. 
Hence $\s D_{(n)}^\c(\x)= 0$. In the case where $n-r+1\le \s(i)<\s(j)\le n$, 
we can prove $\s D_{(n)}^\c(\x)= 0$ in the same way as above.
\par
Conversely we assume $\s D_{(n)}^\c(\x)= 0$. There exist $i$ and $j$, $i<j$ such that 
$\x_{\s(i)}-q^{-\c}\x_{\s(j)}=0$. 
The case $\s(i)=n-r$ is in contradiction with 
$\x_{\s(j)}=q^{\c}\x_{\s(i)}=x_1 q^{(n-r)\c}$ for generic $x_1$ and $x_2$. 
Thus $\s(i)\ne n-r$. Since $q^\c\x_l=\x_{l+1}$ for $l\ne n-r$ by definition, we have 
$\x_{\s(j)}=q^{\c}\x_{\s(i)}=\x_{\s(i)+1}$, so that $\s(j)=\s(i)+1$. 
Therefore 
$1\le \s(i)<\s(j)\le n-r$ or $n-r+1\le \s(i)<\s(j)\le n$
\hfill\qed 
\begin{lem}
\label{lem:7.3}
For $\s\in \frak S_n$, there exists $i$, $n-s+1\le i$, such that $\s(i)=1$ if and only if 
$$\s \phi_s(\x)=0.$$
\end{lem}
\par\noindent{\bf Proof.}
By recalling that 
$$
\s\phi_s(\x)=
                   \prod_{i=1}^{n-s}
                   (1-q^{\b_2} \x_{\s(i)}/x_2)
                   \prod_{i=n-s+1}^n
                   (1-\x_{\s(i)}/x_1)
\quad \hbox{ and }\quad\x_1=x_1,
$$
the proof easily follows.\hfill\qed
\begin{lem}
\label{lem:7.4}
If $r<s$, then $\A\{\phi_s(\x)D^\c(\x)\}=0$. 
\end{lem}
\par\noindent{\bf Proof.}
Suppose that there exists $\s\in \frak S_n $ such that $\s\phi_s(\x)D^\c(\x)\ne 0$. 
By Lemma \ref{lem:7.2}, there exist $\nu_i$'s, $\nu_1<\cdots<\nu_{n-r}$ such that 
$
n-r=\s(\nu_1)>\s(\nu_2)>\cdots>\s(\nu_{n-r})=1.
$
By Lemma \ref{lem:7.3}, we have $\nu_{n-r}\le n-s$. Thus $1\le \nu_1<\cdots<\nu_{n-r}\le n-s$.
Therefore 
$n-r \le n-s$.
\hfill\qed 
\begin{lem}
\label{lem:7.5}
For $r\ge s$, it follows that 
\begin{align}   
&\A\{\phi_s(\x)D_{(n)}^\c(\x)\}\nonumber\\
&=
   \frac{(q^{-\c};q^{-\c})_{n-s}}
        {(1-q^{-\c})^{n-s}}
\cdot 
  \sgn w\cdot  
w\phi_s(\x) 
 \cdot 
  \left.
w  \left\{
\!
  D_{(n-s)}^0(t_1,\ldots,t_{n-s}) \frac{D_{(n)}^\c(t_1,\ldots,t_n)}
        {D_{(n-s)}^\c(t_1,\ldots,t_{n-s})}
\!
  \right\}
  \right|_{t=\x},
\label{eq:7.8}
\end{align}
where $w=w_{rs}\in \frak S_n$ is the permutation defined by 
$$
w(i)=w_{rs}(i)=
\left\{
\begin{array}{ll}
i  &  \hbox{\rm if} \quad 1\le i\le n-r,    \\
i+s  &  \hbox{\rm if} \quad n-r+1\le i\le n-s,   \\
2n-r+1-i   &  \hbox{\rm if} \quad n-s+1\le i\le n.
\end{array}
\right.
$$
\end{lem}
\par\noindent{\bf Proof.}
We denote by $S_{n-s}$ the subgroup of $\frak S_n$ defined by 
$$
S_{n-s}:
=\{\s\in \frak S_n\,;\, \s(i)=i \quad \hbox{\rm for} \quad n-r+1\le i\le n-r+s\}\simeq \frak S_{n-s}.
$$
By definition, the LHS of \eqref{eq:7.8} is 
$$
\A\{\phi_s(\x)D_{(n)}^\c(\x)\}=\sum_{\s\in \frak S_n}\sgn\s\cdot \s\phi_s(\x)\cdot\s D^\c(\x).
$$
By Lemmas \ref{eq:7.2} and \ref{eq:7.3}, it is enough to consider the summation when $\s$ runs over the set 
$S_{n-s}w=\{\s\in \frak S_n\,;\, \s(i)=2n-r+1-i \quad \hbox{\rm for} \quad n-s\le i\le n\}$. 
Then we have 
\begin{align*}
\A\{\phi_s(\x)D_{(n)}^\c(\x)\}
&=\sum_{\s\in S_{n-s}}\sgn(\s w)\cdot \s w \phi_s(\x)\cdot\s w D_{(n)}^\c(\x)
=\sgn w\cdot w \phi_s(\x)\sum_{\s\in S_{n-s}}\sgn\s\cdot\s w D_{(n)}^\c(\x),
\end{align*}
because $\s w \phi_s(\x) =w \phi_s(\x)$ for $\s\in S_{n-s}$.
Since the following also holds for $\s\in S_{n-s}$:
$$
\s w 
     \left\{
     \!
                                          \frac{D_{(n)}^\c(t_1,\ldots,t_n)}
                                               {D_{(n-s)}^\c(t_1,\ldots,t_{n-s})}
     \!
     \right\}
=
   w 
     \left\{
     \!
                                          \frac{D_{(n)}^\c(t_1,\ldots,t_n)}
                                               {D_{(n-s)}^\c(t_1,\ldots,t_{n-s})}
     \!
     \right\},
$$
it follows that 
\begin{align}
&\A\{\phi_s(\x)D^\c(\x)\}\nonumber\\
&=\sgn w\cdot w \phi_s(\x)
\sum_{\s\in S_{n-s}}\sgn\s\cdot
  \s w  \!\!\left.
     \left\{
     \!
            D_{(n-s)}^\c(t_1,\ldots,t_{n-s})
                                          \frac{D_{(n)}^\c(t_1,\ldots,t_n)}
                                               {D_{(n-s)}^\c(t_1,\ldots,t_{n-s})}
     \!
     \right\}
     \right|_{t=\x}\nonumber\\
&=\sgn w\cdot w \phi_s(\x)\cdot
  w \!\! \left.
     \left\{
     \!
                                          \frac{D_{(n)}^\c(t_1,\ldots,t_n)}
                                               {D_{(n-s)}^\c(t_1,\ldots,t_{n-s})}
     \!
     \right\}
     \right|_{t=\x}
\cdot \!\!
\sum_{\s\in S_{n-s}}\!\!\sgn\s\cdot
  \s w D_{(n-s)}^\c(\x_1,\ldots,\x_{n-s}). 
\label{eq:7.9}
\end{align}
From Kadell's lemma \eqref{eq:6.2}, we have 
\begin{equation}
\label{eq:7.10}
\sum_{\s\in S_{n-s}}\sgn\s\cdot
  \s \Big(w D_{(n-s)}^\c(t_1,\ldots,t_{n-s})\Big)
= 
    \frac{(q^{-\c};q^{-\c})_{n-s}}
        {(1-q^{-\c})^{n-s}}
\ w D_{(n-s)}^0(t_1,\ldots,t_{n-s}).
\end{equation}
The result follows from \eqref{eq:7.9} and \eqref{eq:7.10}.
\hfill\qed 
\begin{prop}
\label{prop:7.6}
If $r<s$, then $c_{\!A,rs}^+=0$.
For $r\ge s$, 
$c_{\!A,rs}^+$ is expressed as follows:
\begin{align*}
c_{\!A,rs}^+ &=\x_2^{-2\c}\x_3^{-4\c}\cdots \x_n^{-2(n-1)\c}\cdot 
         (\x_1^{n-1}\x_2^{n-2}\cdots \x_{n-1})^{-1}\cdot \vP'(\x)\\
&\qquad\cdot 
   \frac{(q^{-\c};q^{-\c})_{n-s}}
        {(1-q^{-\c})^{n-s}}
\cdot 
  \sgn w\cdot  
w\phi_s(\x) 
 \cdot 
  \left.
w  \left\{
\!
  D_{(n-s)}^0(t_1,\ldots,t_{n-s}) \frac{D_{(n)}^\c(t_1,\ldots,t_n)}
        {D_{(n-s)}^\c(t_1,\ldots,t_{n-s})}
\!
  \right\}
  \right|_{t=\x}.
\end{align*}
\end{prop}
\par\noindent{\bf Proof.}
Proposition \ref{prop:7.6} follows from Lemmas \ref{lem:7.4} and  \ref{lem:7.5}. 
\hfill\qed
\begin{prop}
\label{prop:7.7}
We have
\begin{align} 
&(t_1^{n-1}t_2^{n-2}\cdots t_{n-1})^{-1}
 \cdot 
  \sgn w\cdot  
w \left\{
  D_{(n-s)}^0(t_1,\ldots,t_{n-s}) \frac{D_{(n)}^\c(t_1,\ldots,t_n)}
        {D_{(n-s)}^\c(t_1,\ldots,t_{n-s})}
  \right\}\nonumber\\
&=q^{-s(2r-s-1)\c/2}\prod_{1\le i<j\le n}(1-d_{ij}{t_j\over t_i}),
\label{eq:7.11}
\end{align}
where
$$
d_{ij}=
\left\{
\begin{array}{ll}
q^{-\c}  &  \hbox{\rm if} \quad 1\le i\le n-r,\quad n-r+1\le j\le n-r+s,   \\
q^{\c}  &  \hbox{\rm if} \quad n-r+1\le i\le n-r+s,\quad i<j\le n,  \\
1  &  \hbox{\rm otherwise}.
\end{array}
\right.
$$
\end{prop}
\par
We denote the RHS of \eqref{eq:7.11} by $D'_{rs}(t)$. \\
\par\noindent{\bf Proof.}
Since $\sgn w =(-1)^{s(2r-s-1)/2}$, we have to show the following identity instead of \eqref{eq:7.11}: 
\begin{equation}
\label{eq:7.12}
w \left\{
\!
  D_{(n-s)}^0(t_1,\ldots,t_{n-s}) \frac{D_{(n)}^\c(t_1,\ldots,t_n)}
         {D_{(n-s)}^\c(t_1,\ldots,t_{n-s})}
\!
  \right\}
=(-q^{-\c})^{s(2r-s-1)/2}
\!\!\!\!
  \prod_{1\le i<j\le n}
                      (t_i-d_{ij}t_j).
\end{equation}
We expand the LHS of \eqref{eq:7.12} without $w$:
\begin{align*}   
&  D_{(n-s)}^0(t_1,\ldots,t_{n-s}) \frac{D_{(n)}^\c(t_1,\ldots,t_n)}
         {D_{(n-s)}^\c(t_1,\ldots,t_{n-s})}\\
&=  \prod_{1\le i<j\le n-s}
\!\!\!\!
                      (t_i-t_j)
\ \cdot\ 
   \prod_{i=1}^{n-s}
   \prod_{j=n-s+1}^{n}
                      (t_i-q^{-\c}t_j)
\ \cdot \!\!\!\!
   \prod_{n-s+1 \le i<j\le n}
\!\!\!\!
                      (t_i-q^{-\c}t_j)\\
&=  \prod_{1\le i<j\le n-r}
\!\!\!\!
                      (t_i-t_j)
\ \cdot\ 
   \prod_{i=1}^{n-r}
   \prod_{j=n-r+1}^{n-s}
                      (t_i-t_j)
\ \cdot \!\!\!\!
   \prod_{n-r+1 \le i<j\le n-s}
\!\!\!\!
                      (t_i-t_j)
\\
&\quad \cdot\ 
   \prod_{i=1}^{n-r}
   \prod_{j=n-s+1}^{n}
                      (t_i-q^{-\c}t_j)
\ \cdot\ 
   \prod_{i=n-r+1}^{n-s}
   \prod_{j=n-s+1}^{n}
                      (t_i-q^{-\c}t_j)
\ \cdot \!\!\!\!
   \prod_{n-s+1 \le i<j\le n}
\!\!\!\!
                      (t_i-q^{-\c}t_j).
\end{align*}
Therefore we can express the LHS of \eqref{eq:7.12} as follows:
\begin{align*} 
&w \left\{
\!
  D_{(n-s)}^0(t_1,\ldots,t_{n-s}) \frac{D_{(n)}^\c(t_1,\ldots,t_n)}
         {D_{(n-s)}^\c(t_1,\ldots,t_{n-s})}
\!
  \right\}\\
&=  \prod_{1\le i<j\le n-r}
\!\!\!\!
                      (t_{w(i)}-t_{w(j)})
\ \cdot\ 
   \prod_{i=1}^{n-r}
   \prod_{j=n-r+1}^{n-s}
                      (t_{w(i)}-t_{w(j)})
\ \cdot \!\!\!\!
   \prod_{n-r+1 \le i<j\le n-s}
\!\!\!\!
                      (t_{w(i)}-t_{w(j)})
\\
&\hskip 5mm  \cdot\ 
   \prod_{i=1}^{n-r}
   \prod_{j=n-s+1}^{n}
                      (t_{w(i)}-q^{-\c}t_{w(j)})
\cdot\!\!\!\!\!\!
   \prod_{i=n-r+1}^{n-s}
   \prod_{j=n-s+1}^{n}
                      (t_{w(i)}-q^{-\c}t_{w(j)})
\ \cdot \!\!\!\!
   \prod_{n-s+1 \le i<j\le n}
\!\!\!\!
                      (t_{w(i)}-q^{-\c}t_{w(j)})\\
&=  \prod_{1\le i<j\le n-r}
\!\!\!\!
                      (t_i-t_j)
\ \cdot\ 
   \prod_{i=1}^{n-r}
   \prod_{j=n-r+s+1}^{n}
                      (t_i-t_j)
\ \cdot \!\!\!\!
   \prod_{n-r+s+1 \le i<j\le n}
\!\!\!\!
                      (t_i-t_j)
\\
&\hskip 5mm  \cdot\ 
   \prod_{i=1}^{n-r}
   \prod_{j=n-r+1}^{n-r+s}
                      (t_i-q^{-\c}t_j)
\cdot
\Bigg[
   \prod_{i=n-r+s+1}^{n}
   \prod_{j=n-r+1}^{n-r+s}
                      (t_i-q^{-\c}t_j)
\ \cdot \!\!\!\!\!\!\!\!
   \prod_{n-r+s\ge i>j\ge n-r+1}
\!\!\!\!\!\!\!\!\!\!\!\!
                      (t_i-q^{-\c}t_j)
\Bigg]
\\
&=  \prod_{1\le i<j\le n-r}
\!\!\!\!
                      (t_i-t_j)
\ \cdot\ 
   \prod_{i=1}^{n-r}
   \prod_{j=n-r+s+1}^{n}
                      (t_i-t_j)
\ \cdot \!\!\!\!
   \prod_{n-r+s+1 \le i<j\le n}
\!\!\!\!
                      (t_i-t_j)
\\
&\hskip 5mm  \cdot\ 
   \prod_{i=1}^{n-r}
   \prod_{j=n-r+1}^{n-r+s}
                      (t_i-q^{-\c}t_j)
\cdot
   \prod_{i=n-r+1}^{n-r+s}
\ \prod_{i<j\le n}
                    (-q^{-\c})(t_i-q^{\c}t_j).
\end{align*}
This coincides with the RHS of \eqref{eq:7.12}.
\hfill\qed\\
\par\noindent{\bf Proof of Proposition \ref{prop:3.1}.}
From Proposition \ref{eq:7.6}, it follows that
\begin{equation}
\label{eq:7.13}
\frac{c_{\!A,rs}^+}{c_{\!A,rr}^+}
=\Big(\frac{c_{\!A,rs}^+}{c_{\!A,rr}^+}\Big)_0
=
\frac{w_{rs}\phi_s(\x)\cdot D'_{rs}(\x)\cdot (q^{-\c};q^{-\c})_{n-s}}
     {w_{rr}\phi_r(\x)\cdot D'_{rr}(\x)\cdot (q^{-\c};q^{-\c})_{n-r}}
\cdot\frac{1}{(1-q^{-\c})^{r-s}}.
\end{equation}
Since, by definition of $w_{rs}$, $w_{rs}\phi_s(\x)$ is 
$$
w_{rs}\phi_s(\x)=
                   \prod_{i=1}^{n-r}
                   (1-q^{\b_2} \x_i/x_2)
                   \prod_{i=n-r+s}^{n}
                   (1-q^{\b_2} \x_i/x_2)
                   \prod_{i=n-r+1}^{n-r+s}
                   (1-\x_i/x_1),
$$
we have 
\begin{equation}
\label{eq:7.14}
\frac{w_{rs}\phi_s(\x)}{w_{rr}\phi_r(\x)}
=\prod_{i=n-r+s+1}^{n}
 \frac{(1-q^{\b_2} \x_i/x_2)}
      {(1-\x_i/x_1)}
= \frac{(q^{\b_2+s\c};q^\c)_{r-s}}
      {(\X21 q^{s\c};q^\c)_{r-s}}.
\end{equation}
From \eqref{eq:7.11} we have 
\begin{equation}
\label{eq:7.15}
\frac{D'_{rs}(\x)}{D'_{rr}(\x)}
=q^{-(r-s)(r-s-1)\c/2}
\prod_{i=1}^{n-r}\prod_{j=n-r+s+1}^{n}
     \frac{1-\x_j/\x_i}{1-q^{-\c}\x_j/\x_i}
\cdot
\prod_{n-r+s+1\le i<j\le n}\frac{1-\x_j/\x_i}{1-q^{\c}\x_j/\x_i}.
\end{equation}
We calculate factors in \eqref{eq:7.15} as follows:
\begin{equation}
\label{eq:7.16}
\prod_{i=1}^{n-r}\prod_{j=n-r+s+1}^{n}
     \frac{1-\x_j/\x_i}{1-q^{-\c}\x_j/\x_i}
=
\prod_{j=n-r+s+1}^{n}
     \frac{1-\x_j/\x_1}{1-q^{-\c}\x_j/\x_{n-r}}
= \frac{(\X21 q^{s\c};q^\c)_{r-s}}
       {(\X21 q^{(s+r-n)\c};q^\c)_{r-s}}
\end{equation}
and
\begin{equation}
\label{eq:7.17}
\prod_{n-r+s+1\le i<j\le n}\frac{1-\x_j/\x_i}{1-q^{\c}\x_j/\x_i}
=\prod_{n-r+s+1\le i\le n-1}\frac{1-\x_{i+1}/\x_i}{1-q^{\c}\x_n/\x_i}
= \frac{(1-q^\c)^{r-s}}
       {(q^\c;q^\c)_{r-s}}.
\end{equation}
Thus, by using \eqref{eq:7.16} and \eqref{eq:7.17} for \eqref{eq:7.15}, we have 
\begin{align}
\frac{D'_{rs}(\x)}{D'_{rr}(\x)}
&=q^{-(r-s)(r-s-1)\c/2}
 \frac{(\X21 q^{s\c};q^\c)_{r-s}}
       {(\X21 q^{(s+r-n)\c};q^\c)_{r-s}}
\frac{(1-q^\c)^{r-s}}
       {(q^\c;q^\c)_{r-s}}\nonumber\\
&= \frac{(\X21 q^{s\c};q^\c)_{r-s}}
        {(\X21 q^{(s+r-n)\c};q^\c)_{r-s}}
   \frac{(1-q^{-\c})^{r-s}}
        {(q^{-\c};q^{-\c})_{r-s}}.
\label{eq:7.18}
\end{align}
Hence, from \eqref{eq:7.13}, \eqref{eq:7.14} and \eqref{eq:7.18}, it follows that   
\begin{align*}
\frac{c_{\!A,rs}^+}{c_{\!A,rr}^+}
&=\Big(\frac{c_{\!A,rs}^+}{c_{\!A,rr}^+}\Big)_0
=
\frac{(q^{\b_2+s\c};q^\c)_{r-s}}
     {(\X21 q^{(s+r-n)\c};q^\c)_{r-s}}
\frac{(q^{-\c};q^{-\c})_{n-s}}
     {(q^{-\c};q^{-\c})_{r-s}\cdot (q^{-\c};q^{-\c})_{n-r}}\\
&=q^{-(r-s)(n-r)\c}
\frac{(q^{\b_2+s\c};q^\c)_{r-s}}
     {(\X21 q^{(s+r-n)\c};q^\c)_{r-s}}
\frac{(q^{\c};q^{\c})_{n-s}}
     {(q^{\c};q^{\c})_{r-s}\cdot (q^{\c};q^\c)_{n-r}}.
\end{align*}
\par
Next we show the former part of Proposition \ref{prop:3.1}. From Propositions \ref{prop:7.6} and \ref{prop:7.7}, we have 
\begin{align}
(c_{\!A,rr}^+)_0
&=\big(\vP'(\x)\big)_0 
\ w_{rr}\phi_r(\x)
\ D'_{rr}(\x)
\ \frac{(q^{-\c};q^{-\c})_{n-r}}
               {(1-q^{-\c})^{n-r}}\nonumber\\
&=\big(\vP'(\x)\big)_0 
\ w_{rr}\phi_r(\x)
        \prod_{1\le i<j\le n}
(1-d'_{ij}\x_j/\x_i)
\cdot
         q^{-r(r-1)\c/2} 
          \frac{(q^{-\c};q^{-\c})_{n-r}}
               {(1-q^{-\c})^{n-r}}\nonumber\\
&=q^{-r(r-1)\c/2-(n-r)(n-r-1)\c/2} 
          \frac{(q^{\c};q^{\c})_{n-r}}
               {(1-q^{\c})^{n-r}}\nonumber\\
&\hskip 2cm\cdot
\big(\vP'(\x)\big)_0 
\ w_{rr}\phi_r(\x)
        \prod_{1\le i<j\le n}
(1-d'_{ij}\x_j/\x_i),
\label{eq:7.19}
\end{align}
where
$$
d'_{ij}=
\left\{
\begin{array}{ll}
1  &  \hbox{\rm if}\quad 1\le i<j\le n-r,  \\
q^{-\c}  &  \hbox{\rm if} \quad 1\le i\le n-r,\quad n-r+1\le j\le n,   \\
q^{\c}  &  \hbox{\rm if} \quad n-r+1\le i<j\le n.
\end{array}
\right.
$$
Since $\big(\vP'(\x)\big)_0$ and $w_{rr}\phi_r(\x)$ in \eqref{eq:7.19} are calculated as follows:
\begin{align*}
\big(\vP'(\x)\big)_0
&=\Bigg(\prod_{j=1}^n 
             \frac{(q \x_j/x_1)_\8}{(\x_j q^{\b_1}/x_1)_\8}
             \frac{(q \x_j/x_2)_\8}{(\x_j q^{\b_2}/x_2)_\8}
       \prod_{1\le i<j\le n}
             \frac{(q^{1-\c}\x_j/\x_i)_\8}{(q^\c \x_j/\x_i)_\8}
  \Bigg)_0\\
&=\prod_{j=1}^n 
             \frac{1}{(1-\x_j q^{\b_1}/x_1)_\8\cdot(1-\x_j q^{\b_2}/x_2)_\8}
       \prod_{1\le i<j\le n}
             \frac{1}{(1-q^\c \x_j/\x_i)_\8}\\
&=\frac{1}
       {(q^{\b_1};q^\c)_{n-r}\cdot (\X21 q^{\b_1};q^\c)_r\cdot 
        (q^{\b_2}/\X21;q^\c)_{n-r} \cdot (q^{\b_2};q^\c)_r }
       \prod_{1\le i<j\le n}
             \frac{1}{(1-q^\c \x_j/\x_i)_\8}
\end{align*} 
and
$$
w_{rr}\phi_r(\x)
=
                   \prod_{i=1}^{n-r}
                   (1-q^{\b_2} \x_i/x_2)
                   \prod_{i=n-r+1}^{n}
                   (1-\x_i/x_1)
=(q^{\b_2}/\X21;q^\c)_{n-r}\cdot (\X21 ;q^\c)_r,
$$
we have 
\begin{align}
(c_{\!A,rr}^+)_0
&=q^{-r(r-1)\c/2-(n-r)(n-r-1)\c/2} 
          \frac{(q^{\c};q^{\c})_{n-r}}
               {(1-q^{\c})^{n-r}}\nonumber\\
&\hskip 5mm\cdot
\frac
     {(\X21 ;q^\c)_r}
     {(q^{\b_1};q^\c)_{n-r}\cdot (\X21 q^{\b_1};q^\c)_r\cdot (q^{\b_2};q^\c)_r}
\prod_{1\le i<j\le n}
\frac
     {1-d'_{ij}\x_j/\x_i}
     {1-q^{\c}\x_j/\x_i}.
\label{eq:7.20}
\end{align}
A factor 
$\prod_{1\le i<j\le n}
\frac
     {1-d'_{ij}\x_j/\x_i}
     {1-q^{\c}\x_j/\x_i}
$ in \eqref{eq:7.20} is evaluated as follows:
\begin{align}
\prod_{1\le i<j\le n}
\frac
     {1-d'_{ij}\x_j/\x_i}
     {1-q^{\c}\x_j/\x_i}
&=
\prod_{1\le i<j\le n-r}
\!\!\!\!\!\!
\frac
     {1-\x_j/\x_i}
     {1-q^{\c}\x_j/\x_i}
\ 
\prod_{i=1}^{n-r}\prod_{j=n-r+1}^{n}
\frac
     {1-q^{-\c}\x_j/\x_i}
     {1-q^{\c}\x_j/\x_i}\nonumber\\
&=
\prod_{1\le i<j\le n-r}
\!\!\!\!\!\!
\frac
     {1-\x_j/\x_i}
     {1-q^{\c}\x_j/\x_i}
\cdot 
\prod_{i=1}^{n-r}\prod_{j=n-r+1}^{n}
\frac
     {1-q^{-\c}\x_j/\x_i}
     {1-\x_j/\x_i}
\cdot 
\prod_{i=1}^{n-r}\prod_{j=n-r+1}^{n}
\frac
     {1-\x_j/\x_i}
     {1-q^{\c}\x_j/\x_i}\nonumber\\
&=
\frac{(1-q^{\c})^{n-r}}
     {(q^{\c};q^{\c})_{n-r}}
\cdot
\frac
     { (\X21 q^{-(n-r)\c};q^\c)_{r}}
     {(\X21 ;q^\c)_r}
\cdot
\frac{(\X21 q^{-(n-r-1)\c};q^\c)_{r}}{(\X21 q^{\c};q^\c)_r}.
\label{eq:7.21}
\end{align}
Therefore the proposition follows from \eqref{eq:7.20} and \eqref{eq:7.21}.\hfill\qed
%
\section{Proof of Proposition \ref{prop:3.2}}
\label{section:8}
In this section, we abbreviate the point $\y_{F^{n-r}_r}\in (\Bbb C^*)^n$ by $\y=(\y_1,\ldots,\y_n)$.
For the cycle $\<\y\>$, we have already given the definition of the regularized Jackson integral
as follows:
\begin{align}
\int_{\<\y\>}
\vP(t)\vp(t)\varpi
:&=(1-q)^n\sum_{(\nu_1,\ldots,\nu_n)\in{\Bbb Z}^n}
\!\Res_{\ \ 
       \substack{ 
         \substack{t_1=\eta_1q^{\nu_1}, \\[1.5pt] \cdots\cdots,}
         \\ 
          t_{n}=\eta_{n}q^{\nu_{n}} 
       }
     }
\vP(t_1,\ldots,t_n)\vp(t_1,\ldots,t_n)
\frac{d t_1}{t_1}\wedge\cdots\wedge\frac{d t_n}{t_n}\nonumber\\
&=(1-q)^n\sum_{\nu_1=0}^\8\cdots \sum_{\nu_n=0}^\8
\!\Res_{\ \ 
       \substack{ 
         \substack{t_1=\eta_1q^{-\nu_1}, \\[1.5pt] \cdots\cdots,}
         \\ 
          t_{n}=\eta_{n}q^{-\nu_{n}} 
       }
     }
\vP(t)\vp(t)
\frac{d t_1}{t_1}\wedge\cdots\wedge\frac{d t_n}{t_n}.
%
\label{eq:8.1}
\end{align}
We define $\vP'(t)$, $\phi_s(t)$ and $D^\c_{(n)}(t)$ as in \eqref{eq:7.1}, \eqref{eq:7.2} and \eqref{eq:7.3} 
in 
Appendix \ref{section:7}.
Since the factor with respect to $\ta$ in the function 
$\vP(t)\vp_s(t)={t_1^{\a_1}\cdots t_n^{\a_n}}\vP'(t)\A\{\phi_s(t) D^\c_{(n)}(t)\}$  
is only $(t_1\cdots t_n)^{\ta}$, by \eqref{eq:8.1}, we have 
\begin{equation}
\label{eq:8.2}
\int_{\<\y\>}
\vP(t)\vp_s(t)\varpi
-(1-q)^n
\Res_{t=\y}
\vP(t)\vp_s(t)=O(\y_1^{\ta}\cdots \y_n^{\ta}q^{-\ta}).
\end{equation}
Let $c_{\!A,rs}^-$ be the constant not depending on $q^{\ta}$ defined by \eqref{eq:3.5}:
\begin{align*}
c_{\!A,rs}^-
&=(1-q)^n\cdot 
  \y_2^{-2\c}\y_3^{-4\c}\cdots \y_n^{-2(n-1)\c}
  (\y_1^{n-1}\y_2^{n-2}\cdots \y_{n-1})^{-1}
\cdot 
\Res_{t=\y}\vP'(t)\cdot\A \{\phi_s(\y)D^\c_{(n)}(\y)\}.
\end{align*}
Then  $(1-q)^n\Res_{t=\y}\vP(t)\vp_s(t)$ is written as 
$
(\y_1\cdots \y_n)^{\ta}\cdot c_{\!A,rs}^-
$.
From \eqref{eq:8.2}, the asymptotic behavior of the matrix $Y_\y$ at $\ta\to -\8$ is the following:
$$Y_\y\sim (q^{\ta})^{D_{\!A}^-}C_{\!A}^-,$$
where $(q^{\ta})^{D_{\!A}^-}=(\y_1\cdots \y_n)^{\ta}$ 
and $C_{\!A}^-=\big(c_{\!A,rs}^-\big)_{r,s=0}^n$. 
\vskip 4mm
Before proving Proposition \ref{prop:3.2}, we show four lemmas and a proposition.
\begin{lem}
\label{lem:8.2}
Let $\nu_i$, $1\le i\le n$, be $n$ integers satisfying 
$$\{\nu_1,\nu_2\ldots,\nu_n\}=\{1,2,\ldots,n\},\quad \nu_1<\nu_2<\cdots <\nu_{n-r}
\quad\hbox{and}\quad \nu_{n-r+1}<\nu_{n-r+2}<\cdots <\nu_{n}.$$ 
For $\s\in \frak S_n$, 
we assume $$\{\s(\nu_1),\ldots,\s(\nu_{n-r})\}=\{1,2,\ldots,n-r\}
\quad\hbox{and}\quad
\{\s(\nu_{n-r+1}),\ldots,\s(\nu_n)\}=\{n-r+1,\ldots,n\}.$$ 
Then we have
$$
\s(\nu_i)=i\quad \hbox{for}\quad 1\le i\le n
$$
if and only if 
$$\s D_{(n)}^\c(\y)\ne 0.$$
\end{lem}
\par\noindent{\bf Proof.} We can prove this lemma in same way as Lemma \ref{lem:7.2}.\hfill\qed
\begin{lem}
\label{lem:8.3}
For $\s\in \frak S_n$, there exists $i$, $i\le n-s$, such that $\s(i)=n-r+1$ if and only if 
$$\s \phi_s(\y)=0.$$
\end{lem}
\par\noindent{\bf Proof.}
From 
$$
\s\phi_s(\y)=
                   \prod_{i=1}^{n-s}
                   (1-q^{\b_2} \y_{\s(i)}/x_2)
                   \prod_{i=n-s+1}^n
                   (1-\y_{\s(i)}/x_1),
$$
the lemma follows because $\y_{n-r+1}=x_2 q^{-\b_2}$.\hfill\qed 
\begin{lem}
\label{lem:8.4}
If $r>s$, then $\A\{\phi_s(\y)D_{(n)}^\c(\y)\}=0$. 
\end{lem}
\par\noindent{\bf Proof.}
Suppose that there exists $\s\in \frak S_n $ such that $\s\{\phi_s(\y)D_{(n)}^\c(\y)\}\ne 0$. 
By Lemmas \ref{lem:8.2} and \ref{lem:8.3}, we have 
$$
\s^{-1}(n-r+1)<\s^{-1}(n-r+1)<\cdots<\s^{-1}(n)\le n
$$
and
$$
n-s+1\le\s^{-1}(n-r+1). 
$$
Therefore $r \le s$. \hfill\qed 
\begin{lem}
\label{lem:8.5}
For $r\le s$, 
it follows that 
$$
\A\{\phi_s(\y)D_{(n)}^\c(\y)\}
=
   \frac{(q^{-\c};q^{-\c})_{s}}
        {(1-q^{-\c})^{s}}
\cdot 
\phi_s(\y)
 \cdot
  \left.
  \left\{
   \frac{D_{(n)}^\c(t_1,\ldots,t_n)}
        {D_{(s)}^\c(t_{n-s+1},\ldots,t_{n})}
 D_{(s)}^0(t_{n-s+1},\ldots,t_{n})
  \right\}
  \right|_{t=\y}
.
$$
\end{lem}
\par\noindent{\bf Proof.}
We can prove Lemma \ref{lem:8.5} in the same way as Lemma \ref{lem:7.5}.
\hfill\qed
\begin{prop}
\label{prop:8.6}
If $r>s$, then $c_{\!A,rs}^-=0$. 
For $r\le s$, 
$c_{\!A,rs}^-$ is expressed as 
\begin{align*}
c_{\!A,rs}^-
&=(1-q)^n
\cdot
\y_2^{-2\c}\y_3^{-4\c}\cdots \y_n^{-2(n-1)\c}
\cdot 
\Res_{t=\y}\vP'(t) \cdot \phi_s(\y) 
\cdot
         (\y_1^{n-1}\y_2^{n-2}\cdots \y_{n-1})^{-1}\\
&
\hskip 8mm \cdot
   \frac{(q^{-\c};q^{-\c})_{s}}
        {(1-q^{-\c})^{s}}
\cdot
  \left.
   \left\{
   \frac{D_{(n)}^\c(t_1,\ldots,t_n)}
        {D_{(s)}^\c(t_{n-s+1},\ldots,t_{n})}
 D_{(s)}^0(t_{n-s+1},\ldots,t_{n})
   \right\}
  \right|_{t=\y}.
\end{align*}
\end{prop}
\par\noindent{\bf Proof.}
Proposition \ref{prop:8.6} follows from Lemmas \ref{lem:8.4} and \ref{lem:8.5}. 
\hfill\qed\\
\par\noindent{\bf Proof of Proposition \ref{prop:3.2}.}
We set 
$$
D^*_s(t):=
           \prod_{i=1}^{n-s}
           \prod_{j=i+1}^{n}
                           (1-q^{-\c}t_j/t_i)
           \prod_{n-s+1\le i<j\le n}
                           (1-t_j/t_i),
$$
which is equal to 
$$
         (t_1^{n-1}t_2^{n-2}\cdots t_{n-1})^{-1}\cdot 
   \frac{D_{(n)}^\c(t_1,\ldots,t_n)}
        {D_{(s)}^\c(t_{n-s+1},\ldots,t_{n})}
 D_{(s)}^0(t_{n-s+1},\ldots,t_{n}).
$$
In order to evaluate the constant $(c_{\!A,rs}^-/c_{\!A,rr}^-)_0$ explicitly, we calculate 
$D^*_s(\y)/D^*_r(\y)$ and 
$\phi_s(\y)/\phi_r(\y)$
first:
\begin{align}
\frac{D^*_s(\y)}{D^*_r(\y)}
&=
\prod_{n-s+1\le i<j\le n-r}
\frac{(1-\y_j/\y_i)}{(1-q^{-\c}\y_j/\y_i)}
\cdot
           \prod_{i=n-s+1}^{n-r}
           \prod_{j=n-r+1}^{n}
\frac{(1-\y_j/\y_i)}{(1-q^{-\c}\y_j/\y_i)}\nonumber\\
&=           \prod_{i=n-s+1}^{n-r}
\frac{(1-\y_{i+1}/\y_i)}{(1-q^{-\c}\y_{n-r}/\y_i)}
\cdot
           \prod_{i=n-s+1}^{n-r}
\frac{(1-\y_{n-r+1}/\y_i)}{(1-q^{-\c}\y_n/\y_i)}\nonumber\\
&=\frac{(1-q^{-\c})^{s-r}}{(q^{-\c};q^{-\c})_{s-r}}
\cdot
  \frac{(\X21 q^{\b_1-\b_2+(n-s)\c};q^{\c})_{s-r}}
       {(\X21 q^{\b_1-\b_2+(n-s-r)\c};q^{\c})_{s-r}},
\label{eq:8.3}\\[5pt]
     \frac{\phi_s(\y)}{\phi_r(\y)}
&=
           \prod_{i=n-s+1}^{n-r}
\frac{(1-\y_i/x_1)}{(1-q^{\b_2}\y_i/x_2)}
=(\X21  q^{-\b_2})^{s-r}
           \prod_{i=n-s+1}^{n-r}
\frac{(1-x_1/\y_i)}{(1-q^{-\b_2}x_2/\y_i)}\nonumber\\
&=(\X21  q^{-\b_2})^{s-r}
\cdot
  \frac{(q^{\b_1+(n-s)\c};q^{\c})_{s-r}}
       {(\X21 q^{\b_1-\b_2+(n-s)\c};q^{\c})_{s-r}}.
\label{eq:8.4}
\end{align}
From Proposition \ref{prop:8.6}, by using \eqref{eq:8.3} and \eqref{eq:8.4}, it follows that
\begin{align*}
\Big(
     \frac{c_{\!A,rs}^-}{c_{\!A,rr}^-}
\Big)_0
&=
 \frac{c_{\!A,rs}^-}{c_{\!A,rr}^-}
=
     \frac{\phi_s(\y)}{\phi_r(\y)}
\cdot
     \frac{D^*_s(\y)}{D^*_r(\y)}
\cdot
     \frac{(q^{-\c};q^{-\c})_s}{(q^{-\c};q^{-\c})_r\cdot (1-q^{-\c})_r}\\
&=
(\X21  q^{-\b_2})^{s-r}
  \frac{(q^{\b_1+(n-s)\c};q^{\c})_{s-r}}
       {(\X21 q^{\b_1-\b_2+(n-s-r)\c};q^{\c})_{s-r}}
\!
\cdot
\!
     \frac{(q^{-\c};q^{-\c})_s}{(q^{-\c};q^{-\c})_r\cdot (q^{-\c};q^{-\c})_{s-r}}\\
&=
(\X21  q^{-\b_2-r\c})^{s-r}
  \frac{(q^{\b_1+(n-s)\c};q^{\c})_{s-r}}
       {(\X21 q^{\b_1-\b_2+(n-s-r)\c};q^{\c})_{s-r}}
\cdot
     \frac{(q^{\c};q^{\c})_s}{(q^{\c};q^{\c})_r\cdot (q^{\c};q^{\c})_{s-r}}.
\end{align*}
\par
Next we show the former part of Proposition 3.2. By Proposition \ref{prop:8.6}, we have 
\begin{equation}
\label{eq:8.5}
(c_{\!A,rr}^-)_0=\Big(\Res_{t=\y}\vP'(t) \cdot \phi_r(\y)\Big)_0\cdot
D^*_r(\y)
\cdot 
   \frac{(q^{-\c};q^{-\c})_{r}}
        {(1-q^{-\c})^{r}}.
\end{equation}
An explicit form of $\Res_{t=\y}\vP'(t) \cdot \phi_r(\y)$ is expressed as follows: 
\begin{align*}
\Res_{t=\y}\vP'(t)\phi_r(t)
&=
             \frac{(q \y_1/x_1)_\8}{-(q)_\8}
       \prod_{i=2}^{n-r} 
             \frac{(q \y_i/x_1)_\8}{(\y_i q^{\b_1}/x_1)_\8}
       \prod_{i=n-r+1}^n
             \frac{(\y_i/x_1)_\8}{(\y_i q^{\b_1}/x_1)_\8}\\
&\hskip 8mm
\cdot
           \frac{(q \y_{n-r+1}/x_2)_\8}{-(q)_\8}
            \prod_{i=1}^{n-r}
             \frac{(q \y_i/x_2)_\8}{(\y_i q^{1+\b_2}/x_2)_\8}
            \prod_{i=n-r+2}^{n}
             \frac{(q t_i/x_2)_\8}{(t_i q^{\b_2}/x_2)_\8}\\
&\hskip 8mm
\cdot
             \frac{(q^{1-\c}\y_{n-r+1}/\y_{n-r})_\8}{(q^\c \y_{n-r+1}/\y_{n-r})_\8}
       \prod_{i\ne n-r}
             \frac{(q^{1-\c}\y_{i+1}/\y_i)_\8}{-(q)_\8}
       \prod_{1\le i<j\le n-1}
             \frac{(q^{1-\c}\y_{j+1}/\y_i)_\8}{(q^\c \y_{j+1}/\y_i)_\8},
\end{align*}
so that
\begin{align}
&\Big(\Res_{t=\y}\vP'(t)\phi_r(t)\Big)_0\nonumber\\
&=(-1)^n
       \prod_{i=2}^{n-r} 
             \frac{1}{1-\y_i q^{\b_1}/x_1}
       \prod_{i=n-r+1}^n
             \frac{1-\y_i/x_1}{1-\y_i q^{\b_1}/x_1}
            \prod_{i=n-r+2}^{n}
             \frac{1}{1-\y_i q^{\b_2}/x_2}\nonumber\\
&
\qquad\cdot 
\!\!\!\!\!\!\!\!\!
           \prod_{1\le i<j\le n-r-1}
\!\!\!
                     \frac{1}{(1-q^{\c}t_{j+1}/t_i)}
           \prod_{i=1}^{n-r}
           \prod_{j=n-r+1}^{n}
                     \frac{1}{(1-q^{\c}t_j/t_i)}
\!\!\!
           \prod_{n-r+1\le i<j\le n-1}
\!\!\!
                     \frac{1}{(1-q^{\c}t_{j+1}/t_i)}\nonumber\\
&=(-1)^n
       \prod_{i=2}^{n-r} 
             \frac{1}{1-\y_i q^{\b_1}}
       \prod_{i=n-r+1}^n
             \frac{1-\y_i}{1-\y_i q^{\b_1}}
            \prod_{i=n-r+2}^{n}
             \frac{1}{1-\y_i q^{\b_2}/x}\nonumber\\
&
\qquad\cdot 
           \prod_{1\le i<j\le n-r-1}
                     \frac{1}{(1-\y_{j}/\y_i)}
           \prod_{i=1}^{n-r}
           \prod_{j=n-r+1}^{n}
                     \frac{1}{(1-q^{\c}\y_j/\y_i)}
           \prod_{n-r+1\le i<j\le n-1}
                     \frac{1}{(1-\y_{j}/\y_i)}.
\label{eq:8.6}
\end{align}
We calculate $D^*_r(\y)$ as 
\begin{equation}
\label{eq:8.7}
D^*_r(\y)=
           \prod_{1\le i<j\le n-r}
\!\!
                           (1-q^{-\c}\y_j/\y_i)
           \prod_{i=1}^{n-r}
           \prod_{j=n-r+1}^{n}
                           (1-q^{-\c}\y_j/\y_i)
           \prod_{n-r+1\le i<j\le n}
                           (1-\y_j/\y_i).
\end{equation}
From \eqref{eq:8.5}, \eqref{eq:8.6} and \eqref{eq:8.7}, we have 
\begin{align}
(c_{\!A,rr}^-)_0
&=(-1)^n
       \prod_{i=2}^{n-r} 
             \frac{1}{1-\y_i q^{\b_1}}
       \prod_{i=n-r+1}^n
             \frac{1-\y_i}{1-\y_i q^{\b_1}}
            \prod_{i=n-r+2}^{n}
             \frac{1}{1-\y_i q^{\b_2}/x}\nonumber\\
&
\qquad\cdot 
           \prod_{1\le i<j\le n-r-1}
                     \frac{1}{(1-\y_{j}/\y_i)}
           \prod_{i=1}^{n-r}
           \prod_{j=n-r+1}^{n}
                     \frac{1}{(1-q^{\c}\y_j/\y_i)}
           \prod_{n-r+1\le i<j\le n-1}
                     \frac{1}{(1-\y_{j}/\y_i)}\nonumber\\
&\qquad\cdot
           \prod_{1\le i<j\le n-r}
\!\!
                           (1-q^{-\c}\y_j/\y_i)
           \prod_{i=1}^{n-r}
           \prod_{j=n-r+1}^{n}
                           (1-q^{-\c}\y_j/\y_i)
           \prod_{n-r+1\le i<j\le n}
                           (1-\y_j/\y_i)\nonumber\\
&\qquad\cdot 
   \frac{(q^{-\c};q^{-\c})_{r}}
        {(1-q^{-\c})^{r}}.
\label{eq:8.8}
\end{align}
We calculate the factors appearing in \eqref{eq:8.8} as follows:
\begin{align}
&       \prod_{i=2}^{n-r} 
             \frac{1}{1-\y_i q^{\b_1}/x_1}
       \prod_{i=n-r+1}^n
             \frac{1-\y_i/x_1}{1-\y_i q^{\b_1}/x_1}
            \prod_{i=n-r+2}^{n}
             \frac{1}{1-\y_i q^{\b_2}/x_2}\nonumber\\
&=\frac{1}{(q^{-\c};q^{-\c})_{n-r-1}}
\cdot
  \frac{(\X21 q^{-\b_2-(r-1)\c};q^{\c})_{r}}{(\X21 q^{\b_1-\b_2-(r-1)\c};q^{\c})_{r}}
\cdot
  \frac{1}{(q^{-\c};q^{-\c})_{r-1}},
\label{eq:8.9}
\end{align}
\begin{align}
&           \prod_{1\le i<j\le n-r-1}
                     \frac{1}{(1-\y_{j}/\y_i)}
           \prod_{1\le i<j\le n-r}
                           (1-q^{-\c}\y_j/\y_i)\nonumber\\
&=
           \prod_{1\le i<j\le n-r-1}
                     \frac{1}{(1-\y_{j}/\y_i)}
           \prod_{1\le i<j\le n-r-1}
                           (1-\y_{j+1}/\y_i)
           \prod_{i=1}^{n-r-1}
                           (1-q^{-\c}\y_{n-r}/\y_i)\nonumber\\
&=
           \prod_{i=1}^{n-r-2}
                     \frac{(1-\y_{n-r}/\y_i)}{(1-\y_{i+1}/\y_i)}
           \prod_{i=1}^{n-r-1}
                           (1-q^{-\c}\y_{n-r}/\y_i)\nonumber\\
&=\frac{(q^{-\c};q^{-\c})_{n-r-1}}{(1-q^{-\c})^{n-r-1}}\cdot
(q^{-2\c};q^{-\c})_{n-r-1}
=(q^{-\c};q^{-\c})_{n-r-1}\cdot
\frac{(q^{-\c};q^{-\c})_{n-r}}{(1-q^{-\c})^{n-r}},
\label{eq:8.10}
\end{align}
\begin{equation}
\label{eq:8.11}
           \prod_{n-r+1\le i<j\le n-1}
                     \frac{1}{(1-\y_{j}/\y_i)}
           \prod_{n-r+1\le i<j\le n}
                           (1-\y_j/\y_i)
=
           \prod_{i=n-r+1} ^{n-1}
                           (1-\y_n/\y_i)
=(q^{-\c};q^{-\c})_{r-1}
\end{equation}
and 
\begin{align}
           \prod_{i=1}^{n-r}
           \prod_{j=n-r+1}^{n}
                     \frac{(1-q^{-\c}\y_j/\y_i)}{(1-q^{\c}\y_j/\y_i)}
&=
           \prod_{i=1}^{n-r}
                     \frac{(1-\y_n/\y_i)(1-q^{-\c}\y_n/\y_i)}
                          {(1-q^{\c}\y_{n-r+1}/\y_i)(1-q^{\c}\y_{n-r+1}/\y_i)}\nonumber\\
&=
                     \frac{(\X21 q^{\b_1-\b_2-(r-1)\c};q^\c)_{n-r}\cdot (\X21 q^{\b_1-\b_2-r\c};q^\c)_{n-r}}
                          {(\X21 q^{\b_1-\b_2};q^\c)_{n-r}\cdot (\X21 q^{\b_1-\b_2+\c};q^\c)_{n-r}}\nonumber\\
&=\frac{1}{{}_n(\X21 q^{\b_1-\b_2-(r-1)\c};q^\c)_{r}\cdot {}_n(\X21 q^{\b_1-\b_2-r\c};q^\c)_{r}}.
\label{eq:8.12}
\end{align}
From \eqref{eq:8.8}--\eqref{eq:8.12}, it finally follows that 
\begin{align*}
\big(c_{\!A,rr}^-\big)_0
&=(-1)^{n}
\frac{1}{(q^{-\c};q^{-\c})_{n-r-1}}
\cdot
  \frac{(\X21 q^{-\b_2-(r-1)\c};q^{\c})_{r}}{(\X21 q^{\b_1-\b_2-(r-1)\c};q^{\c})_{r}}
\cdot
  \frac{1}{(q^{-\c};q^{-\c})_{r-1}}\\
&\hskip 8mm\cdot
(q^{-\c};q^{-\c})_{n-r-1}
\cdot
\frac{(q^{-\c};q^{-\c})_{n-r}}{(1-q^{-\c})^{n-r}}
\cdot
(q^{-\c};q^{-\c})_{r-1}\\
&\hskip 8mm\cdot
\frac{1}{{}_n(\X21 q^{\b_1-\b_2-(r-1)\c};q^\c)_{r}\cdot {}_n(\X21 q^{\b_1-\b_2-r\c};q^\c)_{r}}
\cdot
\frac{(q^{-\c};q^{-\c})_{r}}{(1-q^{-\c})^{r}}\\
&=(-1)^n\frac{(q^{-\c};q^{-\c})_{n-r}\cdot(q^{-\c};q^{-\c})_{r}}{(1-q^{-\c})^{n}}
\cdot
  \frac{(\X21 q^{-\b_2-(r-1)\c};q^{\c})_{r}}{(\X21 q^{\b_1-\b_2-(r-1)\c};q^{\c})_{r}}\\
&\hskip 8mm\cdot
\frac{1}{{}_n(\X21 q^{\b_1-\b_2-(r-1)\c};q^\c)_{r}\cdot {}_n(\X21 q^{\b_1-\b_2-r\c};q^\c)_{r}}.
\end{align*}
\hfill\qed
\section{Proof of Theorems \ref{thm:1.3.3} and \ref{thm:1.3.4}} 
\label{section:9}
In Section \ref{section:1}, we write
\begin{equation}
\label{eq:9.0}
(\vp_0(t),\vp_1(t),\ldots,\vp_n(t))
R_{2,1}(\textstyle{x_2\over x_1}\displaystyle)=
(\psi_n(t),\psi_{n-1}(t),\ldots,\psi_0(t)),
\end{equation}
so that
\begin{equation}
\label{eq:9.1}
(\wt\vp_0(\x),\wt\vp_1(\x),\ldots,\wt\vp_n(\x))
R_{2,1}(\textstyle{x_2\over x_1}\displaystyle)=
(\wt\psi_n(\x),\wt\psi_{n-1}(\x),\ldots,\wt\psi_0(\x))
\end{equation}
for any $\x\in (\Bbb C)^*$.
We set 
\begin{align*}
&Y_\x:=\big(\wt\vp_s(\x_{F^{n-r}_r})\big)_{r,s=0}^n,\quad
Y'_\x:=\big(\wt\psi_{n-s}(\x_{F^{n-r}_r})\big)_{r,s=0}^n,\\
&Y_\y:=\big(\wt\vp_s(\y_{F^{n-r}_r})\big)_{r,s=0}^n,\quad
Y'_\y:=\big(\wt\psi_{n-s}(\y_{F^{n-r}_r})\big)_{r,s=0}^n.
\end{align*}
By the relation (9.1), we have 
\begin{align}
Y_\x\ R_{2,1}(\textstyle{x_2\over x_1}\displaystyle)&=Y'_\x,
\label{eq:9.2}
\\
Y_\y\ R_{2,1}(\textstyle{x_2\over x_1}\displaystyle)&=Y'_\y.
\label{eq:9.3}
\end{align} 
We set $(y)^{\ta}:=y_1^{\ta}\cdots y_n^{\ta}$ if $y=(y_1,\ldots, y_n)$. 
The matrices 
$Y_\x$, $Y'_\x$, $Y_\y$  and $Y'_\y$ 
have the following asymptotic behaviors:
\begin{align}
Y_\x\sim (q^{\ta})^{D_{\!A}^+}C_{\!A}^+,\quad 
Y'_\x\sim (q^{\ta})^{D_{\!A}^+}{C'}_{\!A}^+ \quad\hbox{at}\quad\ta\to +\8,
\label{eq:9.4}\\
Y_\y\sim (q^{\ta})^{D_{\!A}^-}C_{\!A}^-,\quad
Y'_\y\sim (q^{\ta})^{D_{\!A}^-}{C'}_{\!A}^- \quad\hbox{at}\quad\ta\to -\8,
\label{eq:9.5}
\end{align}
where 
$$
(q^{\ta})^{D_{\!A}^+}
=\diag[(\x_{F^{n-r}_r} )^{\ta}]_{r=0}^n,\quad
(q^{\ta})^{D_{\!A}^-}
=\diag[(\y_{F^{n-r}_r} )^{\ta}]_{r=0}^n,
$$
and
$C_{\!A}^+=(c_{\!A,rs}^+)_{r,s=0}^n$, 
$C_{\!A}^-=(c_{\!A,rs}^-)_{r,s=0}^n$, 
${C'}_{\!A}^+=({c'}_{\!A,rs}^+)_{r,s=0}^n$ 
and ${C'}_{\!A}^-=({c'}_{\!A,rs}^-)_{r,s=0}^n$ 
are matrices not depending on $q^{\ta}$ defined by 
$$
c_{\!A,rs}^+
:=(1-q)^n \frac{\vP(\x_{F^{n-r}_r} )\vp_s(\x_{F^{n-r}_r} )}
              {(\x_{F^{n-r}_r} )^{\ta}},\quad
c_{\!A,rs}^-
:=(1-q)^n \!\!\!
\Res_{t=\y_{F^{n-r}_r}}\!\!\!
         \frac{\vP(t)\vp_s(t)}
              {(t )^{\ta}}
\frac{d t_1}{t_1}\wedge\cdots\wedge\frac{d t_n}{t_n},
$$
$$
{c'}_{\!A,rs}^{+}
:=(1-q)^n \frac{\vP(\x_{F^{n-r}_r} )\psi_{n-s}(\x_{F^{n-r}_r} )}
               {(\x_{F^{n-r}_r} )^{\ta}},\ \ 
{c'}_{\!A,rs}^-
:=(1-q)^n \!\!\!
\Res_{t=\y_{F^{n-r}_r}}\!\!\!
         \frac{\vP(t)\psi_{n-s}(t)}
              {(t )^{\ta}}
\frac{d t_1}{t_1}\wedge\cdots\wedge\frac{d t_n}{t_n}.
$$
From \eqref{eq:9.2} and \eqref{eq:9.3}, it follows 
\begin{align*}
C_{\!A}^+\ R_{2,1}(\textstyle{x_2\over x_1}\displaystyle)&={C'}_{\!A}^+,
\\
C_{\!A}^-\ R_{2,1}(\textstyle{x_2\over x_1}\displaystyle)&={C'}_{\!A}^-,
\end{align*}
so that we have 
\begin{prop}
\label{prop:9.1}
The matrix $R_{2,1}(\textstyle{x_2\over x_1}\displaystyle)$ is expressed as 
\begin{align}
R_{2,1}(\textstyle{x_2\over x_1}\displaystyle)&=(C_{\!A}^+)^{-1}_0({C'}_{\!A}^+)_0,
\label{eq:9.6}
\\
R_{2,1}(\textstyle{x_2\over x_1}\displaystyle)&=(C_{\!A}^-)^{-1}_0({C'}_{\!A}^-)_0.
\label{eq:9.7}
\end{align}
\end{prop}
\begin{lem}
\label{lem:9.2}
The relations between $Y_\x$ and $Y'_\x$ or between $Y_\y$ and $Y'_\y$ 
are expressed as 
\begin{align}
\wt\psi_s(\x_{F^{r}_{n-r}})
&=\t \left[\wt\vp_s(\x_{F^{n-r}_{r}})\cdot\sgn {\s_r}\cdot U_{\s_r}(\x_{F^{n-r}_r})\right],
\label{eq:9.8}\\
\wt\psi_s(\y_{F^{r}_{n-r}})
&=\t \left[\wt\vp_s(\y_{F^{n-r}_{r}})\cdot\sgn {\s_r}\cdot U_{\s_r}(\y_{F^{n-r}_r})\right],
\label{eq:9.9}
\end{align}
where $\s_r \in \frak S_n$ is 
$$
\s_r=
\left(
\begin{matrix}
        1&2&\cdots&r&r+1&r+2&\cdots&n\\
        n-r+1&n-r+2&\cdots&n&1&2&\cdots&n-r
\end{matrix}
\right).
$$
\end{lem}
\par\noindent{\bf Proof.}
By definition, we have 
\begin{equation}
\label{eq:9.10}
\wt\psi_s(\x_{F^{r}_{n-r}})=\t \left[\t\wt\psi_s(\t\x_{F^{r}_{n-r}})\right],
\end{equation}
where 
$$
\t\x_{F^{r}_{n-r}}=(x_2,x_2q^\c,\ldots,x_2q^{(r-1)\c},x_1,x_1q^\c,\ldots,x_1q^{(n-r-1)\c})
=\s_r^{-1} \x_{F^{n-r}_{r}}.
$$
The RHS of \eqref{eq:9.10} without $\t$ is 
\begin{align*}
\t\wt\psi_s(\t\x_{F^{r}_{n-r}})
&=\wt\vp_s(\s_r^{-1} \x_{F^{n-r}_{r}})
  \quad\quad\hbox{(by \eqref{eq:1.3.1})}\\
&=\int_{\<\s_r^{-1} \x_{F^{n-r}_{r}}\>}\vP(t)\vp_s(t)\varpi
 =\int_{\<\x_{F^{n-r}_{r}}\>}\vP(\s_r^{-1}t)\vp_s(\s_r^{-1}t)\varpi\\
&=\int_{\<\x_{F^{n-r}_{r}}\>}\s_r \vP(t)\cdot \s_r\vp_s(t)\varpi\\
&=U_{\s_r}(\x_{F^{n-r}_r})\cdot \sgn {\s_r}\cdot 
 \int_{\<\x_{F^{n-r}_{r}}\>}\vP(t)\vp_s(t)\varpi\quad\hbox{(by \eqref{eq:1.1.1} and \eqref{eq:1.3.8})}\\
&=U_{\s_r}(\x_{F^{n-r}_r})\cdot \sgn {\s_r}\cdot \wt\vp_s(\x_{F^{n-r}_{r}}).
\end{align*}
Thus we have \eqref{eq:9.8}. 
We also have the relation \eqref{eq:9.9} in the same way as above. 
\hfill\qed
\begin{lem}
\label{lem:9.3}
{\it 
The relations between $C_{\!A}^+$ and ${C'}_{\!A}^+$ or between $C_{\!A}^-$ and ${C'}_{\!A}^-$ 
are expressed as 
\begin{align*}
{c'}_{\!A,n-r,n-s}^{+}
&=\t \left[ c_{\!A,rs}^{+}\cdot\sgn {\s_r}\cdot U_{\s_r}(\x_{F^{n-r}_r})\right],
\\
{c'}_{\!A,n-r,n-s}^{-}
&=\t \left[ c_{\!A,rs}^{-}\cdot\sgn {\s_r}\cdot U_{\s_r}(\y_{F^{n-r}_r})\right].
\end{align*}
}
\end{lem}
\par\noindent{\bf Proof.} It is easily deduced from \eqref{eq:9.4}, \eqref{eq:9.5} and Lemma \ref{lem:9.2}. \hfill\qed
\subsection{{Proof of Theorems \ref{thm:1.3.3} and \ref{thm:1.3.4}}}
From Proposition \ref{prop:3.1}, we have that the matrix $C_{\!A}^-$ is upper triangular.
By Lemma \ref{lem:9.3} the matrix ${C'}_{\!A}^-$ is lower triangular. 
Therefore, from the expression \eqref{eq:9.7} and Lemma \ref{lem:9.3}, 
we have the Gauss decomposition of the matrix $R_{2,1}(\textstyle{x_2\over x_1}\displaystyle)$ as follows: 
$$R_{2,1}(\textstyle{x_2\over x_1}\displaystyle)=U_R\cdot D_R\cdot L_R ,$$
where $U_R^{-1}=(u_{\!R,rs}^*)_{r,s=0}^n$, $D_R=\diag[d_{\!R,0},\ldots,d_{\!R,n}]$
and $L_R=(l_{\!R,rs})_{r,s=0}^n$ are the matrices such that 
\begin{equation}
\label{eq:9.11}
u_{\!R,rs}^*=
            \Big(
                  \frac{c_{\!A,rs}^-}{c_{\!A,rr}^-}
            \Big)_0, 
\quad l_{\!R,rs}=
            \Big(
                  \frac{{c'}_{\!A,rs}^-}{{c'}_{\!A,rr}^-}
            \Big)_0
=\t
            \Big(
                  \frac{c_{\!A,n-r,n-s}^-}{c_{\!A,n-r,n-r}^-}
            \Big)_0,
\end{equation}
$$
d_{\!R,r}=(c_{\!A,rr}^-)_0^{-1}({c'}_{\!A,rr}^-)_0
=(c_{\!A,rr}^-)_0^{-1}
  \t\left[( c_{\!A,n-r,n-r}^-)_0\cdot \sgn {\s_{n-r}}\cdot \big(U_{\s_{n-r}}(\y_{F^{r}_{n-r}})\big)_0
  \right].
$$
Hence Theorem \ref{thm:1.3.3} follows from above Gauss decomposition, Proposition \ref{prop:3.2} and 
\begin{equation}
\label{eq:9.12}
\big(U_{\s_r}(\y_{F^{n-r}_r})\big)_0
  =(-q^{-\c})^{r(n-r)}{}_n(\X21 q^{\b_1-\b_2-(r-1)\c};q^\c)_{r}\cdot {}_n(\X21 q^{\b_1-\b_2-r\c};q^\c)_{r},
\end{equation}
\begin{equation}
\label{eq:9.13}
\sgn {\s_r}=(-1)^{r(n-r)}.
\end{equation}
By using Proposition \ref{prop:3.1} and 
\begin{equation}
\label{eq:9.14}
\big(U_{\s_r}(\x_{F^{n-r}_r})\big)_0
=(-q^{-\c})^{r(n-r)}{}_n(\X21 q^{-(n-r-1)\c};q^\c)_{r}\cdot {}_n(\X21 q^{-(n-r)\c};q^\c)_{r},
\end{equation}
we can prove Theorem \ref{thm:1.3.4} in the same way as Theorem \ref{thm:1.3.3}. \hfill\qed\\
\begin{rmk}
\label{rmk:9.4}
{\rm
From \eqref{eq:9.12}, \eqref{eq:9.13}, \eqref{eq:9.14} and Lemma \ref{lem:9.3}, the matrices $({C'}_{\!A}^+)_0$ and $({C'}_{\!A}^-)_0$ 
are written by 
$$
({C'}_{\!A}^+)_0=J\t\diag[g_0^+,g_1^+,\ldots,g_n^+]\t(C_{\!A}^+)_0 J
=J\t\diag[g_0^+,g_1^+,\ldots,g_n^+]J J\t(C_{\!A}^+)_0 J
$$
and
$$
({C'}_{\!A}^-)_0
=J\t\diag[g_0^-,g_1^-,\ldots,g_n^-]\t(C_{\!A}^-)_0 J 
=J\t\diag[g_0^-,g_1^-,\ldots,g_n^-]J J\t(C_{\!A}^-)_0 J 
$$
where
\begin{align*}
g_r^+
&=q^{-r(n-r)\c}{}_n(\X21 q^{-(n-r-1)\c};q^\c)_{r}\cdot {}_n(\X21 q^{-(n-r)\c};q^\c)_{r}
,\\
g_r^-
&=q^{-r(n-r)\c}{}_n(\X21 q^{\b_1-\b_2-(r-1)\c};q^\c)_{r}\cdot {}_n(\X21 q^{\b_1-\b_2-r\c};q^\c)_{r}
.
\end{align*} 
Thus, from Proposition \ref{prop:9.1}, the $R$-matrix $R_{2,1}(\textstyle{x_2\over x_1}\displaystyle)$ are given by 
\begin{align}
R_{2,1}(\textstyle{x_2\over x_1}\displaystyle)&=(C_{\!A}^+)^{-1}_0\ \diag[\t g_n^+,\t g_{n-1}^+,\ldots,\t g_0^+]\ J \t(C_{\!A}^+)_0 J,
\label{eq:9.15}\\
R_{2,1}(\textstyle{x_2\over x_1}\displaystyle)&=(C_{\!A}^-)^{-1}_0\ \diag[\t g_n^-,\t g_{n-1}^-,\ldots,\t g_0^-]\ J \t(C_{\!A}^-)_0 J. 
\label{eq:9.16}
\end{align}
}\end{rmk}

\section{Examples.}
\label{section:10}
\subsection{The case when $m=2$, $n=1$.}
we explain the simplest case when $n=1$. 
In Theorem \ref{thm:1.3.3}, the $R$-matrix $R_{2,1}(\textstyle{x_2\over x_{1}}\displaystyle)$ is written by 
\begin{align*}
R_{2,1}({\textstyle\frac{x_2}{x_1}}
)
&=
\left(
\begin{array}{cc}
{\X21
-1\over\X21
-q^{\b_2}}
&
{(1-q^{\b_2})\X21\over\X21
-q^{\b_2}}
\\[10pt]
{1-q^{\b_2}\over\X21
-q^{\b_2}}
&
{\X21
q^{\b_1}-q^{\b_2}\over\X21
-q^{\b_2}}
\end{array}
\right)
=U_R\ D_R\ L_R=L'_R\ D'_R\ U'_R
\end{align*}
where 
\begin{align*}
U_R&=
\left(
\begin{array}{cc}
1&-{q^{-\b_2}(1-q^{\b_1})\X21\over 1-\X21
q^{\b_1-\b_2}}
\\[10pt]
0&1
\end{array}
\right),\quad
D_R=
\left(
\begin{array}{cc}
{q^{-\b_2}(1-\X21
q^{\b_1})\over 1-\X21
q^{\b_1-\b_2}}
&0
\\
0&{1-\X21
q^{\b_1-\b_2}\over 1-\X21
q^{-\b_2}}
\end{array}
\right), \\
L_R&=
\left(
\begin{array}{cc}
1&0\\[5pt]
\!\!-{(1-q^{\b_2})q^{-\b_2}\over 1-\X21
q^{\b_1-\b_2}}
&1
\end{array}\right),
\end{align*} 
and 
$$
L'_R=
\left(
\begin{array}{cc}
1&0\\
\!\!-{1-q^{\b_2}\over 1-\X21}
&1
\end{array}\right),\quad 
D'_R=
\left(
\begin{array}{cc}
{q^{-\b_2}(1-\X21) \over 1-q^{-\b_2}\X21 }
&0\\
0&{1-q^{\b_1}\X21 \over 1-\X21}
\end{array}\right)
,\quad 
U'_R=
\left(
\begin{array}{cc}
1&{1-q^{\b_1}\over 1-\textstyle{x_1\over x_2}}\\
0&1
\end{array}\right).
$$
In Theorem \ref{thm:5.0.1}, the Gauss decomposition of the matrix $A(q^{\ta})$ is given by 
\begin{align*}
A(q^{\ta})
&=
\left(
\begin{array}{cc}
1&-{q^{\ta}(1-q^{\b_1})\over 1-q^{\ta+\b_1}}
\\[8pt]
0&1
\end{array}\right)
\left(
\begin{array}{cc}
{(1-q^{\ta})x_1\over 1-q^{\ta+\b_1}}
&0
\\
0&
{(1-q^{\ta+\b_1})x_2\over 1-q^{\ta+\b_1+\b_2}}
\end{array}\right)
\left(
\begin{array}{cc}
1&0\\[5pt]
-{(1-q^{\b_2})x_1\over (1-q^{\ta+\b_1})x_2}
&1
\end{array}\right)\\
&=
\left(
\begin{array}{cc}
{(1-q^{\ta+\b_2})x_1\over 1-q^{\ta+\b_1+\b_2}}
&
-{q^{\ta}(1-q^{\b_1})x_2\over 1-q^{\ta+\b_1+\b_2}}
\\[5pt]
-{(1-q^{\b_2})x_1\over 1-q^{\ta+\b_1+\b_2}}
&
{(1-q^{\ta+\b_1})x_2\over 1-q^{\ta+\b_1+\b_2}}
\end{array}\right),
\end{align*}
so that we have 
$$
A(0)
=
\left(
\begin{array}{cc}
x_1
&
0
\\
-(1-q^{\b_2})x_1
&
x_2
\end{array}\right)
,\quad
A(\8)
=
\left(
\begin{array}{cc}
x_1 q^{-\b_1}
&
-(1-q^{-\b_1}) x_2 q^{-\b_2}
\\
0
&
x_2 q^{-\b_2}
\end{array}\right).
$$
From \eqref{eq:3.6} and \eqref{eq:3.7}, taking the unipotent matrices $C^+$ and $C^-$ such that 
$$
C^+
=
\left(
\begin{array}{cc}
1
&
0
\\[0.5pt]
{1-q^{\b_2}
\over 
1-\X21}
&
1
\end{array}\right),\quad
C^-
=
\left(
\begin{array}{cc}
1
&
{q^{-\b_2}(1-q^{\b_1})\X21
\over 1-q^{\b_1-\b_2}\X21
}
\\
0
&
1
\end{array}\right),
$$
we can diagonalize $A(0)$ and $A(\8)$ as follows:
$$
A(0)
=
(C^+)^{-1}
\left(
\begin{array}{cc}
x_1
&
0
\\
0&x_2
\end{array}\right)
C^+,
\quad
A(\8)
=
(C^-)^{-1}
\left(
\begin{array}{cc}
x_1q^{-\b_1}
&
0
\\
0&x_2q^{-\b_2}
\end{array}\right)
C^-
$$
From Propositions \ref{prop:3.1} and \ref{prop:3.2}, the matrices $(C_{\!A}^+)_0$ and $(C_{\!A}^-)_0$ are written 
as a product of diagonal and unipotent matrices as follows:
$$
(C_{\!A}^+)_0
=\diag[(c_{\!A,00}^+)_0,(c_{\!A,11}^+)_0]\ C^+,\quad 
(C_{\!A}^-)_0
=\diag[(c_{\!A,00}^-)_0,(c_{\!A,11}^-)_0]\ C^-,
$$
where 
$$
\diag[(c_{\!A,00}^+)_0,(c_{\!A,11}^+)_0]
=
\left(
\begin{array}{cc}
{1\over 1-q^{\b_1}}
&
0
\\
0
&
{1-\X21\over (1-q^{\b_1}\X21
)(1-q^{\b_2})}
\end{array}\right)
$$
and
$$
\diag[(c_{\!A,00}^-)_0,(c_{\!A,11}^-)_0]
=
\left(
\begin{array}{cc}
-1
&
0
\\
0
&
-{1-q^{-\b_2}\X21
\over 1-q^{\b_1-\b_2}\X21
}
\end{array}\right).
$$
From \eqref{eq:9.15} and \eqref{eq:9.16},
the $R$-matrix $R_{2,1}(\textstyle{x_2\over x_{1}}\displaystyle)$ is determined 
from the matrices $(C_{\!A}^+)_0$ or $(C_{\!A}^-)_0$ as follows:
$$
R_{2,1}(\textstyle{x_2\over x_{1}}\displaystyle)
=(C_{\!A}^+)_0^{-1} 
J \t(C_{\!A}^+)_0 J
=(C_{\!A}^-)_0^{-1} 
J \t(C_{\!A}^-)_0 J,
$$
and it follows that 
$$
D'_R
=
\diag[(c_{\!A,00}^+)_0^{-1}\t(c_{\!A,11}^+)_0,
     \ (c_{\!A,11}^+)_0^{-1}\t(c_{\!A,00}^+)_0]
$$
and 
$$
D_R
=
\diag[(c_{\!A,00}^-)_0^{-1}\t(c_{\!A,11}^-)_0,
     \ (c_{\!A,11}^-)_0^{-1}\t(c_{\!A,00}^-)_0].
$$
Thus, in particular, the upper and lower triangular matrices 
as factors of the Gauss matrix decomposition of $R_{2,1}(\textstyle{x_2\over x_{1}}\displaystyle)$
are determined from the matrices $C^+$ and $C^-$ as follows:
$$L'_R=(C^+)^{-1},\quad U'_R=J\t C^+ J,$$
and 
$$U_R=(C^-)^{-1},\quad L_R=J\t C^- J.$$

%
\subsection{The case when $m=2$, $n=2$.}
\par
In Theorem \ref{thm:5.0.1}, the Gauss decomposition of the matrix $A(q^{\ta})$ is given by 
\begin{align*}
A(q^{\ta})
&=
\left(
\begin{array}{ccc}
1
&-{q^{\ta-\c}(1-q^{\b_1+\c})\over 1-q^{\ta+\b_1}}
&{q^{2\ta-3\c}(1-q^{\b_1})(1-q^{\b_1+\c})\over (1-q^{\ta+\b_1-2\c})(1-q^{\ta+\b_1-\c})}
\\[5pt]
0&1
&-{q^{\ta-2\c}(1-q^{2\c})(1-q^{\b_1})\over (1-q^{\c})(1-q^{\ta+\b_1-2\c})}
\\[5pt]
0&0&1
\end{array}\right)\\
&\quad\cdot 
\left(
\begin{array}{ccc}
{x_1^2q^\c(1-q^{\ta-2\c})(1-q^{\ta-\c})\over (1-q^{\ta+\b_1-\c})(1-q^{\ta+\b_1})}&0&0\\
0&{x_1x_2(1-q^{\ta+\b_1})(1-q^{\ta-2\c})\over (1-q^{\ta+\b_1+\b_2})(1-q^{\ta+\b_1-2\c})}&0\\
0&0&{x_2^2q^\c(1-q^{\ta+\b_1-2\c})(1-q^{\ta+\b_1-\c})\over (1-q^{\ta+\b_1+\b_2-\c})(1-q^{\ta+\b_1+\b_2})}
\end{array}\right)\\
&\quad\cdot 
\left(
\begin{array}{ccc}
1&0&0\\[1pt]
-{(1-q^{2\c})(1-q^{\b_2})x_1\over (1-q^\c)(1-q^{\ta+\b_1})x_2}&1&0\\[6pt]
{q^{-\c}(1-q^{\b_2})(1-q^{\b_2+\c})x_1^2\over (1-q^{\ta+\b_1-2\c})(1-q^{\ta+\b_1-\c})x_2^2}
&-{q^{-\c}(1-q^{\b_2+\c})x_1\over (1-q^{\ta+\b_1-2\c})x_2}&1
\end{array}\right),
\end{align*}
so that we have 
\begin{align*}
A(0)
&=
\left(
\begin{array}{ccc}
x_1^2q^\c&0&0\\
0&x_1x_2&0\\
0&0&x_2^2q^\c
\end{array}\right)
\left(
\begin{array}{ccc}
1&0&0\\
-{(1-q^{2\c})(1-q^{\b_2})x_1\over (1-q^\c)x_2}
&1&0\\[5pt]
{q^{-\c}(1-q^{\b_2})(1-q^{\b_2+\c})x_1^2\over x_2^2}
&-{q^{-\c}(1-q^{\b_2+\c})x_1\over x_2}&1
\end{array}\right)\\
&=\left(
\begin{array}{ccc}
q^\c x_1^2&0&0\\
-{(1+q^{\c})(1-q^{\b_2})x_1^2}&x_1x_2&0\\
{(1-q^{\b_2})(1-q^{\b_2+\c})x_1^2}
&-{(1-q^{\b_2+\c})x_1x_2}&q^\c x_2^2
\end{array}\right)
\end{align*}
and
\begin{align*}
A(\8)
&=
\left(\begin{array}{ccc}
1
&{q^{-\b_1-\c}(1-q^{\b_1+\c})}
&{q^{-2\b_1}(1-q^{\b_1})(1-q^{\b_1+\c})}
\\
0&1
&{q^{-\b_1}(1+q^{\c})(1-q^{\b_1})}
\\
0&0&1
\end{array}\right)\\
&\qquad \cdot 
\left(\begin{array}{ccc}
x_1^2q^{-2\b_1-\c}&0&0\\
0&x_1x_2 q^{-\b_1-\b_2}&0\\
0&0&x_2^2q^{-2\b_2-\c}
\end{array}\right)\\
&=
\left(\begin{array}{ccc}
x_1^2q^{-2\b_1-\c}
&
-x_1x_2 q^{-\b_1-\b_2}(1-q^{-\b_1-\c})
&
x_2^2q^{-2\b_2}(1-q^{-\b_1})(1-q^{-\b_1-\c})
\\
0 & x_1x_2 q^{-\b_1-\b_2} & -x_2^2 q^{-2\b_2}(1+q^{-\c})(1-q^{-\b_1})
\\
0&0&x_2^2q^{-2\b_2-\c}
\end{array}\right).
\end{align*}
From \eqref{eq:3.6} and \eqref{eq:3.7}, taking the unipotent matrices $C^+$ and $C^-$ such that 
\begin{align*}
C^+
&=
\left(\begin{array}{ccc}
1&0&0\\
{q^{-\c}(1-q^{2\c})(1-q^{\b_2})\over (1-q^\c)(1-\X21
q^{-\c})}&1&0\\[8pt]
{(1-q^{\b_2})(1-q^{\b_2+\c})\over (1-\X21
 )(1-\X21
q^{\c})}
&{(1-q^{\b_2+\c})\over (1-\X21
q^{\c})}&1
\end{array}\right),\\[5pt]
C^-
&=
\left(\begin{array}{ccc}
1
&{\X21
q^{-\b_2}(1-q^{\b_1+\c})\over 1-\X21
q^{\b_1-\b_2+\c}}
&{(\X21
q^{-\b_2})^2(1-q^{\b_1})(1-q^{\b_1+\c})
  \over (1-\X21
q^{\b_1-\b_2})(1-\X21
q^{\b_1-\b_2+\c})}
\\[10pt]
0&1
&{\X21
q^{-\b_2-\c}(1-q^{2\c})(1-q^{\b_1})\over (1-q^{\c})(1-\X21
q^{\b_1-\b_2-\c})}
\\[8pt]
0&0&1
\end{array}\right),
\end{align*}
we can diagonalize $A(0)$ and $A(\8)$ as follows:
\begin{align*}
A(0)
&=
(C^+)^{-1}
\left(\begin{array}{ccc}
x_1^2q^{\c}&0&0\\
0&x_1 x_2&0\\
0&0&x_2^2q^{\c}
\end{array}\right)
C^+,
\\
A(\8)
&=
(C^-)^{-1}
\left(\begin{array}{ccc}
x_1^2q^{-2\b_1-\c}&0&0\\
0&x_1x_2 q^{-\b_1-\b_2}&0\\
0&0&x_2^2q^{-2\b_2-\c}
\end{array}\right)
C^-.
\end{align*}
\par
From Propositions \ref{prop:3.1} and \ref{prop:3.2}, the matrices $(C_{\!A}^+)_0$ and $(C_{\!A}^-)_0$ are written 
as a product of diagonal and unipotent matrices as follows:
\begin{align*}
(C_{\!A}^+)_0
&=\diag[(c_{\!A,00}^+)_0,(c_{\!A,11}^+)_0,(c_{\!A,22}^+)_0]\ C^+,\\
(C_{\!A}^-)_0
&=\diag[(c_{\!A,00}^-)_0,(c_{\!A,11}^-)_0,(c_{\!A,22}^-)_0]\ C^-,
\end{align*}
where 
\begin{align*}
&\diag[(c_{\!A,00}^+)_0,(c_{\!A,11}^+)_0,(c_{\!A,22}^+)_0]\\
&\!=
\hskip -3pt
\left(\!\begin{array}{ccc}
{q^\c \over (1-q^{\b_1})(1-q^{\b_1+\c})}\!\!\!\!&0&0\\
0&\!\!\!\!{(1-\X21
q^{-\c})(1-\X21
) 
    \over (1-q^{\b_1})(1-\X21
    q^{\b_1})(1-q^{\b_2})(1-\X21
    q^{\c})}\!\!\!\!&0\\
0&0&\!\!\!\!{q^\c(1-\X21
)(1-\X21
q^{\c})
    \over (1-\X21
    q^{\b_1})(1-\X21
    q^{\b_1+\c})(1-q^{\b_2})(1-q^{\b_2+\c})}
\end{array}\!\right)
\end{align*}
\begin{align*}
&\diag[(c_{\!A,00}^-)_0,(c_{\!A,11}^-)_0,(c_{\!A,22}^-)_0]\\
&=
\left(\begin{array}{ccc}
{1-q^{-2\c} \over 1-q^{-\c}}&0&0\\
0&{(1-\X21
q^{-\b_2})(1-\X21
q^{\b_1-\b_2-\c})
    \over (1-\X21
    q^{\b_1-\b_2})(1-\X21
    q^{\b_1-\b_2+\c})}&0\\
0&0&{(1-q^{-2\c})(1-\X21
q^{-\b_2-\c})(1-\X21
q^{-\b_2})
    \over (1-q^{-\c})(1-\X21
    q^{\b_1-\b_2-\c})(1-\X21
    q^{\b_1-\b_2})}
\end{array}\right)
\end{align*}
\par
From \eqref{eq:9.15} and \eqref{eq:9.16}, 
the $R$-matrix $R_{2,1}(\textstyle{x_2\over x_{1}}\displaystyle)$ is determined 
from the matrices $(C_{\!A}^+)_0$ or $(C_{\!A}^-)_0$ as follows:
\begin{align*}
R_{2,1}(\textstyle{x_2\over x_{1}}\displaystyle)
&=(C_{\!A}^+)_0^{-1} 
\left(\begin{array}{ccc}
1&0&0\\
0&\t g_1^+
&0
\\
0&0&1
\end{array}\right)
J \t(C_{\!A}^+)_0 J\\
&=(C_{\!A}^-)_0^{-1} 
\left(\begin{array}{ccc}
1&0&0\\
0&\t g_1^-
&0\\
0&0&1
\end{array}\right)
J \t(C_{\!A}^-)_0 J,
\end{align*}
where 
$$
\t g_1^+
={q^{-\c} (1-\textstyle{x_1\over x_2}\displaystyle q^\c) \over 1-\textstyle{x_1\over x_2}\displaystyle q^{-\c}},\quad
\t g_1^-
={q^{-\c} (1-\textstyle{x_1\over x_2}\displaystyle q^{\b_2-\b_1+\c}) 
    \over 1-\textstyle{x_1\over x_2}\displaystyle q^{\b_2-\b_1-\c}}.$$
If we express the Gauss decomposition of $R_{2,1}(\textstyle{x_2\over x_{1}}\displaystyle)$ by
$$R_{2,1}(\textstyle{x_2\over x_{1}}\displaystyle)=U_R\ D_R\ L_R=L'_R\ D'_R\ U'_R,$$
then we have 
$$
D'_R
=
\diag[(c_{\!A,00}^+)_0^{-1} \t(c_{\!A,22}^+)_0,
     \ (c_{\!A,11}^+)_0^{-1} \t g_1^+ \t(c_{\!A,11}^+)_0,
     \ (c_{\!A,22}^+)_0^{-1} \t(c_{\!A,00}^+)_0]
$$
and 
$$
D_R
=
\diag[(c_{\!A,00}^-)_0^{-1} \t(c_{\!A,22}^-)_0,
     \ (c_{\!A,11}^-)_0^{-1} \t g_1^- \t(c_{\!A,11}^-)_0,
     \ (c_{\!A,22}^-)_0^{-1} \t(c_{\!A,00}^-)_0].
$$
In particular, the upper and lower triangular matrices 
as factors of the Gauss matrix decomposition of $R_{2,1}(\textstyle{x_2\over x_{1}}\displaystyle)$
are determined from the matrices $C^+$ and $C^-$ as follows:
$$L'_R=(C^+)^{-1},\quad U'_R=J\t C^+ J,$$
and 
$$U_R=(C^-)^{-1},\quad L_R=J\t C^- J.$$
Therefore the Gauss decomposition of $R_{2,1}(\textstyle{x_2\over x_{1}}\displaystyle)$ is given by
\begin{align*}
&R_{2,1}(\textstyle{x_2\over x_{1}}\displaystyle)\\
&=
\left(\begin{array}{ccc}
1
&-{\X21
q^{-\b_2}(1-q^{\b_1+\c})\over 1-\X21
q^{\b_1-\b_2+\c}}
&{(\X21
q^{-\b_2})^2 q^{-\c}(1-q^{\b_1})(1-q^{\b_1+\c})
     \over (1-\X21
     q^{\b_1-\b_2-\c})(1-\X21
     q^{\b_1-\b_2})}
\\[10pt]
0&1
&-{\X21
q^{-\b_2-\c}(1-q^{2\c})(1-q^{\b_1})\over (1-q^{\c})(1-\X21
q^{\b_1-\b_2-\c})}
\\[8pt]
0&0&1
\end{array}\right)\\
&
\cdot \!
\left(\!\!\!\begin{array}{ccc}
{q^{-2\b_2}(1-\X21
q^{\b_1})(1-\X21
q^{\b_1+\c})
   \over (1-\X21
   q^{\b_1-\b_2})(1-\X21
   q^{\b_1-\b_2+\c})}\!\!\!&0&0\\
0&\!\!\!{q^{-\b_2-\c}(1-\X21
q^{\b_1})(1-\X21
q^{\b_1-\b_2+\c})
   \over (1-\X21
   q^{-\b_2})(1-\X21
   q^{\b_1-\b_2-\c})}\!\!\!&0\\
0&0&\!\!\!{(1-\X21
q^{\b_1-\b_2-\c})(1-\X21
q^{\b_1-\b_2})
\over (1-\X21
q^{-\b_2-\c})(1-\X21
q^{-\b_2})}
\end{array}\!\!\!\right)\\
&
\cdot \!
\left(\begin{array}{ccc}
1&0&0\\[1pt]
-{q^{-\b_2}(1-q^{2\c})(1-q^{\b_2})\over (1-q^\c)(1-\X21
q^{\b_1-\b_2+\c})}&1&0\\[10pt]
{q^{-2\b_2-\c}(1-q^{\b_2})(1-q^{\b_2+\c})
    \over (1-\X21
    q^{\b_1-\b_2-\c})(1-\X21
    q^{\b_1-\b_2})}
&-{q^{-\b_2-\c}(1-q^{\b_2+\c})x_1\over 1-\X21
q^{\b_1-\b_2-\c}}&1
\end{array}\right).
\end{align*}
%

\end{document}